\documentclass[journal]{IEEEtran}

\ifCLASSINFOpdf
   \usepackage[pdftex]{graphicx}
\else
   \usepackage[dvips]{graphicx}
\fi
\usepackage{amsmath}
\usepackage{amscd,amssymb}
\usepackage{amsfonts}
\usepackage{graphicx,epsf}
\usepackage{enumerate,color}
\usepackage[mathcal]{euscript}
\usepackage[dvipsnames]{xcolor}
\usepackage{hyperref}
\usepackage{subfigure}
\usepackage{amsmath}
\usepackage{amssymb}
\usepackage{amsfonts}
\usepackage{dsfont}

\usepackage{afterpage}

\usepackage{epstopdf}

\usepackage{cite}
\usepackage{enumitem}

\usepackage{tcolorbox}

\usepackage[utf8]{inputenc}
\usepackage{tikz}
\usetikzlibrary{shapes.geometric, arrows}

\tikzstyle{startstop} = [rectangle, rounded corners, minimum width=3cm, minimum height=1cm,text centered, text width=5cm, draw=black, fill=orange!10]
\tikzstyle{io} =  [rectangle, rounded corners, minimum width=3cm, minimum height=1cm,text centered, text width=5cm, draw=black, fill=blue!20]
\tikzstyle{process} = [rectangle, rounded corners, minimum width=3cm, minimum height=1cm, text width=5cm, text centered, draw=black, fill=yellow!30]
\tikzstyle{decision} = [diamond, minimum width=3cm, minimum height=1cm, text centered, text width=3cm,draw=black, fill=green!30]
\tikzstyle{arrow} = [thick,->,>=stealth]

\newcommand{\ND}{\mathcal{R}}
\newcommand{\expikz}{e^{ikz}}

\newcommand{\expiconjkz}{e^{i\bar{k}\bar{z}}}

\DeclareMathOperator{\de}{\partial}

\DeclareMathOperator{\dez}{\de_z}

\DeclareMathOperator{\dbar}{\overline{\partial}}
\DeclareMathOperator{\dbarz}{\dbar_z}
\DeclareMathOperator{\dbark}{\dbar_k}

\DeclareMathOperator{\dnu}{\partial_\nu}

\DeclareMathOperator{\by}{\times}
\DeclareMathOperator{\bndry}{\partial\Omega}
\DeclareMathOperator{\T}{\mathbf{t}}

\DeclareMathOperator{\R}{\mathbb R}
\DeclareMathOperator{\C}{\mathbb C}

\DeclareMathOperator{\muexp}{\mathbf{\mu}^{\mbox{\tiny\textbf{exp}}}}
\DeclareMathOperator{\texp}{\mathbf{t}^{\mbox{\tiny \textbf{exp}}}}

\DeclareMathOperator{\sigexp}{\sigma^{\mbox{\tiny \textbf{exp}}}}

\newcommand{\trev}[1]{{#1}}
\newcommand{\trevNew}[1]{{#1}}

\usepackage{soul,color}
\soulregister\cite7
\soulregister\ref7
\soulregister\pageref7

\begin{document}

\title{Deep D-bar: Real time Electrical Impedance Tomography Imaging with Deep Neural Networks}


\author{ S.~J. Hamilton and  A. Hauptmann

\thanks{Copyright (c) 2017 IEEE. Personal use of this material is permitted. However, permission to use this material for any other purposes must be obtained from the IEEE by sending a request to \url{pubs-permissions@ieee.org.}}  
\thanks{S.~J.~Hamilton is with the Department of Mathematics, Statistics, and Computer Science, Marquette University, Milwaukee, WI, 53233 USA ({e-mail: sarah.hamilton@marquette.edu}).}
\thanks{A. Hauptmann is with the Department of Computer Science; University College London, London, United Kingdom, (email: {a.hauptmann@ucl.ac.uk})}
\thanks{Accompanying codes will be made available: \url{https://github.com/asHauptmann/DeepDbar}}
}

\markboth{}{Hamilton and Hauptmann: Deep D-bar}
%

\maketitle

\begin{abstract}
The mathematical problem for Electrical Impedance Tomography (EIT) is a highly nonlinear ill-posed inverse problem requiring carefully designed reconstruction procedures to ensure reliable image generation. D-bar methods are based on a rigorous mathematical analysis and provide robust direct reconstructions by using a low-pass filtering of the associated nonlinear Fourier data.  Similarly to low-pass filtering of linear Fourier data, only using low frequencies in the image recovery process results in blurred images lacking sharp features such as clear organ boundaries.  Convolutional Neural Networks provide a powerful framework for post-processing such convolved direct reconstructions. In this study, we demonstrate that these CNN techniques lead to sharp and reliable reconstructions even for the highly nonlinear inverse problem of EIT.  The network is trained on data sets of simulated examples and then applied to experimental data without the need to perform an additional transfer training.  Results \trev{for absolute EIT images} are presented \trev{using} experimental EIT data from the ACT4 and KIT4 EIT systems.
\end{abstract}


\IEEEpeerreviewmaketitle

\section{Introduction}\label{sec:intro}
\IEEEPARstart{E}{{lectrical Impedance Tomography}} (EIT) images traditionally display the tissue-dependent conductivity distribution of a patient in the plane of the attached measurement electrodes allowing, e.g.,  visualization of heart and lung function as well as injuries \cite{Cinnella2015, Grant2011, Karagiannidis2015,  Pesenti2016, Reinius2015, Schlibler2013}. The resulting images are of high-contrast and data acquisition is done by harmless electrical measurements without the need for contrast agents or ionizing radiation. However, the image recovery process of forming the EIT image from the current/voltage measurement data is a severely ill-posed nonlinear inverse problem, and thus requires a noise-robust regularization strategy for stability.  The `D-bar method', the only proven regularization strategy for the full nonlinear problem \cite{Knudsen2009}, provides real-time noise-robust image recovery by using a low-pass filter of the associated nonlinear Fourier data.  Unfortunately, this results in images that suffer a loss of sharp features often important in medical imaging applications.  In this work, we propose combining {D-bar} with Deep Learning, specifically with a Convolutional Neural Network, to `learn' and undo the image blurring resulting in real-time sharp EIT images.

EIT reconstructions are typically computed with iterative algorithms that are based on minimizing a penalty functional, such as \cite{Zhou2015,Gonzalez2017}. These methods perform very well in reconstruction quality due to a flexibility of incorporating prior knowledge, but require careful modeling of the boundary shape in the repeated simulation of the forward problem. Possibilities to overcome the boundary sensitivity are proposed in \cite{Darde2013,Nissinen2011}, but tend to be computationally demanding.  On the other hand, direct (non-iterative) reconstruction algorithms do not need the \trev{repeated} simulation of the forward operator.  One such method is known as the D-bar algorithm which is based on a nonlinear Fourier transformation of the measured surface current/voltage data. The method employs a low-pass filtering of this transformed data as a regularization strategy to stabilize the image reconstruction process against noise in the measured data. \trev{Consequently, this filtering results in reconstructed images that suffer from a significant loss of sharpness. It has been shown that the direct D-bar method is robust to incorrect or incomplete knowledge of electrode locations as well as errors in boundary shape, see for instance \cite{Murphy2009} and the discussion in Section \ref{sec:Robustness}. Iterative methods on the other hand are either very sensitive to the correct forward model or are based on sophisticated modelling to cope with uncertainties in the model, such as unknown electrode locations, boundary shape, or contact impedances \cite{Nissinen2011,Kolehmainen2008b,Darde2013}.}

Recent advances in the larger field of image reconstruction have demonstrated the power of Deep Learning and Neural Networks for improving low quality or corrupted images.  In particular, combining fast direct reconstruction procedures with deep neural networks can provide high quality images with low latency, leading to prospective real-time imaging in many applications.  Convolutional Neural Networks (CNN) are especially suitable for post-processing initial reconstructions that come from algorithms based on, or related to, Fourier transforms, as suggested in \cite{Jin2017}.   Such initial reconstructions typically suffer from a loss of spatial resolution, due to some sort of low-pass filtering, as well as additional undersampling artefacts. Training a CNN to remove these artefacts to improve the information content of the reconstructed image has been studied for several linear inverse problems in medical imaging, including CT \cite{Jin2017,Kang2017}, MRI \cite{Sandino2017}, and PAT \cite{Antholzer2017,Hauptmann2018}.  Although the EIT problem is nonlinear in nature, the low-pass filtered images from the low-passed D-bar method naturally fit into this setting.

In this study we formulate a real-time capable reconstruction algorithm that produces high quality sharp \trev{absolute} EIT images by combining the D-bar algorithm with subsequent processing by a CNN. For this task we utilize an established CNN architecture, known as U-net, adjusted to cope with the typical image structures of D-bar EIT reconstructions. We train the network on simulated training data and directly apply the trained network to experimental data with no training on experimental data itself. This successful transition to experimental data underlines the robustness of the D-bar algorithm and is especially important as the need for good training data is often the bottleneck for the success of such network-based approaches for other imaging modalities, \cite{Jin2017, Adler2017, Hauptmann2018}.

This paper is organized as follows.  Section~\ref{sec:EITmath} presents a brief review of the mathematical problem of EIT and the D-bar solution method.  The deep learning CNN for D-bar, coined `Deep D-bar' is introduced in Section~\ref{sec:DeepDbar}.  The experimental setup as well as simulation of training data are described in Section~\ref{sec:ExpSetups} and results presented in Section~\ref{sec:results}.  A discussion of the results is given in Section~\ref{sec:discussion} and conclusions drawn in Section~\ref{sec:conclusion}.  The reader is encouraged to view the manuscript on a computer screen as details in the image contrast may be masked in printed versions.

\section{Electrical Impedance Tomography and the D-bar reconstruction method}\label{sec:EITmath}

%

Electrical impedance tomography is a nonlinear inverse problem in which we aim to determine the interior conductivity from current-to-voltage measurements at the boundary.  
The problem can be formulated as a generalized Laplace equation
\begin{equation}
\label{eqn:condEq}
\left\{
\begin{array}{rcl}
\nabla\cdot\sigma\nabla u &=& 0 \mbox{ in }\Omega,\\
\sigma\dnu u &=& \varphi \mbox{ on } \bndry, \\
\end{array}
\right.
\end{equation}
modeling the electrical potential $u$ inside the domain $\Omega\subset\R^n$ for a given conductivity $\sigma$, with the Neumann boundary condition describing the boundary voltage occurring from the applied mean-free current $\varphi$.  The measurement data consists of pairs of current and voltage measurements and is modeled by the current-to-voltage map $\ND_\sigma$ defined by
\[\ND_\sigma \varphi := \left. u\right|_{\bndry}.\]
This measurement operator is also known as the Neumann-to-Dirichlet (ND) map, and knowledge of it allows one to predict the resulting voltage for any injected current pattern for $n=2,3$.  In practice, an approximation to the ND map is formed by applying a basis of current patterns and tracking the responses of the voltages.   The D-bar algorithm we use below requires the corresponding Dirichlet-to-Neumann (DN) map, which can be obtained as the inverse of the ND map, $\Lambda_{\sigma}=\left(\ND_{\sigma}\right) ^{-1}$, for full (vs. partial) boundary data.  In this work we consider the $n=2$ case as the D-bar reconstruction framework is further developed in 2D.  However, we expect a natural extension to 3D \cite{Delbary2012}.

\subsection{Real-time reconstructions using an approximate D-bar method}\label{sec:Dbar}
By the D-bar method, we refer to the regularized D-bar method \cite{Knudsen2009} based on the theoretical proof given in \cite{Nachman1996}.  The approach uses a nonlinear Fourier transform, called a {\it scattering transform}, tailor-made for the EIT problem which is applied to the measured current/voltage data in the form of the DN map $\Lambda_\sigma$.  That scattering data is then used as input data into a partial differential equation, a $\dbark$ or `D-bar' equation, giving the method its name.  Note that the derivative operators $\dez$ and $\dbarz$ are defined as 
$\dez = \frac12\left(\partial_{z_1} - i\partial_{z_2}\right)$ and $\dbarz = \frac12\left(\partial_{z_1} + i\partial_{z_2}\right)$, where $z=z_1+iz_2\in\C$.  The conductivity $\sigma$ is then recovered directly from the solution to the D-bar equation.  

The D-bar approach \cite{Nachman1996} is to transform the physical conductivity equation $\nabla\cdot\sigma\nabla u=0$ into a nonphysical Schr\"odinger equation,
solve that problem instead using the D-bar methods popularized by Beals and Coifman \cite{Beals1985}, and then transform back to the physical setting.  The change of variables $\tilde{u}=\sigma^{1/2}u$ and $q(z)=\sigma^{-1/2}(z)\Delta \sigma^{1/2}(z)$ produces the desired Schr\"odinger equation $[-\Delta + q(z)]\tilde{u}(z)=0$, where  $z\in\Omega$.
Provided that $\sigma(z)$ is constant in a neighborhood of the boundary, without loss of generality $\sigma=1$ near $\bndry$, the conductivity can be extended from $\Omega$ to the entire plane by setting $\sigma(z)\equiv1$ for $z\in\C\setminus\Omega$.  Note that this gives the potential $q(z)$ compact support in $\Omega$.  We make use of special solutions $\psi(z,k)$ to the Schr\"odinger equation 
\begin{equation}\label{eq:schro}
[-\Delta + q(z)]\psi(z,k)=0, \qquad z\in\C,\;\; k\in\C\setminus\{0\},
\end{equation}
called {\it Complex Geometrical Optics (CGO) solutions}, that have a specific asymptotic behavior for large $|z|$ or $|k|$, {$\psi(z,k)\sim e^{ikz}$}.  Note that we associate $\R^2$ with $\C$ via the mapping $z=(z_1,z_2)\mapsto z_1+iz_2$ and thus $kz=\left(k_1 + ik_2\right)\left(z_1+iz_2\right)$ denotes complex multiplication.
The CGO solutions {$\mu(z,k)=e^{-ikz}\psi(z,k)\sim 1$} solve a D-bar equation in the nonphysical scattering variable $k$
\begin{equation}\label{eq:dbark}
\dbark \mu(z,k) = \frac{1}{4\pi\bar{k}}\T(k)e(z,-k)\overline{\mu(z,k)},
\end{equation}
where $e(z,k):=\exp\{i(kz+\bar{k}\bar{z})\}$ and $\T(k)$ is the nonlinear scattering data defined by 
\begin{equation}\label{eq:t-domain}
\T(k):=\int_{\C} e(z,k)q(z)\mu(z,k)\;dz.
\end{equation}
Note that this scattering data $\T$ can be thought of as nonlinear Fourier data by the following observation.  Replacing the CGO solutions $\mu(z,k)$ in \eqref{eq:t-domain} with the asymptotic behavior $1$ yields
\[\texp(k) = \int_{\C} e(z,k)q(z)(1)\;dz = \hat{q}(-2k_1,2k_2),\]
and thus the `Born' approximation $\texp$ is essentially a shifted Fourier transform of the potential $q$.  A connection to the measurement data $\Lambda_\sigma$ can be established via Alessandrini's identity \cite{Alessandrini1988}
\[\texp(k)=\int_{\C} e(z,k) q(z)  dz= \int_{\bndry} \expiconjkz (\Lambda_\sigma-\Lambda_1) \expikz dz.\]

\begin{figure*}[h!]
\begin{picture}(500,220)
\put(0,0){\includegraphics[width=500pt]{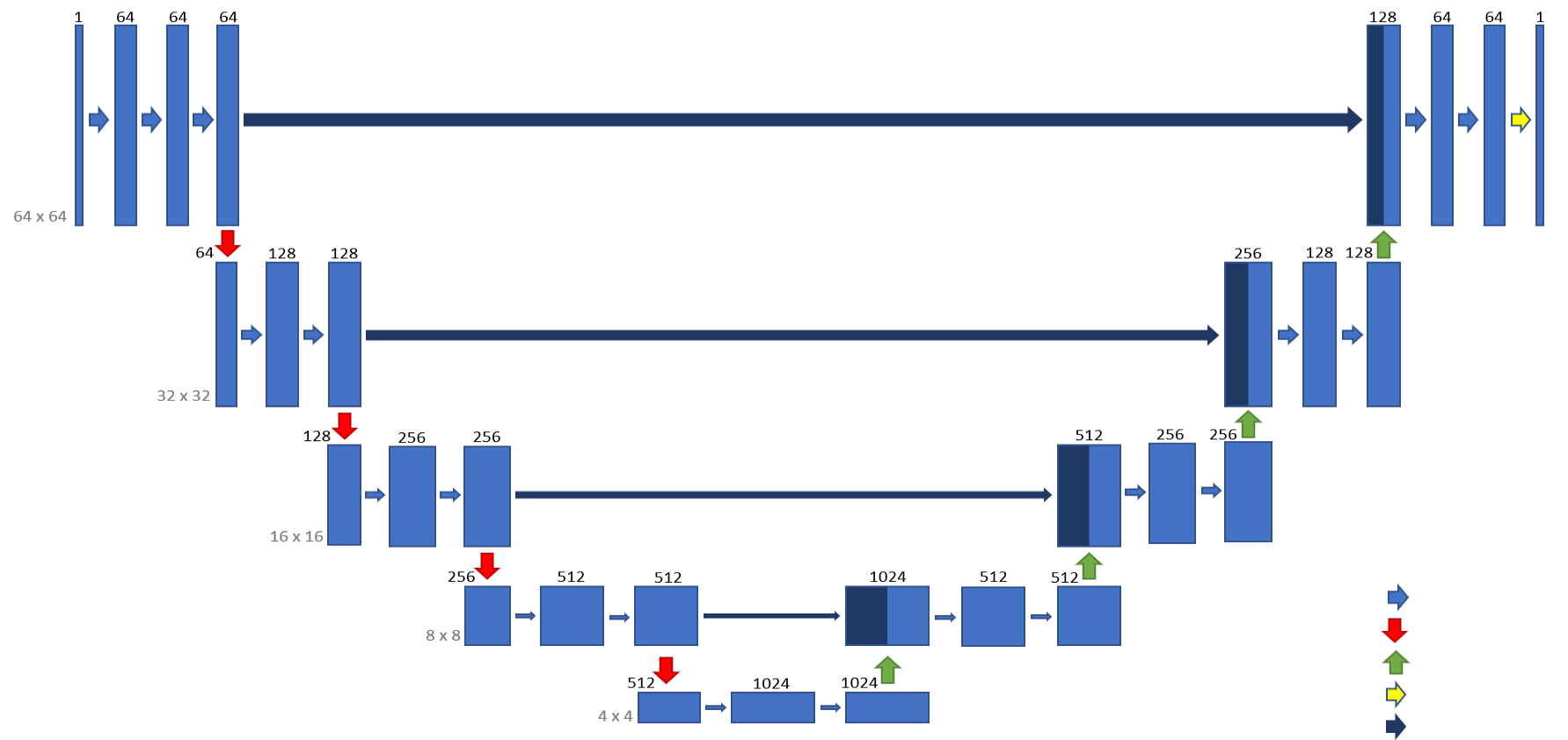}}
\put(452,42.5){\footnotesize ReLU(conv$_{5\times 5}$)}
\put(452,33){\footnotesize maxpool$_{2\times 2}$}
\put(452,21.5){\footnotesize ReLU(convt$_{5\times 5}$)}
\put(452,11.5){\footnotesize ReLU(conv$_{1\times 1}$)}
\put(452,2){\footnotesize  concat}
\put(-3,196){\footnotesize $\sigma^{\exp}=$ }
\put(495,196){\footnotesize $=\widetilde{\sigma}$ }
\end{picture}
\caption{Deep D-bar network structure. The input is given by the D-bar reconstruction $\sigma^{\exp}$ with a resolution of $64\times64$ and the output is denoted by $\widetilde{\sigma}$. The numbers on top of the blue bars denote the channels for each layer. The resolution for each multilevel decomposition is shown in gray on the left. Each convolutional layer is equipped with a Rectified Linear Unit as nonlinearity, given by $\mathrm{ReLU(x)=max(0,x)}$.}
\label{fig:DeepDbar}
\end{figure*}
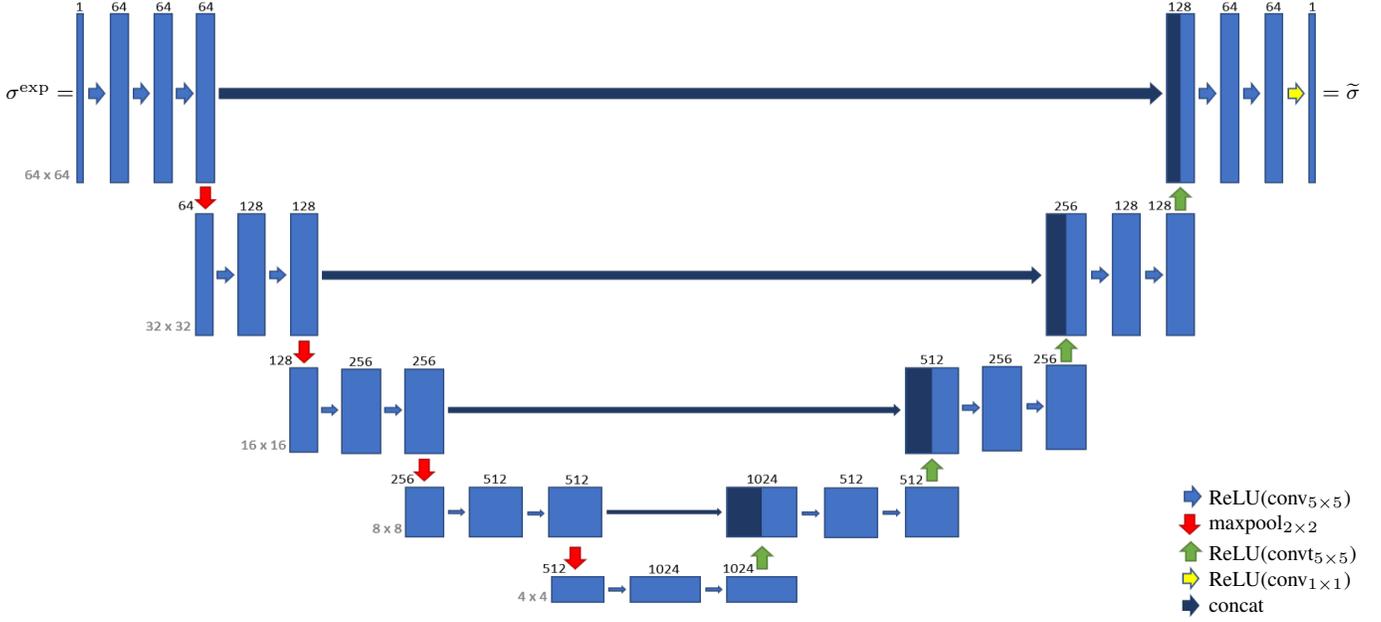

In this work we use this `Born' approximation $\texp$ to the scattering data, first presented in \cite{Isaacson2004}, as it allows the D-bar method to solve the EIT problem fast enough to be considered `real-time' \cite{Dodd2014} and is robust against noisy data.  
The main steps in the algorithm are outlined below:
\begin{tcolorbox}
\[\footnotesize
\begin{array}{c}
\text{Current/Voltage Data}\\
\left(\Lambda_\sigma,\Lambda_1\right)
\end{array} 
\overset{1}{\longrightarrow }
\begin{array}{c}
\text{Scattering Data}\\
\texp(k)
\end{array}
\overset{2}{\longrightarrow }
\begin{array}{c}
\text{Conductivity}\\
\sigma(z)
\end{array}
\]

\vspace{1em}

\noindent{\bf Step 1:} For each $k\in\C\setminus\{0\}$, evaluate the approximate scattering data
\small
\begin{equation}\label{eq:texp}
\texp(k)=\begin{cases}
\int_{\bndry} e^{i\bar{k}\bar{z}}\left(\Lambda_\sigma - \Lambda_1\right) e^{ikz} dS(z), & 0<|k|\leq R\\
0 & |k|>R.
\end{cases}
\end{equation}

\normalsize
\vspace{1em}

\noindent{\bf Step 2:} For each $z\in\Omega$, solve the D-bar equation \eqref{eq:dbark} using the integral equation
\small
\begin{equation}\label{eq:dbark-sol}
\muexp(z,\kappa) = 1+ \frac{1}{4\pi^2}\int_{\C}\frac{\texp(k)e(z,-k)}{(\kappa-k)\bar{k}}\overline{\muexp(z,k)}\;d\kappa_1d\kappa_2,
\end{equation}
\normalsize
and recover the approximate conductivity 
\begin{equation}\label{eq:sigExp}
\sigexp(z)=\left[\muexp(z,0)\right]^2.
\end{equation}
\end{tcolorbox}

%

\trev{
\subsection{Robustness of D-bar Methods for EIT}\label{sec:Robustness}
Recent studies \cite{Murphy2009,Hamilton2017_PhysMeas2} suggest that D-bar based reconstruction methods for 2D EIT are robust to incorrect electrode locations and boundary shape.  This robustness holds for absolute, as well as time-difference, imaging with both images behaving similarly to incorrect boundary shape and electrode locations.  This may be due to the fact that incorrect domain modeling leads to EIT data from a DN map that is only possible for an anisotropic conductivity, even when the true conductivity is isotropic.  While the anisotropic conductivity cannot be recovered uniquely, one can recover a unique isotropization, $\sqrt{\det(\sigma)}$, of the matrix-valued anisotropic conductivity, interpreted as a deformation of the true anisotropic conductivity by isothermal coordinates.  In \cite{Henkin2010, Hamilton2014a}, it is proved that the equations in the D-bar reconstruction methods are identical for anisotropic and isotropic EIT data, helping to explain why D-bar methods have still produced quality images even on anisotropic conductivities and imprecisely known boundary shapes.   Here we focus on absolute images.}

\section{Deep D-bar}\label{sec:DeepDbar}
The aim of this study is to formulate a real-time reconstruction algorithm for electrical impedance tomography that produces sharp and robust \trev{absolute EIT} images. To achieve this we combine the D-bar algorithm, described in Section \ref{sec:Dbar}, with a convolutional neural network (CNN). This idea relies on a network architecture known as U-Net \cite{Ronneberger2015}, originally developed for image segmentation. It has been shown for several linear inverse problems \cite{Jin2017,Kang2017,Sandino2017,Antholzer2017,Hauptmann2018} that this particular network structure can be modified to successfully remove artefacts in medical image reconstructions. The basic recipe is to use a fast and simple reconstruction algorithm to obtain corrupted images and then train the network to remove those artefacts.  A related study for electrical impedance tomography is \cite{Martin2017}, where the authors used artificial neural networks \trev{(ANNs)} to post-process initial reconstructions from one step of a linear Gauss-Newton algorithm \trev{for 3D time-difference EIT imaging}. \trev{Our approach is fundamentally different as it recovers absolute EIT images.}

The network architecture we have chosen relies on the established U-Net \cite{Ronneberger2015}, which consists of a  multilevel decomposition and several skip connections to avoid singularities in the training procedure, see Figure \ref{fig:DeepDbar} for an illustration of our specific architecture. \trev{The original purpose of U-Net was image segmentation. This is very similar to our application, where the main goal is to identify organ boundaries and deconvolve the reconstruction, hence the output of our network is a sharpened image.  Therefore, we believe that the U-Net architecture is a suitable choice for the purpose of EIT imaging, since the multilevel structure can deal efficiently with the non-linearity and sharpening over large image areas. Additionally, as discussed in \cite{Wiatowski2017}, pooling layers leads to translational invariance, which is important to reduce locational bias in the reconstruction process and detect injuries not present in the training set.}  As a modification to the original architecture we needed to increase the convolutional filter size to $5\times5$ (compared to $3\times3$), presumably to deal with the nonlinearity of the reconstructions and enforce consistency of the reconstructions.  We would like to note, that in contrast to the studies in \cite{Jin2017,Kang2017,Sandino2017,Antholzer2017,Hauptmann2018}, where the authors learn a residual update to the initial reconstruction\trev{, we train the network to produce a single sharpened version of the input.}

\subsection{Training of the network}
Given the true conductivity $\sigma$, we simulate measurement data, as will be described in Section \ref{sec:modMeas}, and reconstruct the approximate conductivity $\sigma^{\exp}$ with the D-bar method outlined in \ref{sec:Dbar}. Since the reconstruction step \eqref{eq:dbark-sol} in the D-bar algorithm can be done for any $z\in\R^2$ we reconstruct $\sigma^{\exp}$ on the square $[-1,1]^2$ to obtain a square image as input to the network. The resolution is chosen to be $64\times 64$. The ground truth $\sigma$ is similarly extended to $[-1,1]^2$ by extending the background conductivity. 

Having obtained the training set $\{\sigma_i,\sigma_i^{\exp}\}_i$, we train the Deep D-bar network, denoted by $\mathcal{D}_\theta$, for the set of network parameters $\theta$, i.e. the convolutional filters and biases in each convolutional layer. Given the output of the network $\widetilde{\sigma}=\mathcal{D}_\theta(\sigma^{\exp})$ we seek to minimize the $\ell^2$-error of network output to phantom, given by the loss
\[
\mathrm{loss}(\widetilde{\sigma}):=\|\widetilde{\sigma} - \sigma\|_2^2. 
\]
The network is implemented with the Python library {TensorFlow} and the optimization is performed for \trev{1,000} epochs in batches of 16, with TensorFlow's implementation of the Adam algorithm and an initial learning rate of $10^{-4}$. The training procedure \trev{takes only} 4 hours on a single Titan XP GPU with 12GB memory.  As we will discuss in the following section, we do not need to perform a transfer training to apply the trained Deep D-bar network to experimental data, the training on simulated data proved to be sufficient.

\section{Experimental Setup and Computational Notes}\label{sec:ExpSetups}
We will demonstrate the new {\it Deep D-bar method} using experimental data from two different EIT machines: ACT4 \cite{Liu2005,Saulnier2007} from Rensselaer Polytechnic Institute (RPI) as well as KIT4 \cite{kit4data} from the University of Eastern Finland (UEF).

%

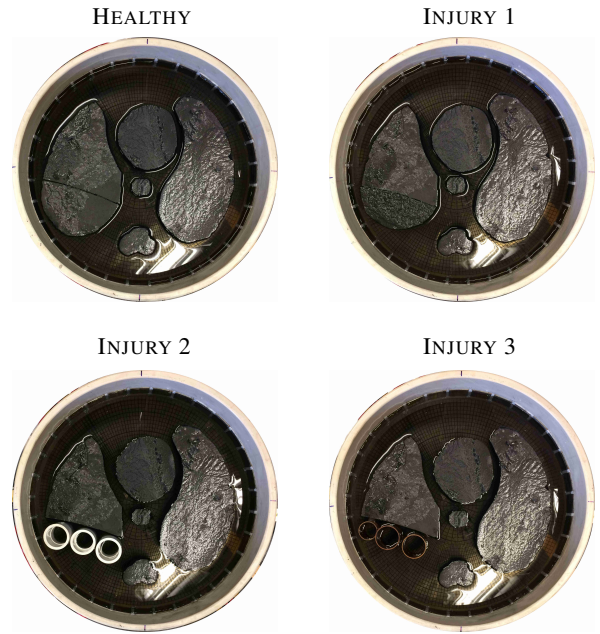
\begin{figure}[h!]
\begin{picture}(250,240)

\put(20,0){\includegraphics[width=100pt]{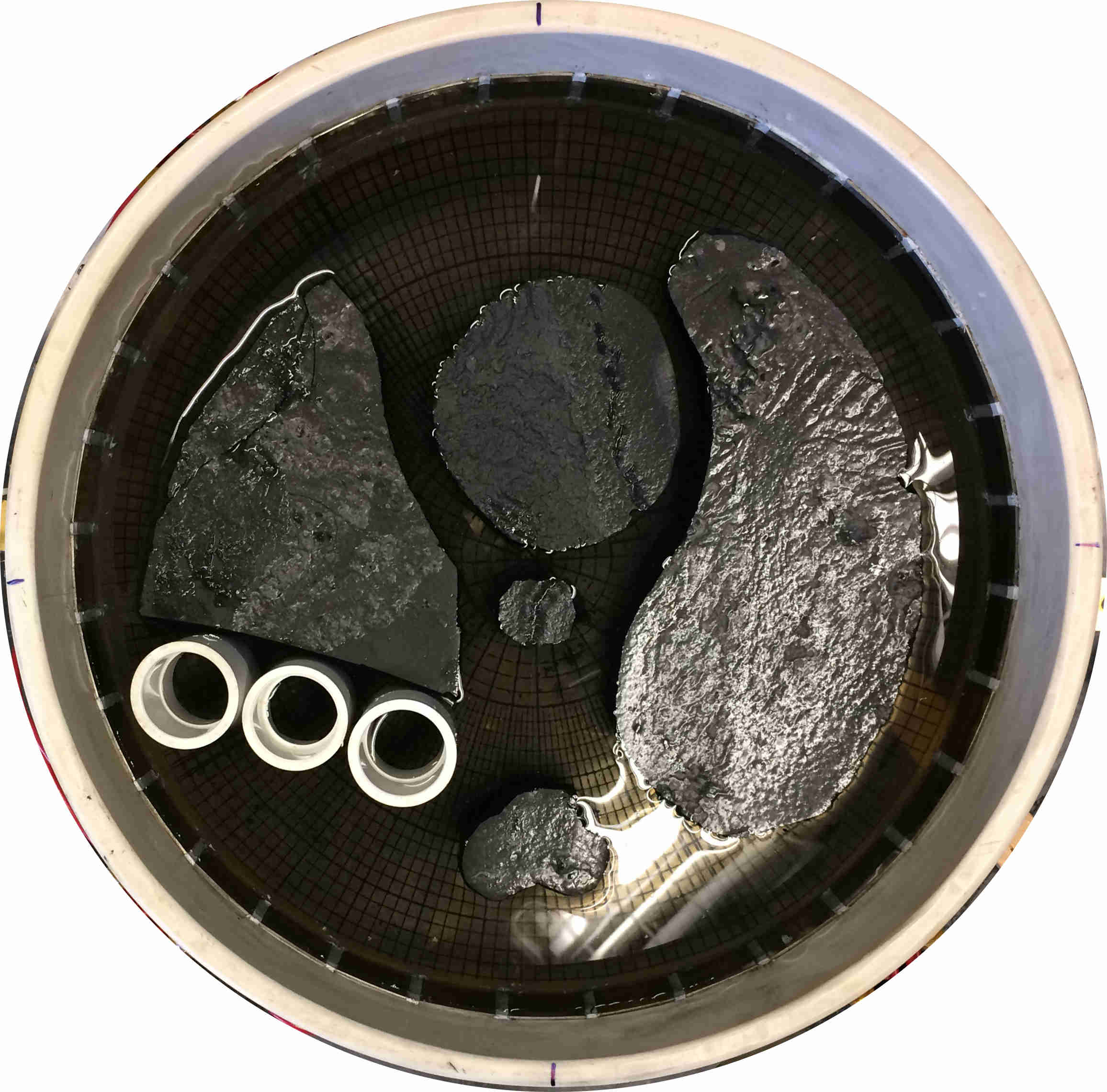}}
\put(140,0){\includegraphics[width=100pt]{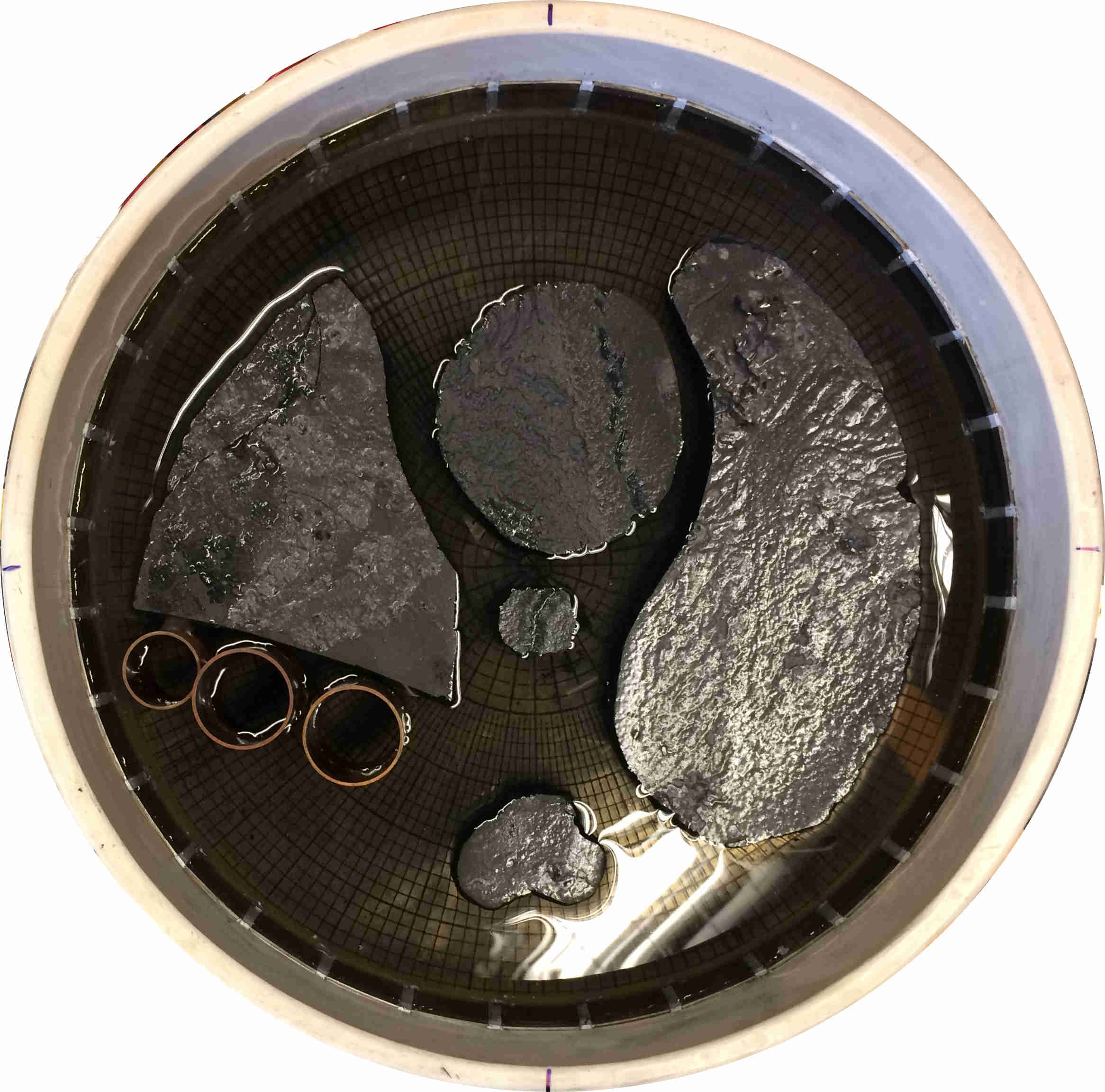}}
\put(20,125){\includegraphics[width=100pt]{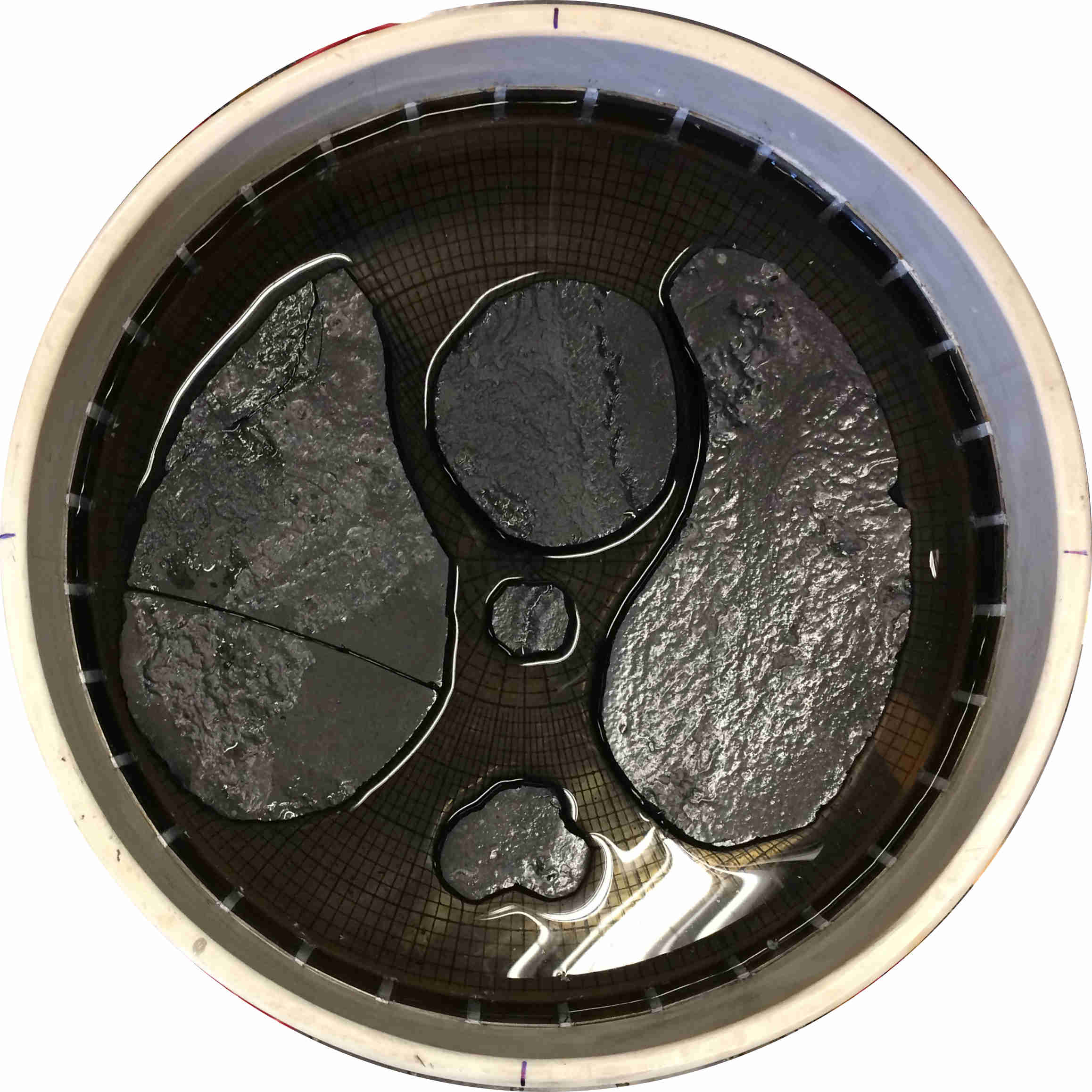}}
\put(140,125){\includegraphics[width=100pt]{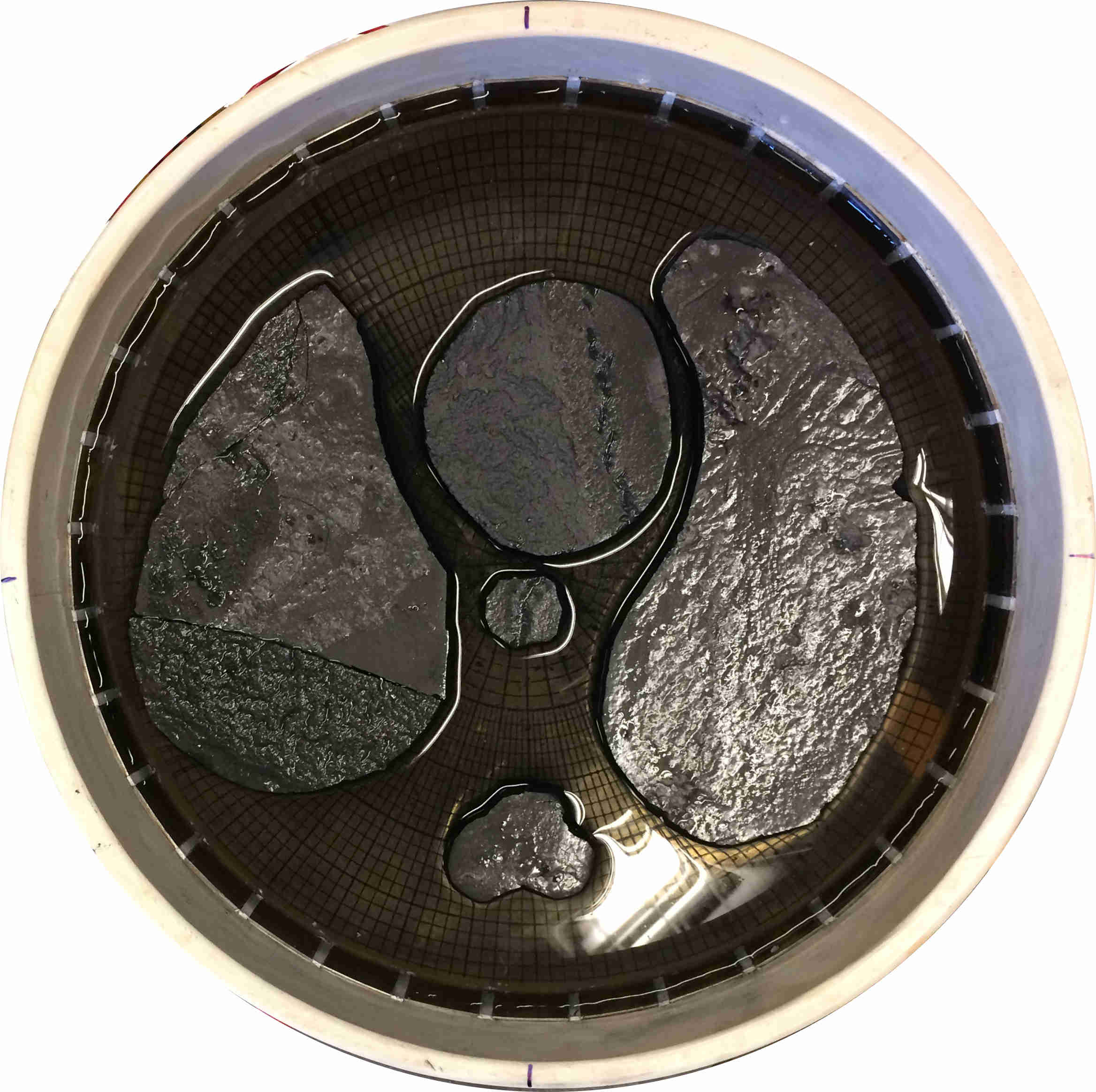}}

\put(50,230){\sc \small {Healthy}}

\put(175,230){\sc \small {Injury 1}}

\put(52,105){\sc  \small {Injury 2}}

\put(175,105){\sc \small {Injury 3}}
\end{picture}
\caption{Experimental Setups for test phantoms taken on the ACT4 system from RPI.  Agar/graphite targets were used to simulate a chest phantom with a heart, two lungs, aorta, and spine.  The first image shows the healthy phantom.  Three injuries are explored:  `Injury 1',  replaced the cut portion of the right lung with agar/graphite of the same conductivity as the heart target to simulate a potential pleural effusion, `Injury 2', replaced the cut portion of the right lung replaced with three plastic tubes, and `Injury 3', replaced the cut portion with three copper tubes.}
\label{fig:RPI_phantoms}
\end{figure}

The ACT4 data uses agar (4\%) based targets with added graphite (10\%) to simulate a heart, two lungs, an aorta, and a spine.  All images are shown in DICOM orientation, meaning that the right lung corresponds to the viewer's left, as if we are looking up through the patient's feet.  Injuries were simulated in the right (DICOM) lung away from the heart by removing a portion of the lung and (1) replacing the missing portion with a piece of agar/graphite with the same conductivity as the heart to simulate an injury such as a pleural effusion, (2) placing three plastic tubes in the missing region to simulate an area of very low conductivity such as a pneumothorax, and (3) replacing the missing portion with three metal tubes.  The experiments are shown in Figure~\ref{fig:RPI_phantoms}.  The approximate conductivities of the targets are displayed in Table~\ref{table:ACT4_setup}.  The admittivity spectrum of the agar/graphite targets were measured on test-cells with Impedimed's SFB-7 bioimpedance meter\footnote{\url{https://www.impedimed.com/products/sfb7-for-body-composition/}}.  Note that the ACT4 system applies voltages and measures currents rather than vice-versa.  In these experiments, trigonometric voltage patterns of maximum amplitude 0.5V (and frequency 3.3kHz) were applied on a circular tank (radius 15cm), with 32 electrodes (width 2.5cm), filled with saline (0.3 S/m) to a height of 2.25cm.  

\begin{table}[h] 
\scriptsize
  \caption{Conductivity Values for ACT4 targets at 3.3kHz}
    \begin{tabular}{l|c|c}
    \hline
&{\sc Measured Values} & {\sc Simulated Values}  \\
& (S/m) & Ranges (S/m) \\
    \hline
    \hline
    {\sc Heart/Aorta} &0.67781 & [0.5, 0.8]\\

    {\sc Lungs/Spine} & 0.056714  & [0.01, 0.2]\\
    
    {\sc Saline Background} & 0.3  & [0.29, 0.31]\\
    
    {\sc Injury 1: Agar/Graphite }&0.67781 & [0.01, 1.5]\\

            {\sc Injury 2: Plastic Tubes} &0  & [0.01, 1.5]\\
            {\sc Injury 3: Copper Tubes} &infinite  & [0.01, 1.5]\\
    \hline
    \hline
    \end{tabular}%
  \label{table:ACT4_setup}%
\end{table}%

\begin{figure}[h!]
\begin{picture}(250,235)

\put(20,0){\includegraphics[width=100pt]{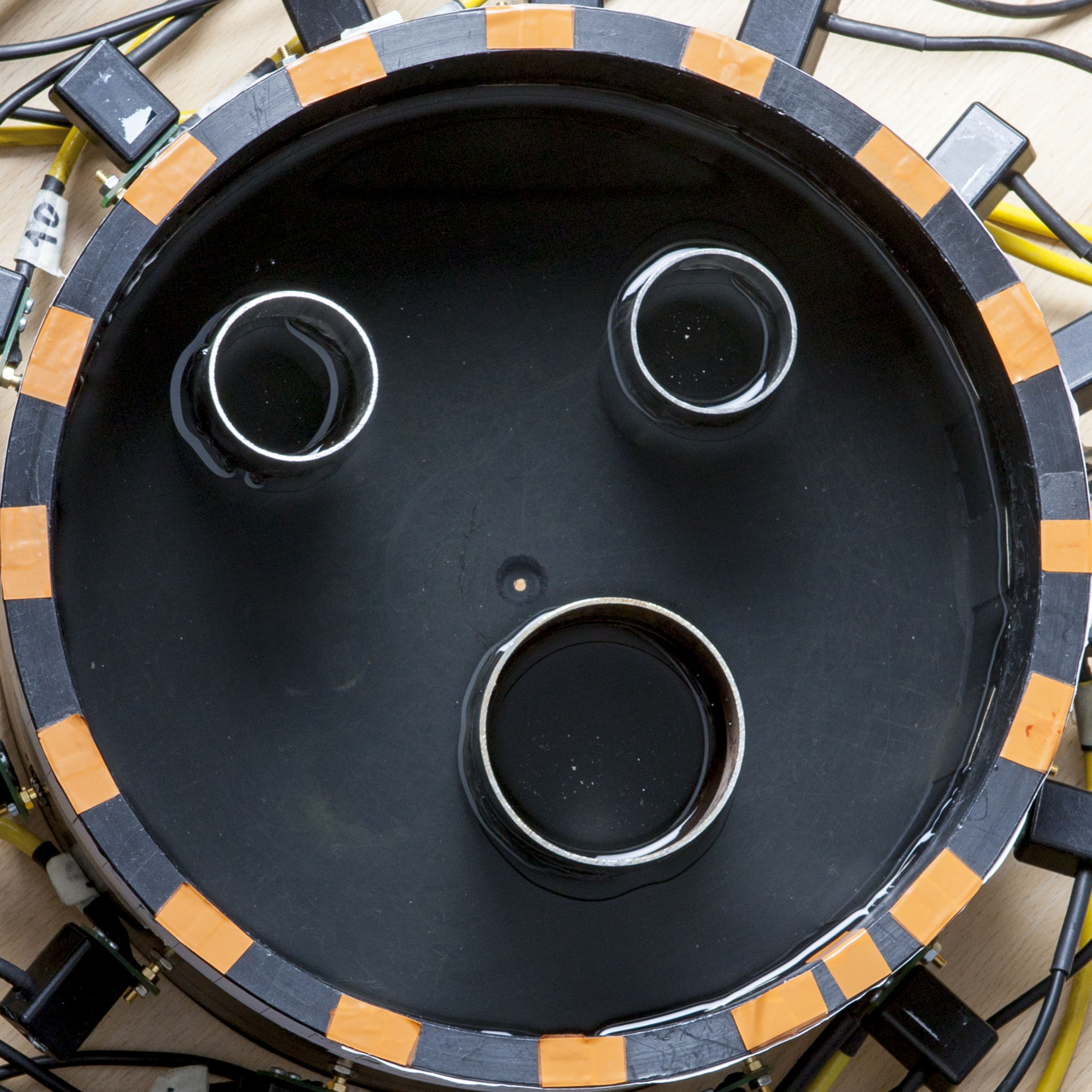}}
\put(140,0){\includegraphics[width=100pt]{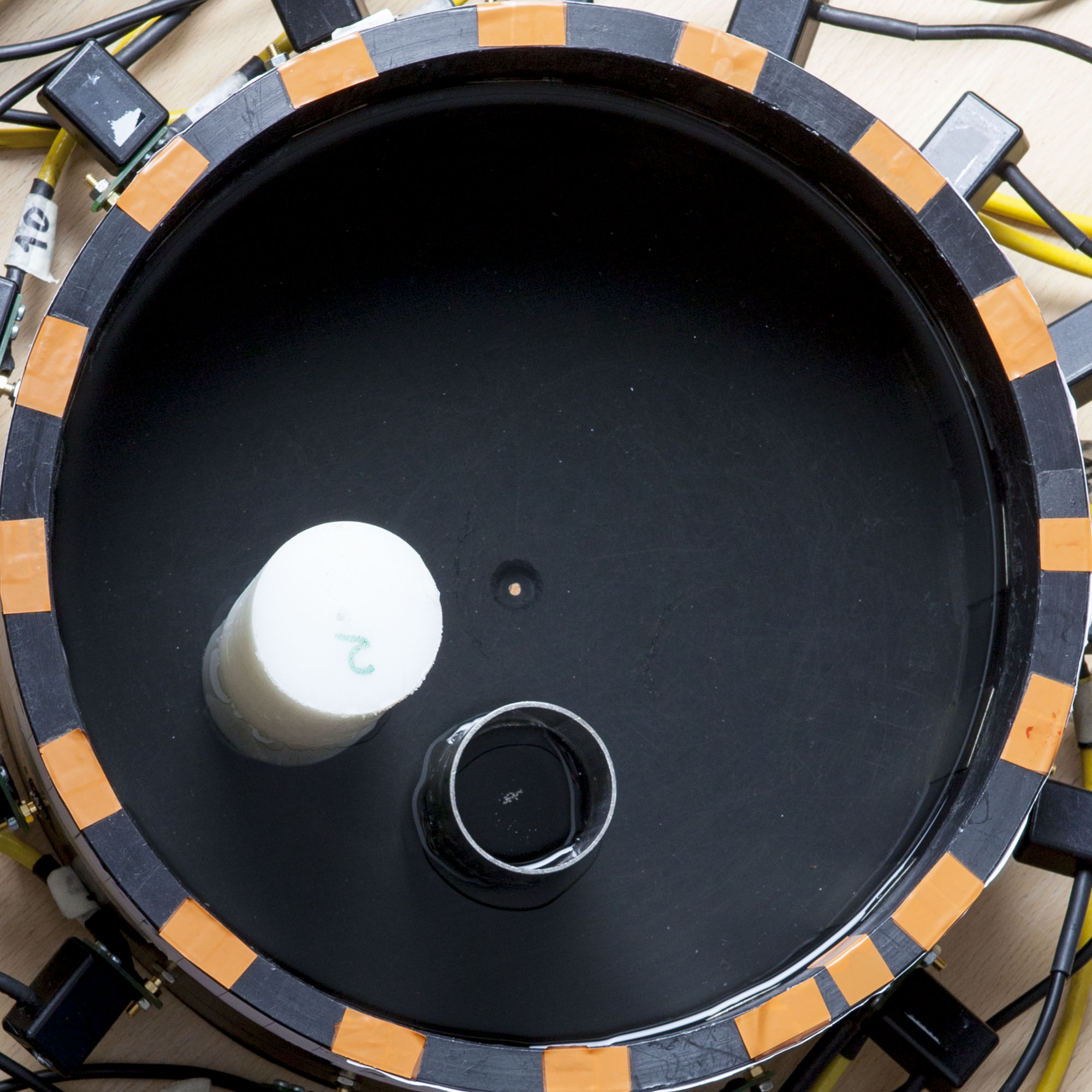}}
\put(20,125){\includegraphics[width=100pt]{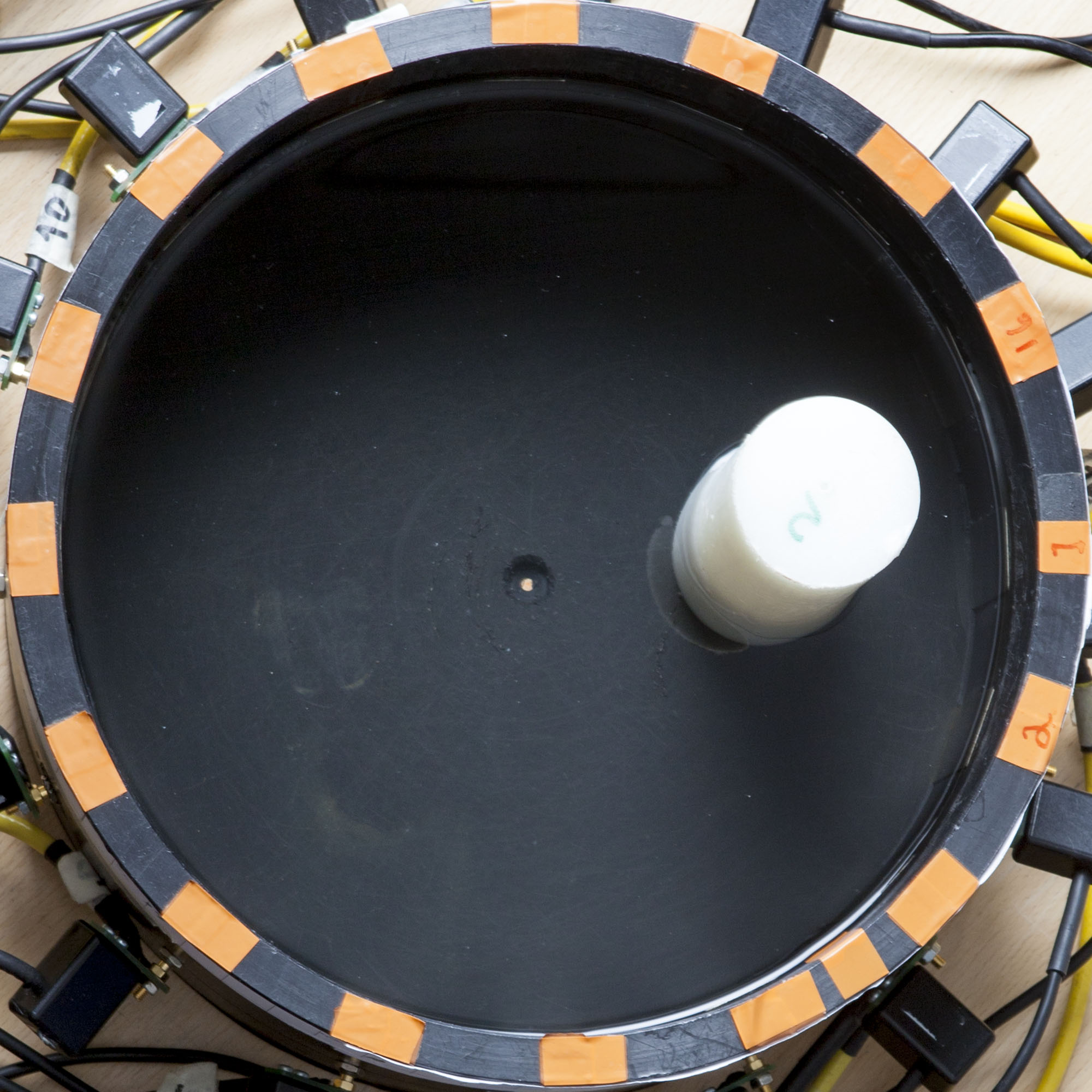}}
\put(140,125){\includegraphics[width=100pt]{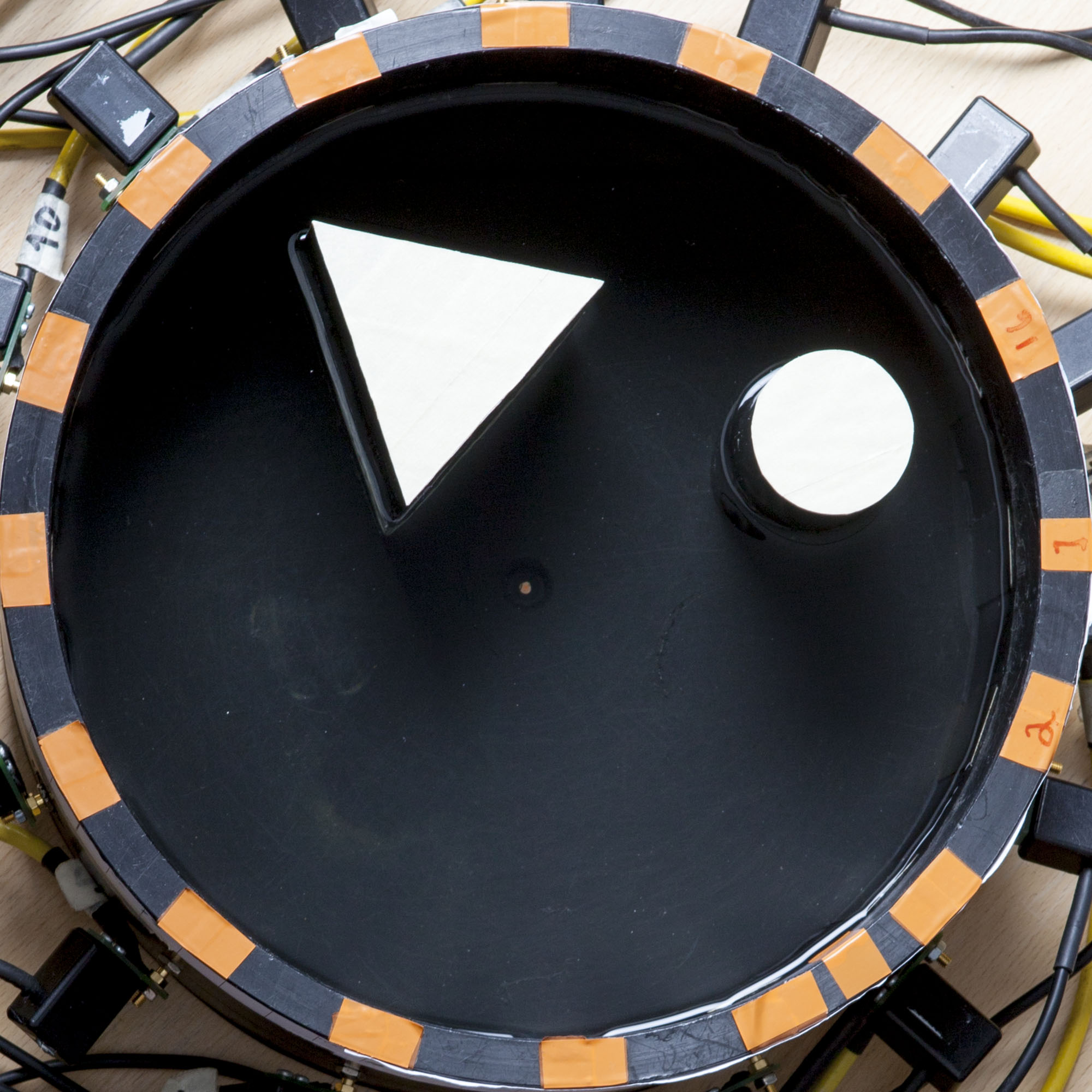}}

\put(42,230){\sc \small {Phantom 1.1}}

\put(165,230){\sc \small {Phantom 2.2}}

\put(42,105){\sc  \small {Phantom 3.4}}

\put(165,105){\sc \small {Phantom 4.4}}
\end{picture}
\caption{Experimental Setups with conductive and resistive targets on the KIT4 EIT system from UEF.  The white objects are made of solid plastic and are resistive.  The hollow circular objects are conductive metal rings.}
\label{fig:KIT4_phantoms}
\end{figure}
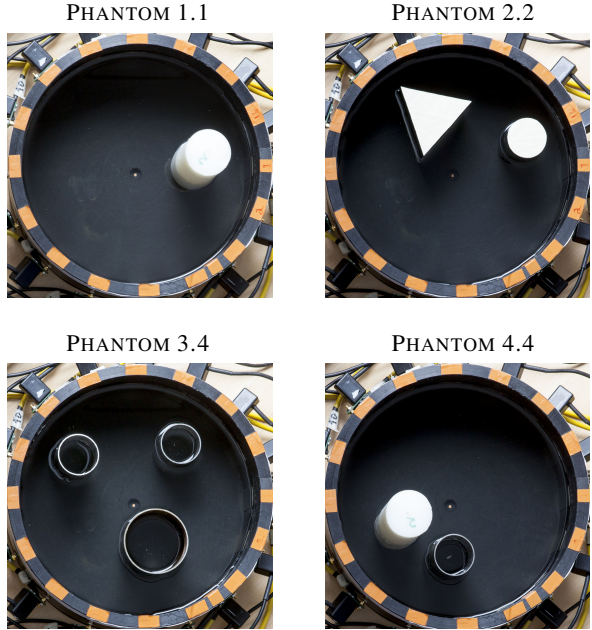

The KIT4 data was taken on a circular tank of radius 14cm with 16 electrodes of width 2.5cm and tap water with conductivity 0.03~S/m filled to a height of 7cm.  Conductive (metal) and resistive (plastic) targets were placed in the tank, as shown in Figure~\ref{fig:KIT4_phantoms}, and adjacent current patterns with amplitude 2mA were applied\trev{ at 1kHz}.  \trevNew{We remark that while this data may not satisfy safety standards for human imaging, it is included for illustrative purposes and potential industrial applications.}

\subsection{Simulation of 2D EIT data}\label{sec:modMeas}
The boundary conditions of \eqref{eqn:condEq} assume a {\it continuum model} for the boundary measurements, completely ignoring \trev{the} discrete positioning of the electrodes.   When simulating the training data, we use a modified version of the continuum model, called the {\it continuum electrode model} introduced in \cite{Hauptmann2017a}, which was developed to simulate realistic electrode data in a continuum setting. In essence, the continuum current/voltage traces are optimally projected onto subsets of the boundary corresponding to the electrode locations.  The training could be done with a more complicated electrode model, such as the {\it Complete Electrode Model} (CEM) \cite{Somersalo1992}, however our simplified continuum electrode model proved sufficient for this proof of concept study.

We aim to represent the ND map as matrix approximation $\mathbf{R}_\sigma$ with respect to an orthonormal basis on the boundary. Let $L$ be an even number of electrodes, then the basis functions are chosen for $n\in\{-L/2,\dots -1,1,\dots,L/2\}$ as
\[
\varphi_n(\theta)=\left\{
\begin{array}{cl}
\frac{1}{\sqrt{\pi}}\sin(n\theta) &\text{ if } n<0,\\
\frac{1}{\sqrt{\pi}}\cos(n\theta) &\text{ if } n>0.
\end{array}\right.
\]
The ACT4 system uses $L=32$ electrodes and the KIT4 system uses $L=16$. The measured voltages are then projected to a continuum trace $g_n$, see \cite{Hauptmann2017a,Hyvoenen2009}, and we obtain the ND matrix ${\mathbf{R}}_\sigma$ by evaluating inner products in $L^2(\bndry)$ as follows
\begin{equation}\label{eqn:matrixApprox}
({\mathbf{R}}_\sigma)_{n,\ell}=(g_n,\varphi_\ell)=\int_{\bndry} g_n(s)\varphi_\ell(s) ds.
\end{equation}
The matrix approximation of the DN map, $\mathbf{L}_\sigma$, is then formed by inverting the ND matrix, i.e. $\mathbf{L}_{\sigma}=\left(\mathbf{R}_{\sigma}\right)^{-1}$.  If the maximal radius $r$ of the domain is not 1, the DN matrix can be scaled by $r$ to correspond to the data that would be obtained if the radius were 1.  Similarly, if $\sigma=\sigma_0\neq1$ near $\bndry$, the DN matrix is scaled by $\frac{1}{\sigma_0}$ to produce the DN matrix that would correspond to $\sigma=1$ near the boundary. If an estimate for  $\sigma_0$ is not available, the best constant conductivity approximation to the data can be formed as described in \cite{Isaacson2004}.  The scaling is undone at the end of the D-bar algorithm by multiplying the conductivity by $\sigma_0$.  \trev{The matrix approximation $\mathbf{L}_1$ to the DN map $\Lambda_1$, required to evaluate the scattering data via \eqref{eq:texp}, is simulated using the constant conductivity $\sigma=1$.}

\subsection{Simulation of Training Data}\label{sec:Training}
Training data for the neural network was created using solely simulated data: one group for the ACT4 data and another group for the KIT4 data.  

\begin{figure}[h!]
\begin{picture}(250,260)

\put(0,0){\includegraphics[height=100pt]{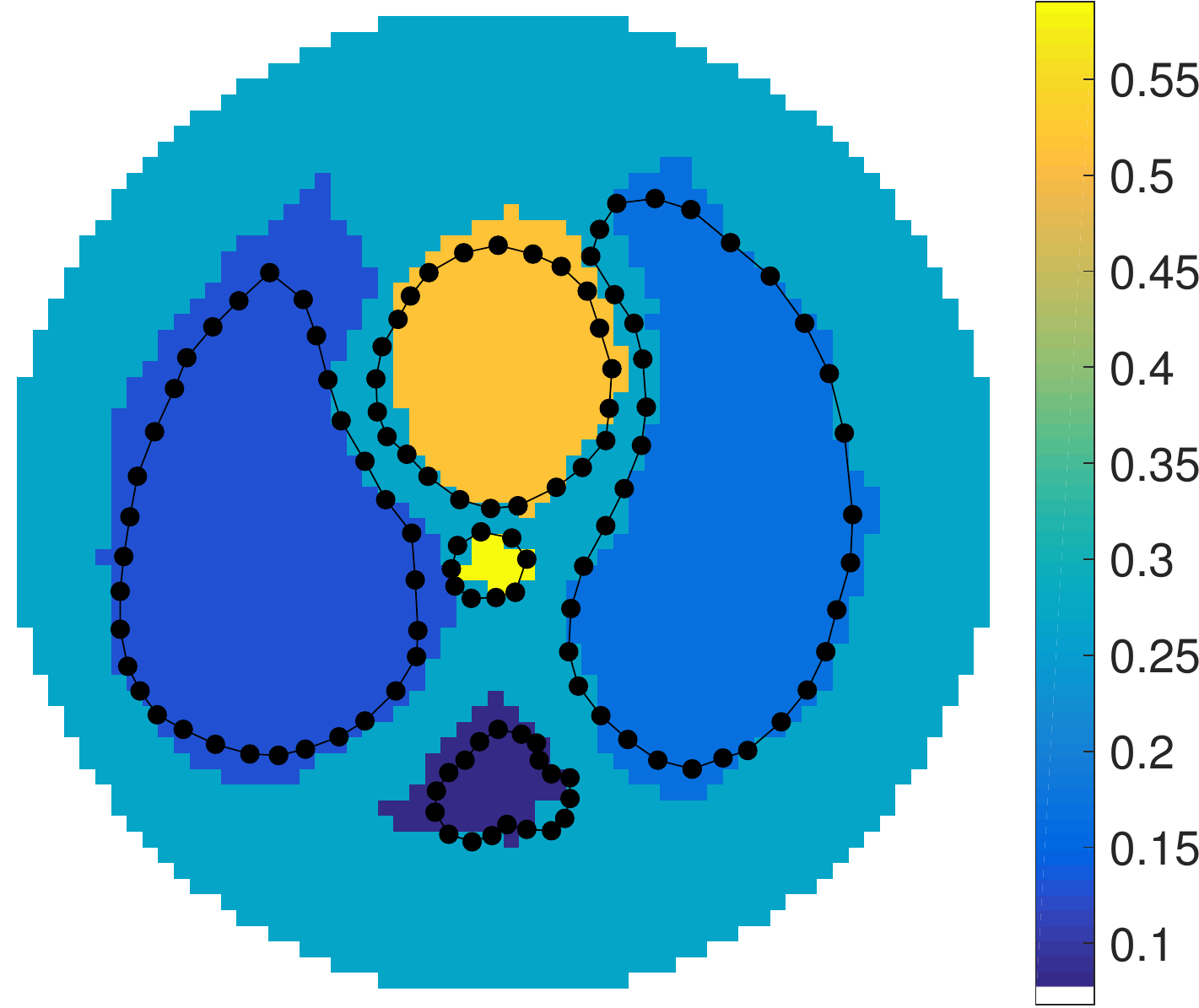}}
\put(130,0){\includegraphics[height=100pt]{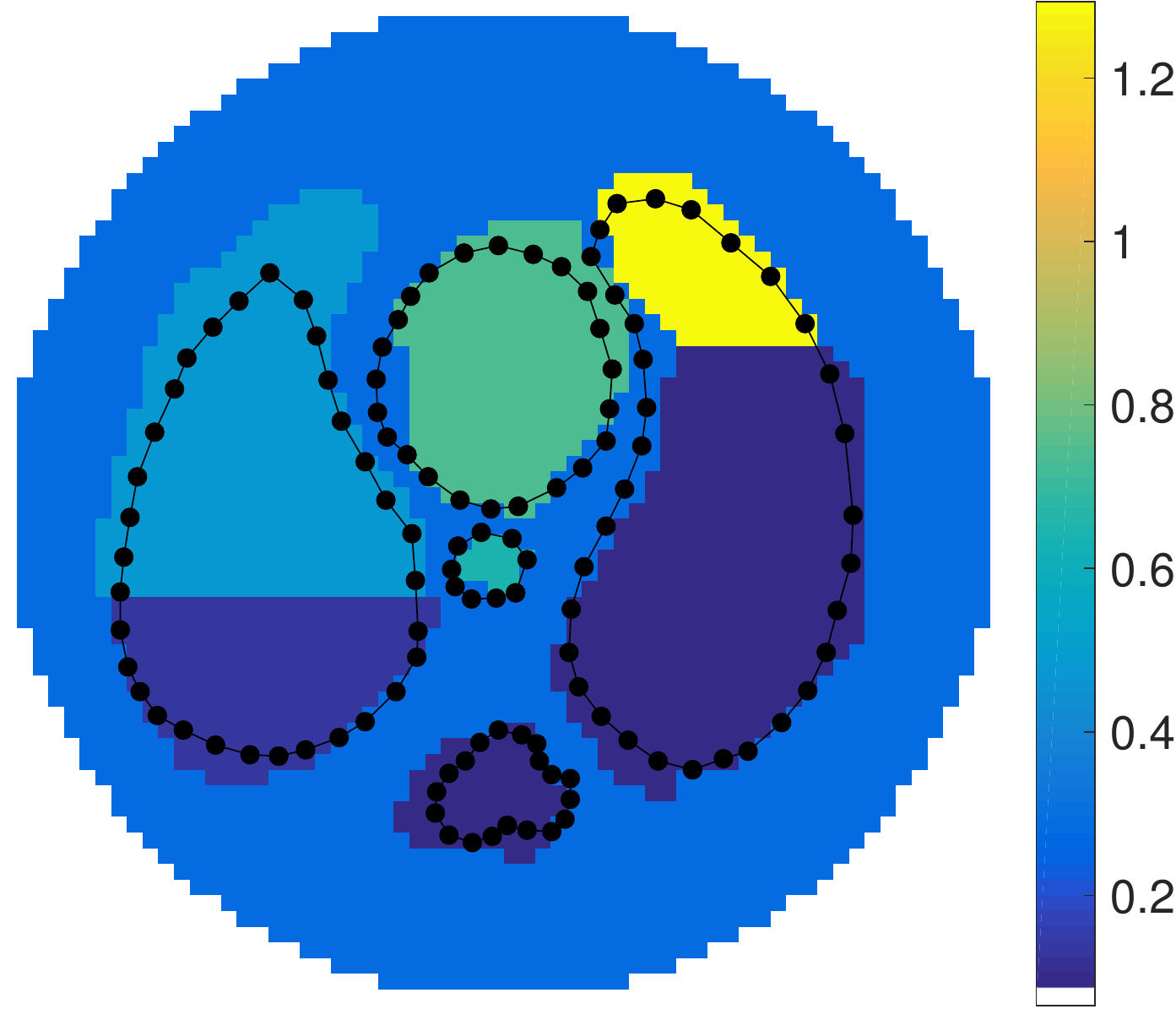}}
\put(0,140){\includegraphics[height=100pt]{HLSA_saline_DeepDbar_DICOM_smaller.jpg}}
\put(130,140){\includegraphics[height=100pt]{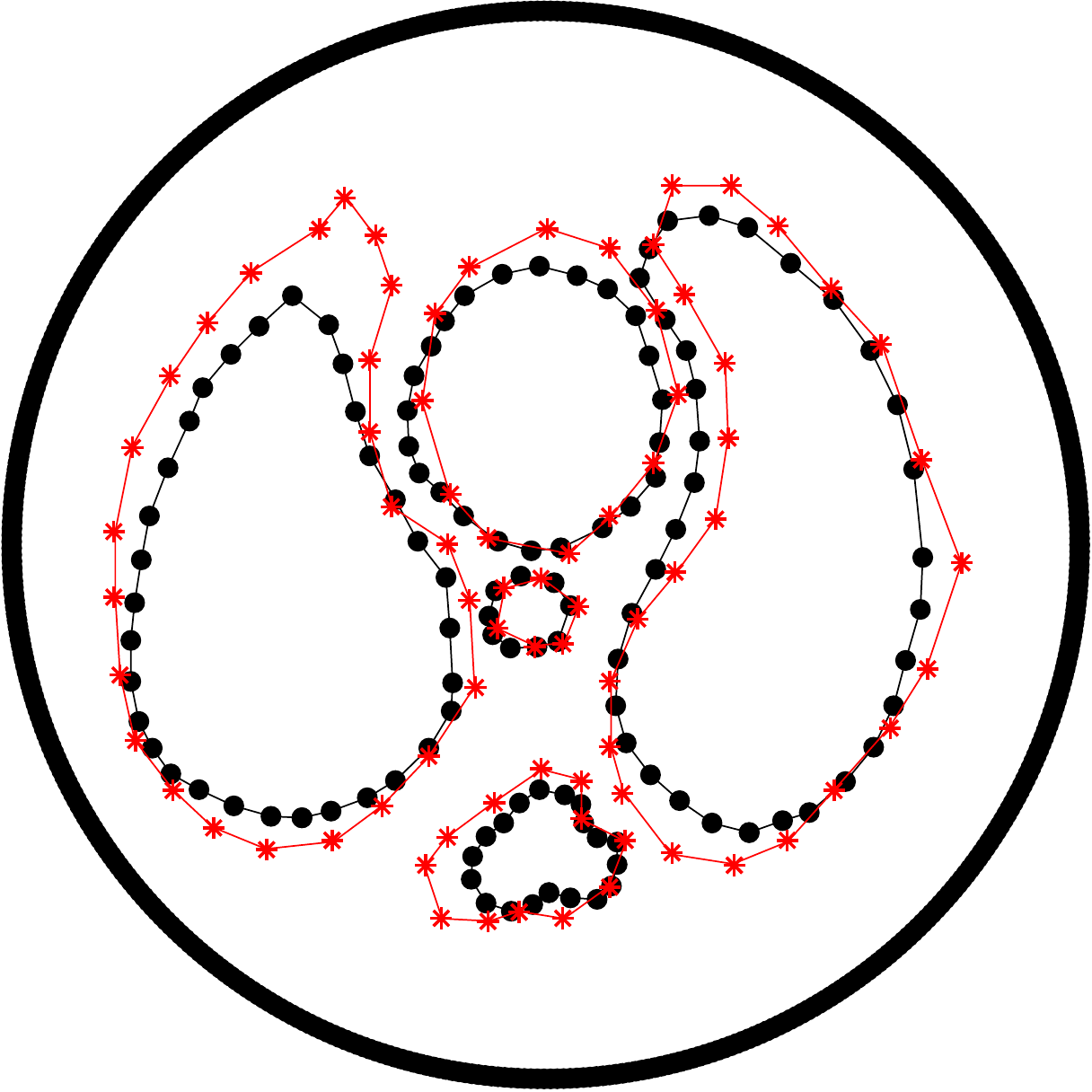}}

\put(20,250){\sc \footnotesize\trev{{ACT4 Healthy}}}

\put(140,255){\sc \footnotesize {True \& Approximate}}
\put(155,245){\sc \footnotesize {Boundaries}}

\put(18,110){\sc \footnotesize{\trev{Noisy Inclusions}}}

\put(145,115){\sc \footnotesize{\trev{Noisy Inclusions}}}
\put(155,105){\sc \footnotesize {with Injuries}}
\end{picture}
\caption{Depiction of the simulation of the training data for the ACT4 experiments of Figure~\ref{fig:RPI_phantoms}.  The first image shows the healthy phantom from which the `true boundary' (black dots) and `approximate boundary' (red stars) were extracted, shown in the second image.  The third and fourth images display sample simulated phantoms using in the training data with and without injuries with the true boundaries overlaid in black dots.}
\label{fig:ACT4demoSims_noisyOrgans}
\end{figure}

The ACT4 training data was created as follows.  Using the `Healthy' image, shown in Figure~\ref{fig:RPI_phantoms} (top left), approximate organ boundaries were extracted by clicking around the targets in the image for the heart, aorta, left lung, right lung, and spine \trev{(Fig.~\ref{fig:ACT4demoSims_noisyOrgans}, top right)}.  Random numbers were generated to decide whether each individual target was included, heart (95\%), aorta (95\%), left lung (90\%), right lung (90\%), spine (100\%).  If a given target was included, white Gaussian noise (25db) was added to the approximate boundary points of the target using the \texttt{awgn} command in {\sc Matlab} to create `noisy' boundary locations.  Figure~\ref{fig:ACT4demoSims_noisyOrgans} \trev{(bottom)} shows the effect of the white noise on the boundary locations.  Noise was added to each target/organ independently.   Conductivities were assigned for each included target by generating a random number from a uniform distribution in the ranges shown in Table~\ref{table:ACT4_setup}, last column.  Elementary injuries were simulated by generating a horizontal dividing line in the lung and assigning randomly generated values in each of the two portions of the divided lung from the uniform distribution of values in [0.01, 1.5]\trev{, see Fig.~\ref{fig:ACT4demoSims_noisyOrgans} bottom right.}  Each lung had an independent chance of such an injury \trev{(30\%)}.  More complex injuries could be simulated but are outside the scope of this study. \trev{A total of 4,096 simulations were performed for the ACT4 training.}

After each conductivity phantom was constructed, the mathematical {\it forward problem} \eqref{eqn:condEq} was solved to recover the corresponding theoretical boundary voltages and currents using a FEM mesh with 65,536 triangular elements using the {\it continuum electrode model} described in Section~\ref{sec:modMeas}.  Relative white noise with variance of $10^{-4}$ was added to the measured voltages. The resulting simulated voltages/currents were used to solve the {\em inverse problem} using the D-bar method described in Section~\ref{sec:Dbar} with a low-pass filtering radius of $R=4.5$ in the scattering domain using the procedure outlined in \cite{Mueller2012} and uniformly spaced $64\by64$ $k$ and $z$-grids on $[-4.5,4.5]^2$ with stepsize $h_k=0.3234$, and  $[-1,1]^2$ with stepsize $h_z=0.0317$, respectively.  A non-uniform cutoff threshold was enforced on the scattering data for frequencies such that $\texp(k)=0$ if either $|\Re(\texp(k)|$ or $|\Im(\texp(k)|$ exceeds 24.  Then, the 4,096 pairs of data in the form of `Truth' and `Low-pass D-bar Reconstruction' were used to train the convolutional neural network described in Section~\ref{sec:DeepDbar}.

\begin{figure}[h!]
\begin{picture}(300,80)

\put(-20,0){\includegraphics[height=70pt]{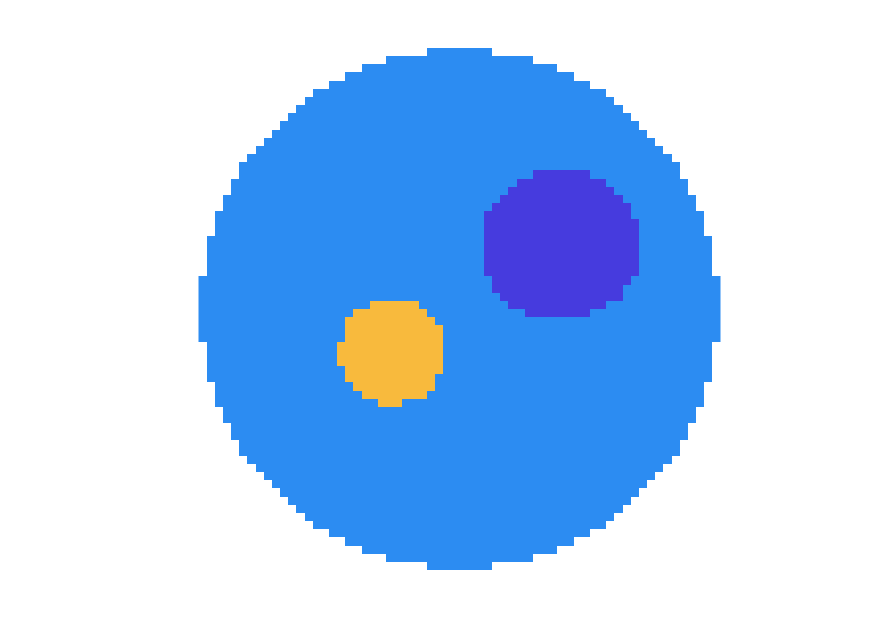}}
\put(60,0){\includegraphics[height=70pt]{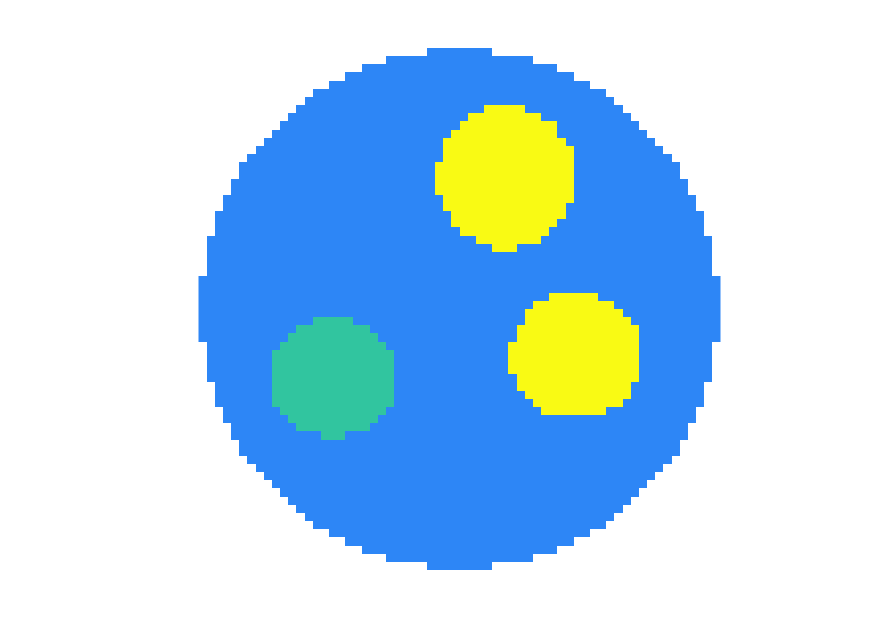}}
\put(140,-10){\includegraphics[height=88pt]{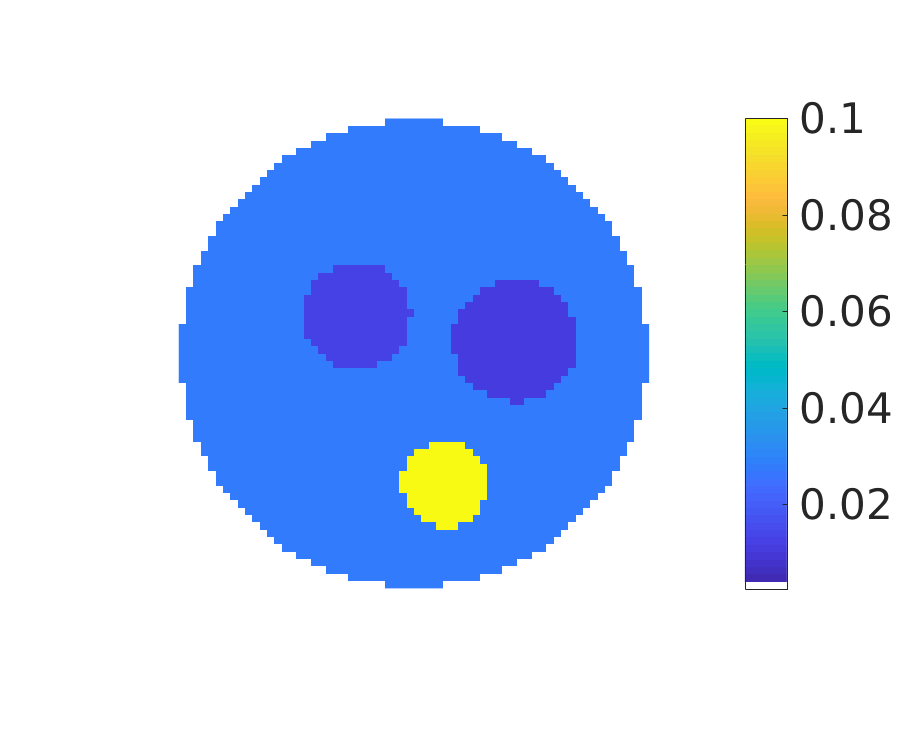}}
\put(65,75){\sc \underline{KIT 4 Sample Phantoms}}

\end{picture}
\caption{Depiction of the simulation of the training data for the KIT4 experiments of Figure~\ref{fig:KIT4_phantoms}.  The images shown are sample simulations of inclusions present in the training data.  Zero to three, non-overlapping, circular inclusions were allowed in each simulation.}
\label{fig:kit4sims}
\end{figure}

Training data for the KIT4 experiments was simulated in a similar manner.  \trev{ In this case, one to three circular inclusions were simulated, with varying radii drawn from the uniform distribution on [0.2, 0.4], with center in [0, 0.6], and an angle in [0, 2$\pi$].  } Inclusions were not allowed to overlap and each inclusion had an equal probability of being  `conductive' or `resistive', and values were assigned from the Uniform distributions [0.05,0.12] and [0.005,0.015] in S/m, respectively.  The conductivity of the background \trev{was} drawn from [0.027, 0.033] S/m.  The conductivity ranges for the inclusions were determined from initial higher-pass (larger filtering radius in the scattering domain) {D-bar} reconstructions of the KIT4 current/voltage data.  While we note that the infinite (metal) conductors should have much higher conductivities, this range was observed in the initial testing and proved sufficient for classifying objects as conductors/resistors in our study.  Note that such infinite conductors/resistors violate the theory of {D-bar} which requires inclusions to have non-negative conductivities bounded away from zero and infinity.  Nevertheless, the method provides useful conductor/resistor information.  The same $k$ and $z$ grids were used in the {D-bar} reconstructions of the KIT4 example as in the ACT4 example.  However, we reduced the non-uniform cutoff threshold of the scattering data from 24 to 8 to reduce oscillatory artefacts that can result from noise in the scattering data for higher frequencies. Figure~\ref{fig:kit4sims} shows sample phantoms used in the training data for the KIT4 example.  A total of 4,096 simulations were performed and pairs of `Truth' vs. `Low-pass D-bar Reconstruction' used to train the network.

\subsection{Computational Notes for D-bar Reconstructions from Experimental EIT Data}\label{sec:DbarCompNotes}
The D-bar reconstructions from the experimental ACT4 and KIT4 data were computed in the same manner as the simulated data case described above in Section~\ref{sec:Training} with the exception of the formation of the DN matrices $\mathbf{L}_\sigma$ and $\mathbf{L}_1$, which now come from discrete vs. continuous measurements.   For convenience, for both the ACT4 and KIT4 data, we synthesized the current/voltage measurements that would have occurred if \trev{orthonormal trigonometric current patterns had been applied,} by using a change of basis.  \trev{Define $t^m_\ell$ to be the value of the $m$-th normalized trigonometric current pattern on the $\ell$-th electrode, following Isaacson et al. \cite{Isaacson2004},
\[\footnotesize t(\ell,m)=t^m_\ell:=\begin{cases}
\sqrt{\frac{2}{L}}\cos(m\theta_\ell), & m=1,\ldots, \frac{L}{2}-1\\
\sqrt{\frac{1}{L}}\cos(m\theta_\ell), & m=\frac{L}{2}\\
\sqrt{\frac{2}{L}}\sin((m-L/2)\theta_\ell), & m=\frac{L}{2}+1,\ldots, L-1
\end{cases}\normalsize\]
for $\ell=1,\ldots, L$.  This ensures that the matrix of current patterns are orthonormal allowing the solution method of DeAngelo and Mueller \cite{DeAngelo2010}.   Alternative methods such as Gram-Schmidt based approaches could also be used to produce a matrix of orthonormal currents.  The corresponding voltages are formed using a change of basis matrix relating the physical currents and the normalized trig currents.   As ACT4 applies voltages and measures currents, we must enforce orthonormality of the currents.  Similarly, the KIT4 data used adjacent current patterns which are not orthogonal and must be converted.}

Using the approach introduced in \cite{Isaacson2004}, the $(m,n)$ entry of the ND matrix $\mathbf{R}_\sigma$ was formed as the discrete inner product
\begin{equation}\label{eq:NDmap_experimental}
\mathbf{R}_\sigma (m,n) = \sum_{\ell=1}^L \frac{\phi_\ell^m v^n_\ell}{|e_{\ell}|},\quad \begin{array}{c}
1\leq m,n\leq L-1\\
1\leq \ell\leq L
\end{array}
\end{equation}
where $\phi^m$ denotes the $m$-th normalized current pattern, $v^n$ the voltage resulting from applying the $m$-th current pattern (normalized such that the voltages sum to zero and are scaled by the $\ell^2$ norms of the applied current patterns), and $|e_\ell|$ denotes the area of the $\ell$-th electrode.  This formula holds for a system with $L$ electrodes where $L-1$ linearly independent current patterns have been applied.  

The discrete DN matrix $\mathbf{L}_\sigma$ was then formed by $\mathbf{L}_{\sigma}=\left(\mathbf{R}_{\sigma}\right)^{-1}$ and scaled by $\frac{r}{\sigma_0}$ as described above.  As the scattering data $\texp$ \eqref{eq:texp} requires the difference in DN matrices $\left(\mathbf{L}_{\sigma}-\mathbf{L}_1\right)$, the discrete DN matrices $\mathbf{L}_1$ must be formed for both the ACT4 and KIT4 system with $L=32$ and $L=16$ electrodes, respectively.   To this end, the conductivity equation in \eqref{eqn:condEq} was solved\trev{, using $\sigma=1$,} with boundary conditions given by the {\it Complete Electrode Model} \cite{Somersalo1992} on a FEM mesh \trev{with triangular elements (ACT4: 4,149; KIT4: 3,493)} simulating injected trigonometric current patterns of amplitude 1mA and non-optimized constant contact impedances of 0.00024 \trevNew{$\Omega\cdot$mm}.  

The scattering data $\texp$ \eqref{eq:texp} was evaluated as a simple Simpson's rule approximation\trev{, following \cite{DeAngelo2010},}
\[\texp(k) \approx  \begin{cases}
\frac{2\pi}{L}\left[e^{i\bar{k}\overline{\mathbf{z}}}\right]^T\phi\left[\mathbf{L}_\sigma-\mathbf{L}_1\right]\mathbf{e}^{\mbox{\tiny $\psi$}}(k)  & 0<|k|\leq R\\
0 & |k|>R
\end{cases}\]
where $\mathbf{z}$ is the vector of the positions of the centers of the electrodes, $T$ denotes the traditional matrix transpose, and 
\[e^{\mbox{\tiny $\psi$}}_\ell(k):=\sum_{\ell}^L a_j(k)\phi^j_\ell \approx {e^{ikz_\ell}} \]
is the expansion of the asymptotic behavior $\psi\sim e^{ikz}$ at the center of the $\ell$-th electrode, $z_\ell$, in the basis of normalized applied current patterns $\phi$.  Then, the d-bar equation \eqref{eq:dbark-sol} can be solved using Fast Fourier Transforms using Vainikko's method \cite{Vainikko2000}.  The interested reader is referred to \cite{Mueller2012} for further details.

\section{Results}\label{sec:results}
We now demonstrate the effect of the {\em Deep D-bar method} on simulated\trev{,} as well as experimental\trev{,} data \trev{for absolute EIT imaging}.
\subsection{Reconstructions from Simulated Data}
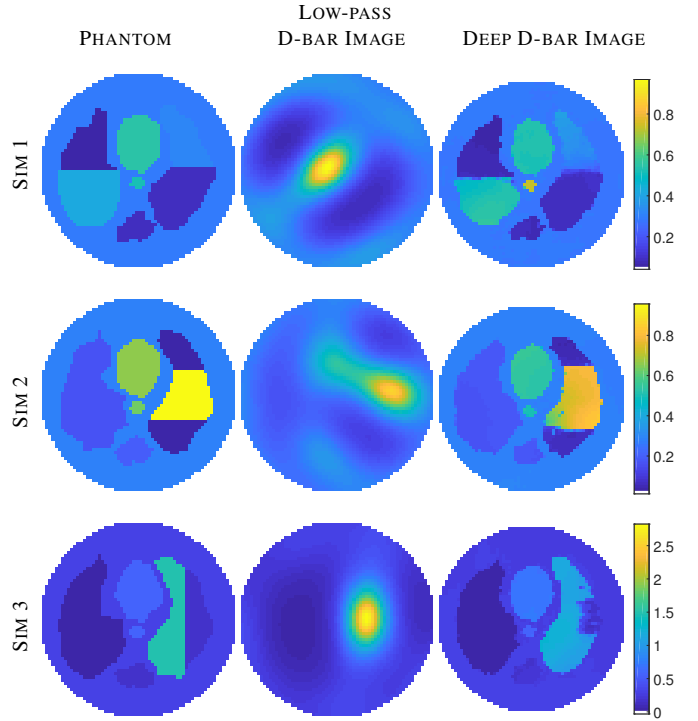
\begin{figure}[h!]
\centering
\begin{picture}(300,265)

\put(10,170){\includegraphics[width=75pt]{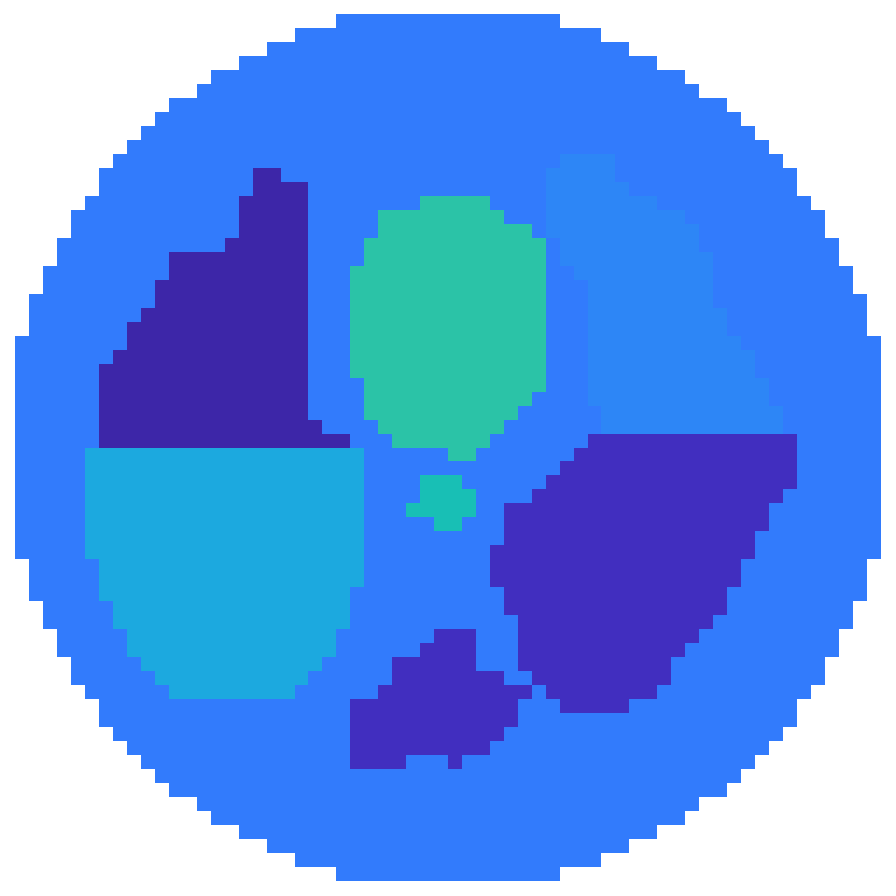}}
\put(85,170){\includegraphics[width=75pt]{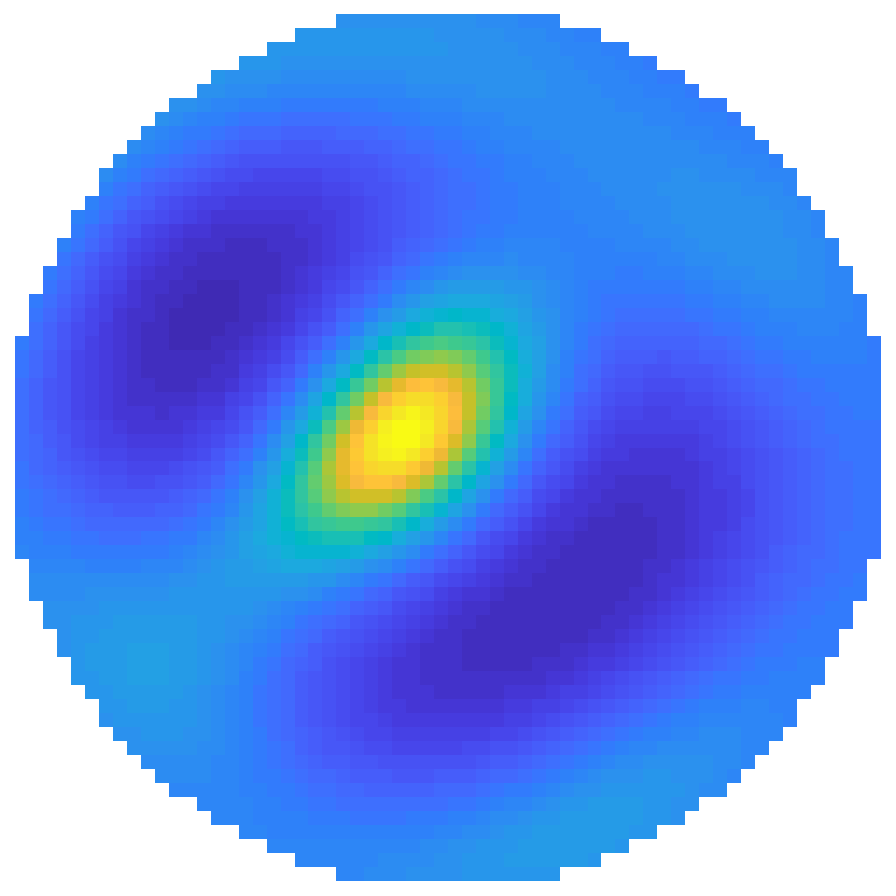}}
\put(160,170){\includegraphics[width=90pt]{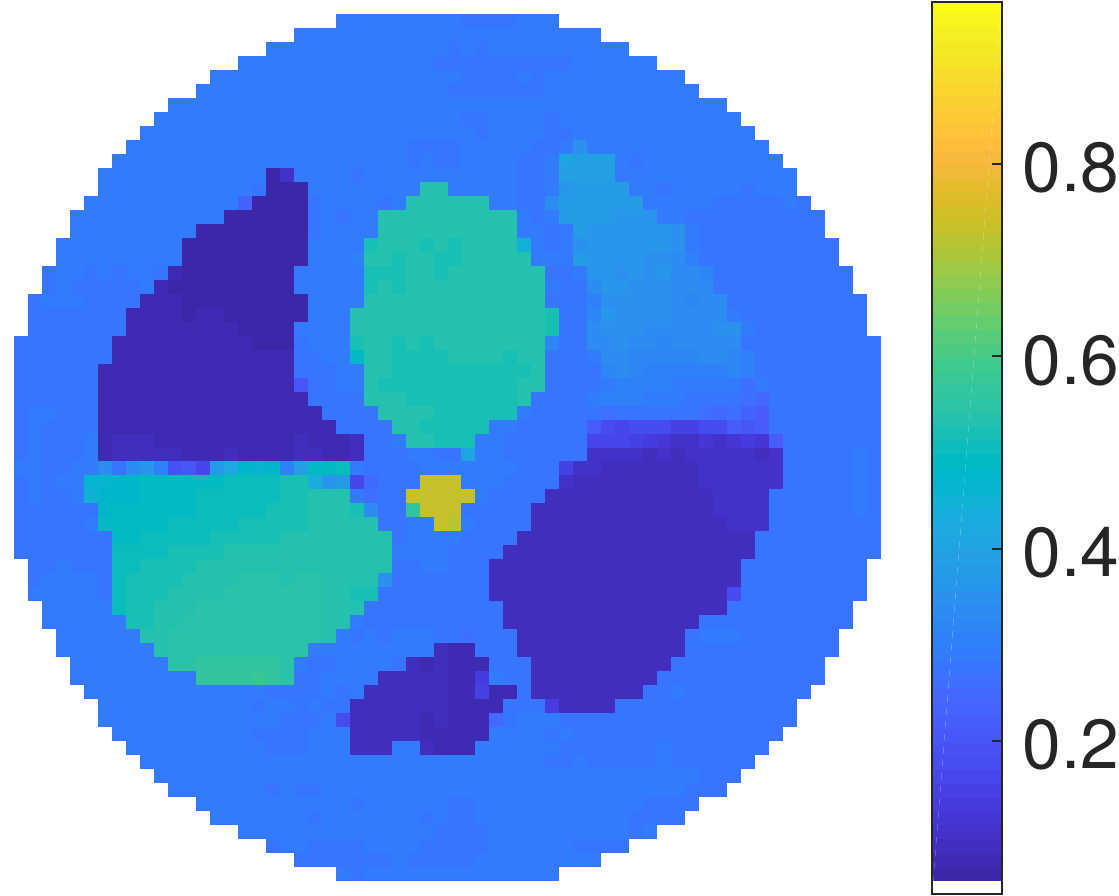}}

%
\put(10,85){\includegraphics[width=75pt]{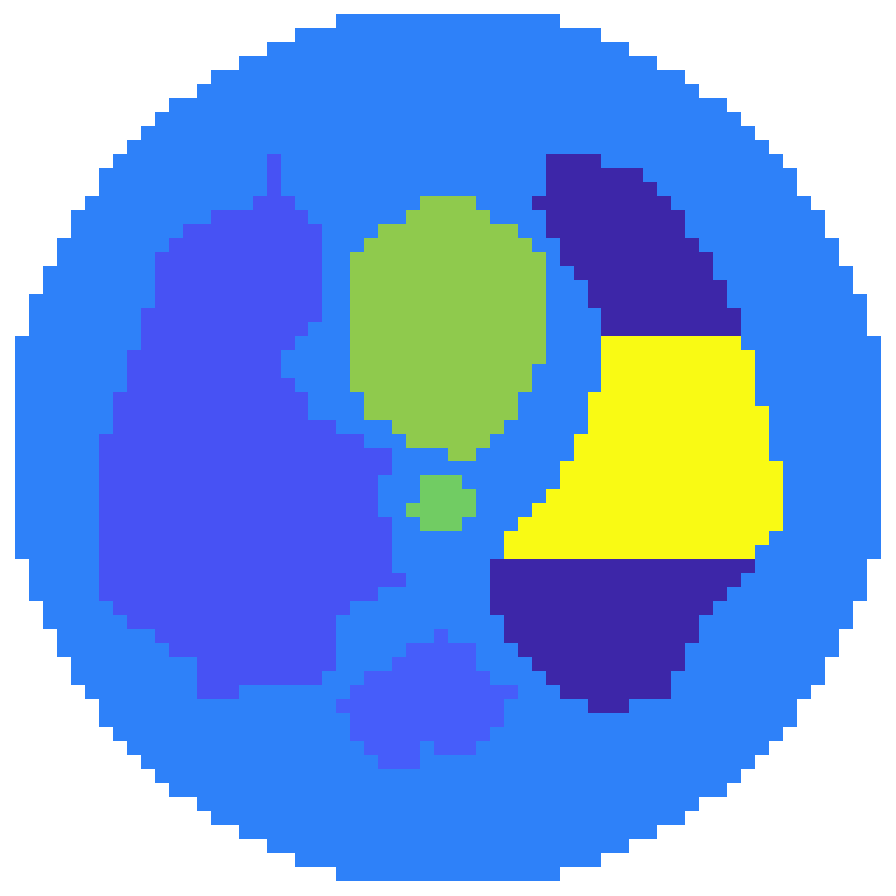}}
\put(85,85){\includegraphics[width=75pt]{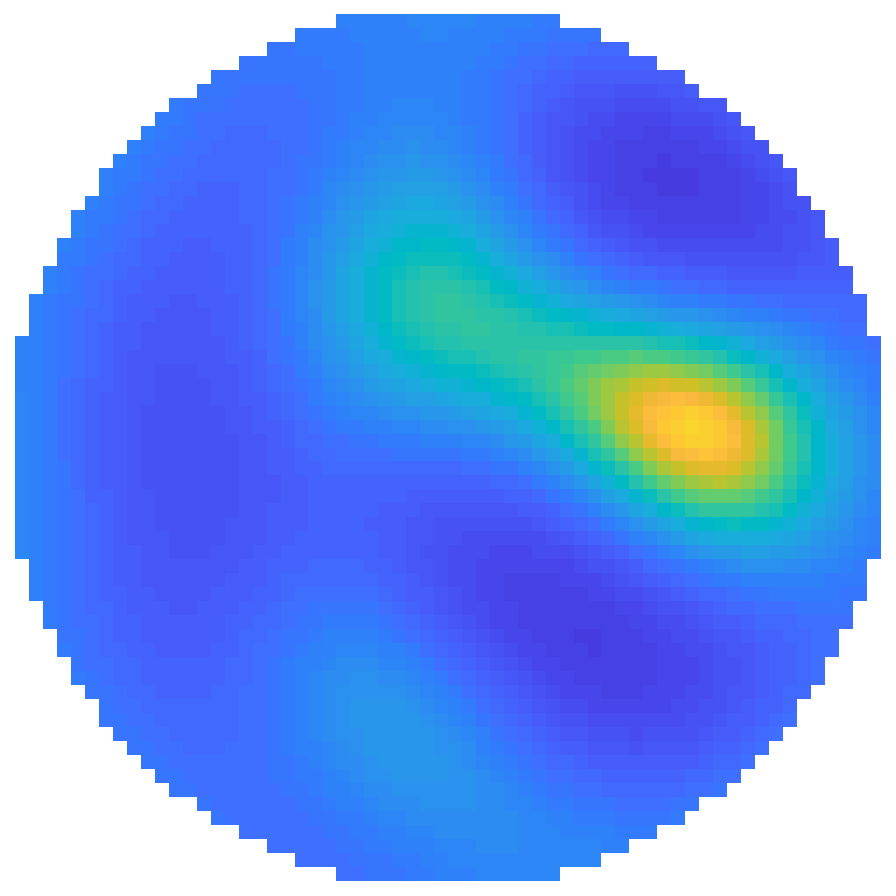}}
\put(160,85){\includegraphics[width=90pt]{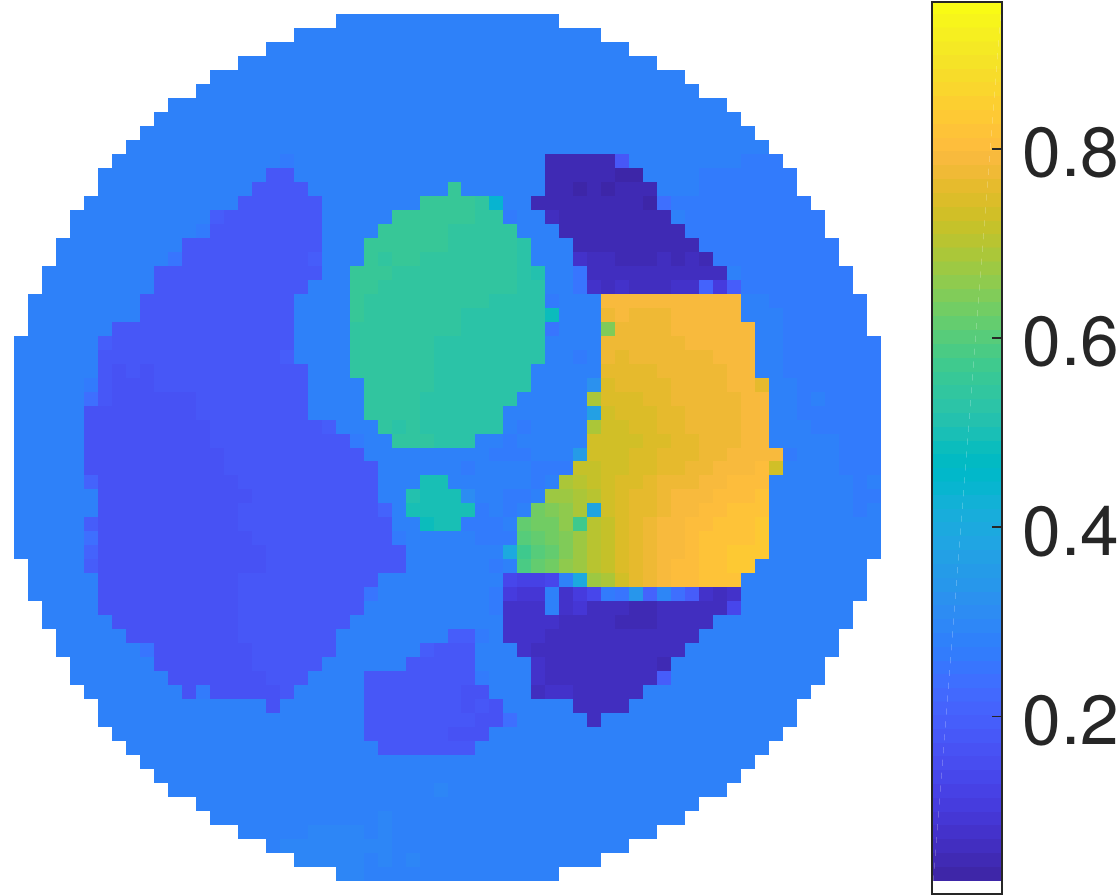}}
%
\put(10,0){\includegraphics[width=75pt]{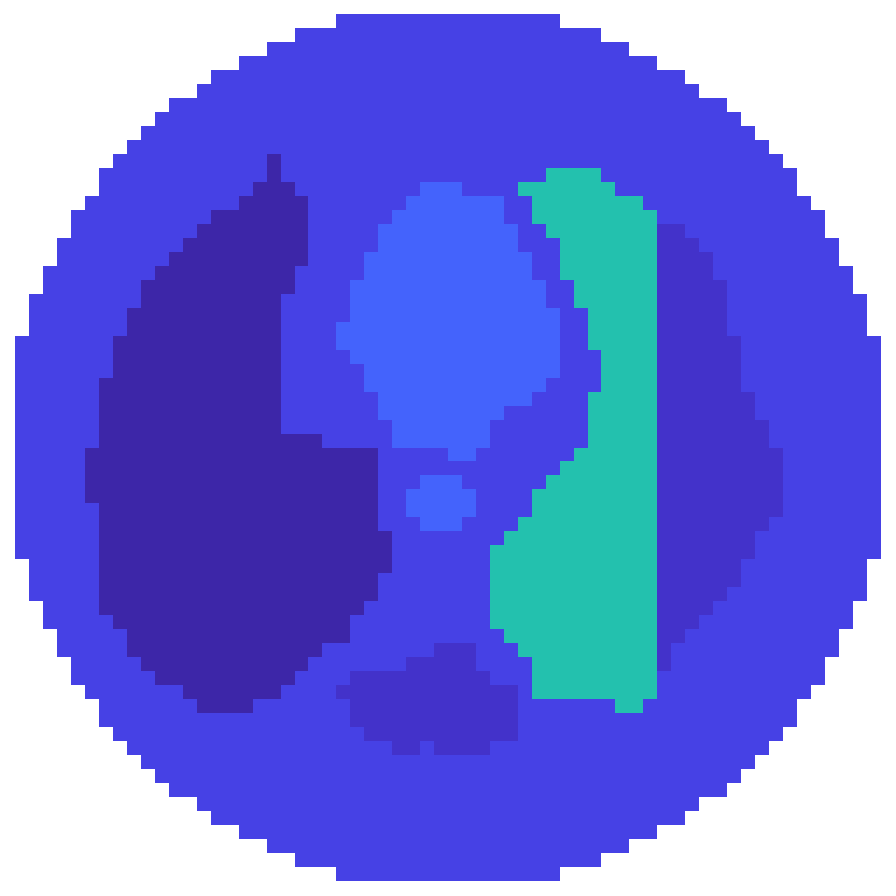}}
\put(85,0){\includegraphics[width=75pt]{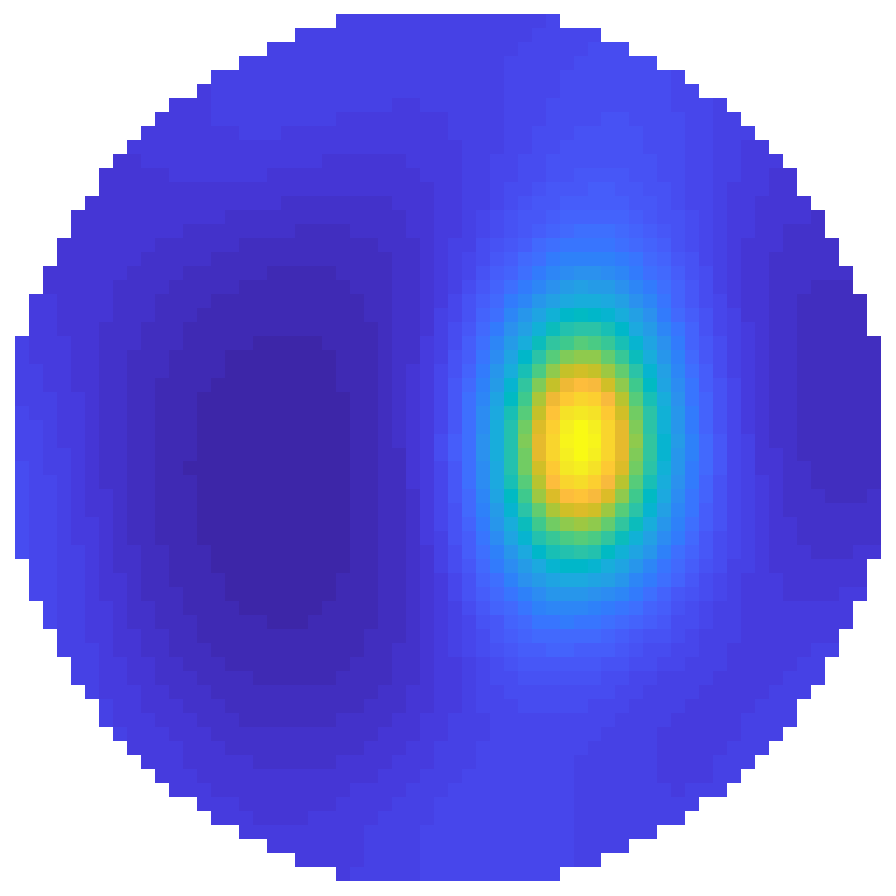}}
\put(160,0){\includegraphics[width=90pt]{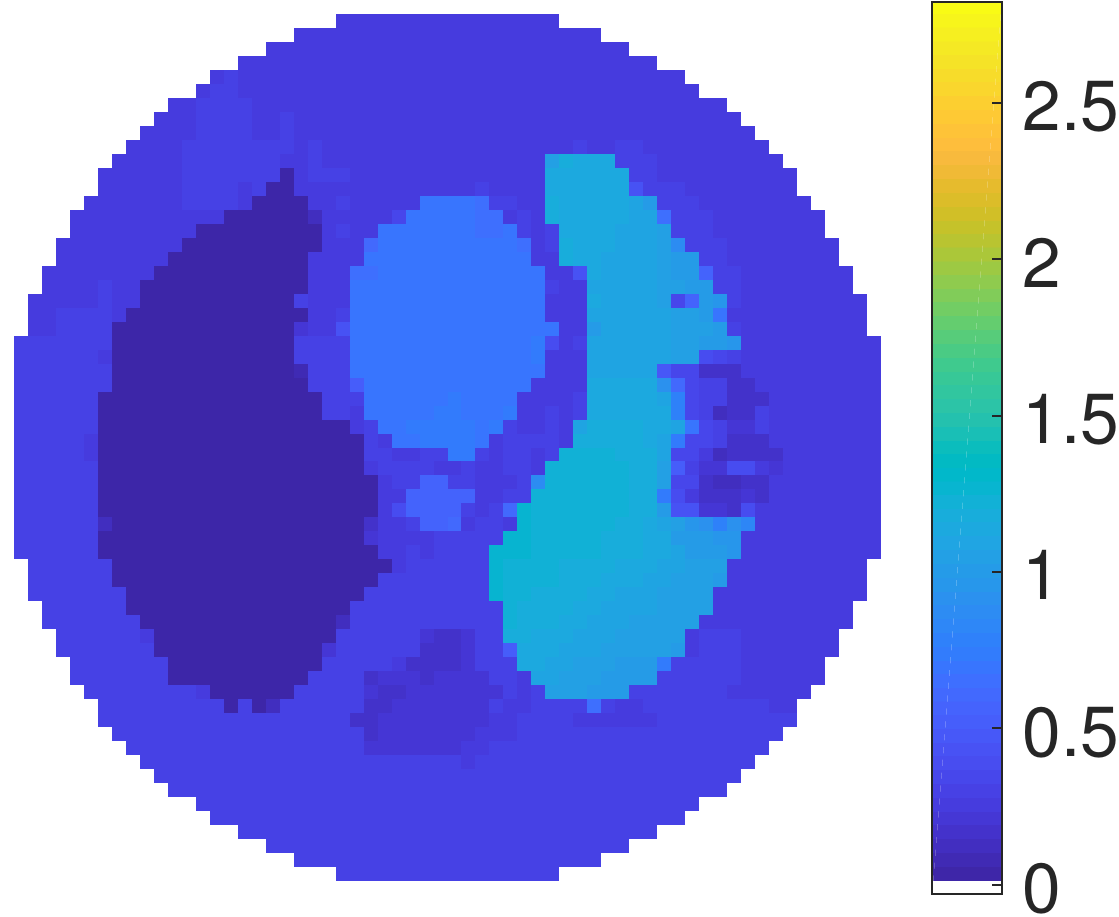}}

\put(25,255){\footnotesize \sc Phantom}

\put(108,265){\footnotesize \sc Low-pass}
\put(100,255){\footnotesize \sc D-bar Image}

\put(170,255){\footnotesize \sc  Deep D-bar Image}

\put(0,25){\rotatebox{90}{{\sc \footnotesize Sim 3}}}
\put(0,110){\rotatebox{90}{{\sc \footnotesize Sim 2}}}
\put(0,200){\rotatebox{90}{{\sc \footnotesize Sim 1}}}

\end{picture}
\caption{\label{fig:ACT4_Results_simData} Results for simulated test data from the ACT4 geometry. The phantom in the first row conforms with the training data and the phantoms in the second and third row\trev{s} include pathologies not supported by the training data. The initial D-bar image is compared to the Deep D-bar image.  The D-bar images, on the full square are used as the `input' images for the CNN.  Images are displayed here on the circular geometry of the tank, for presentation only.  Each row is plotted on its own scale.}
\end{figure}

\begin{figure}[h!]
\centering
\begin{picture}(300,265)

\put(10,170){\includegraphics[width=75pt]{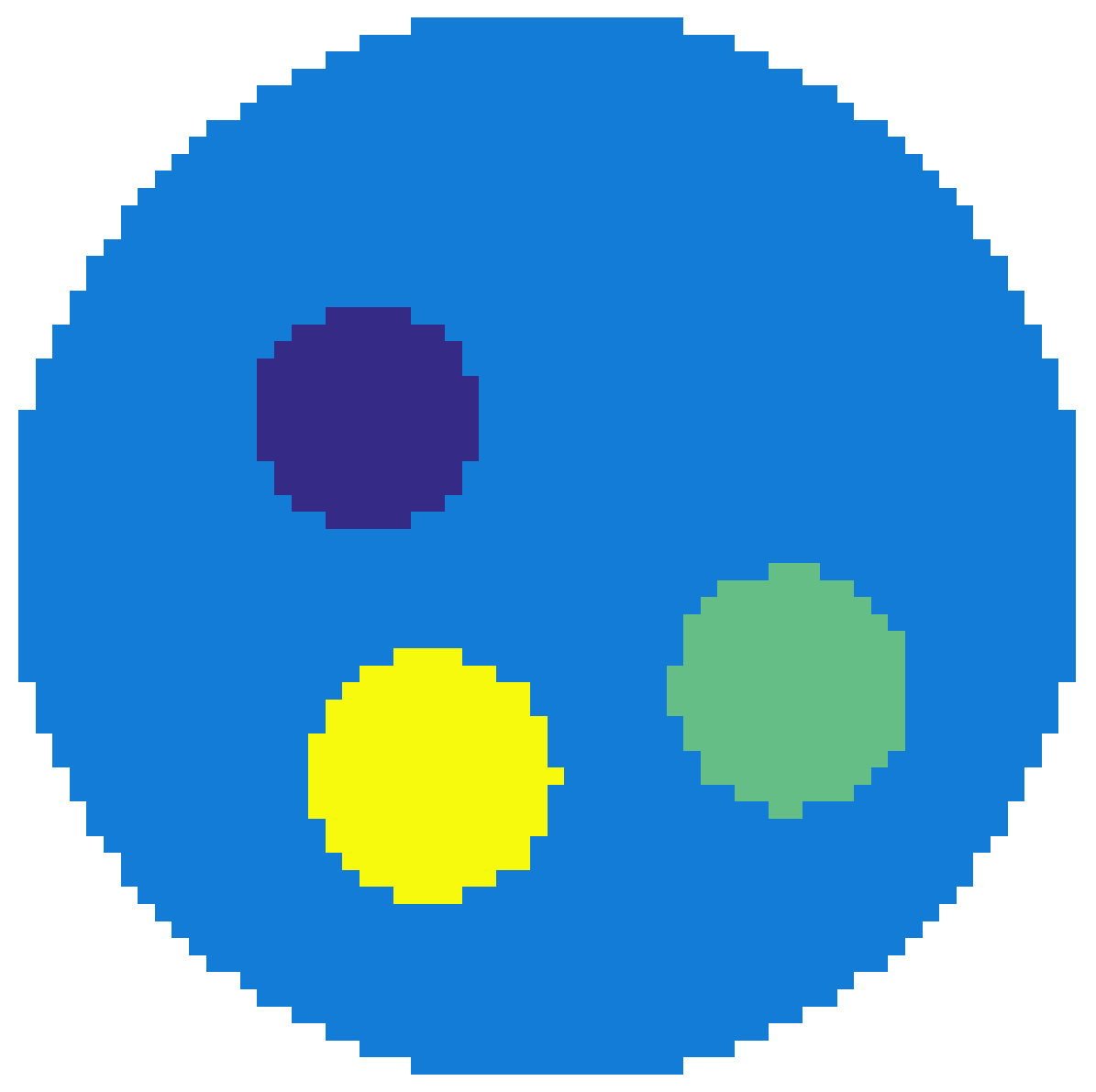}}
\put(85,170){\includegraphics[width=75pt]{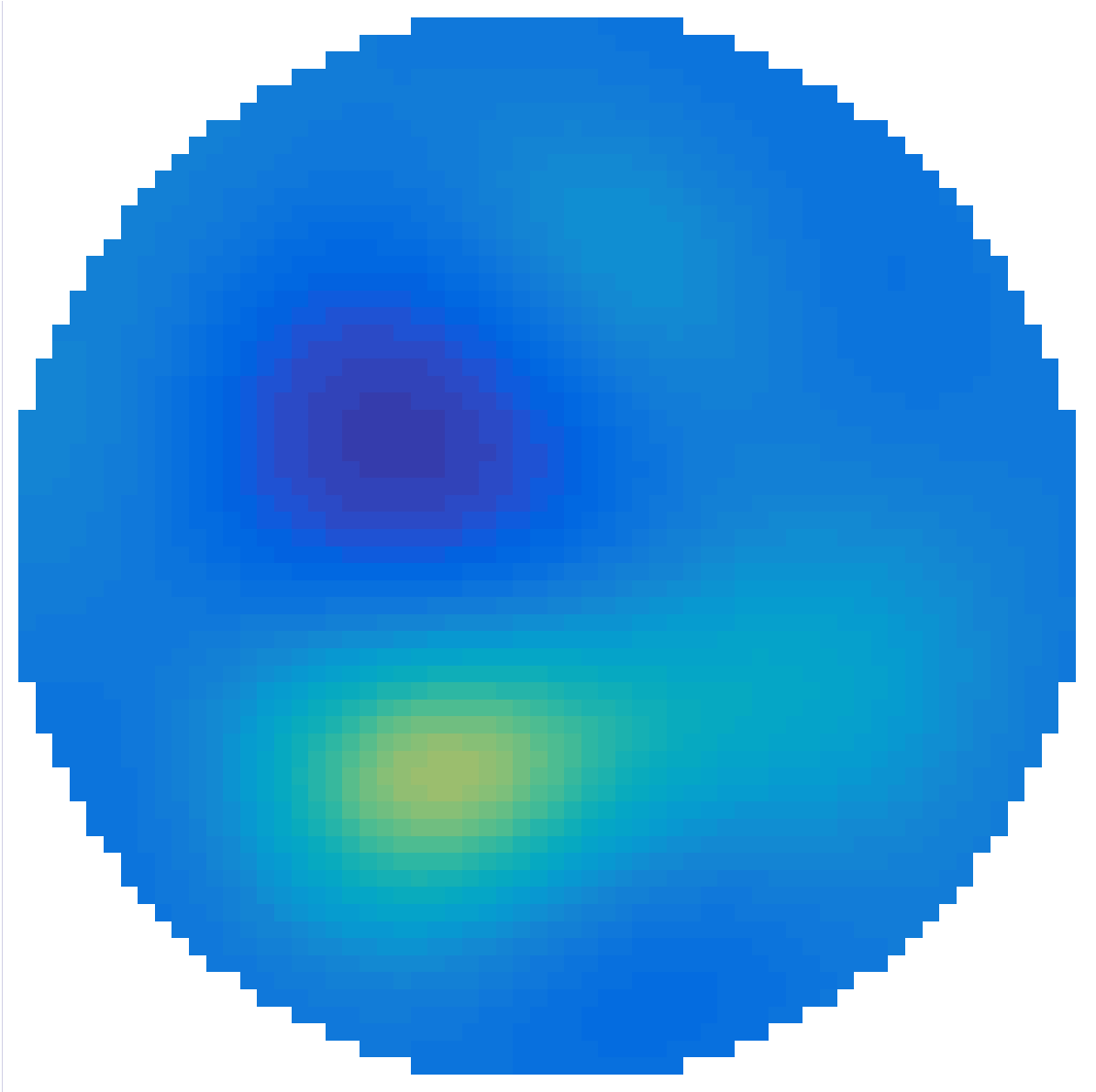}}
\put(160,170){\includegraphics[width=90pt]{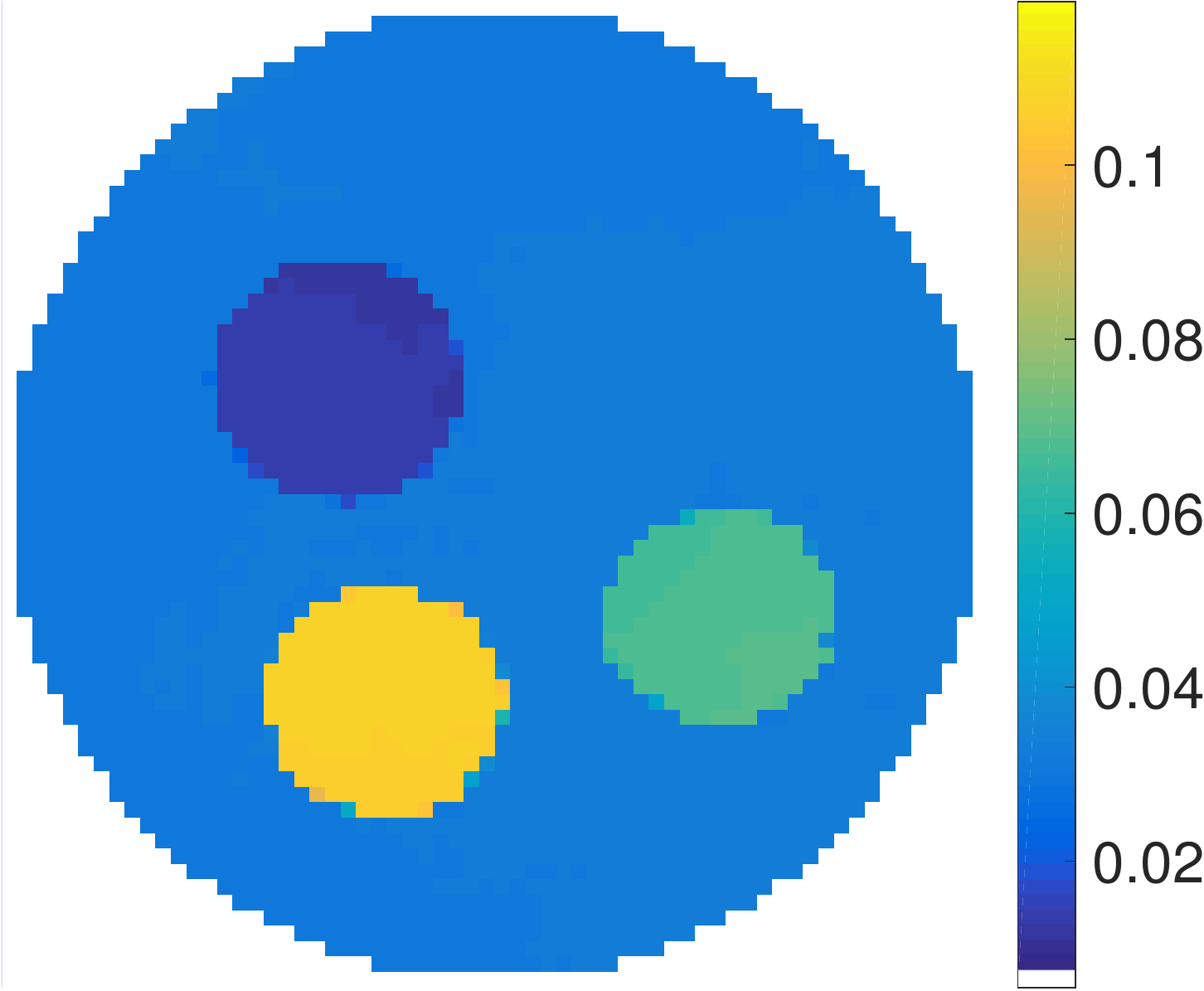}}

%
\put(10,85){\includegraphics[width=75pt]{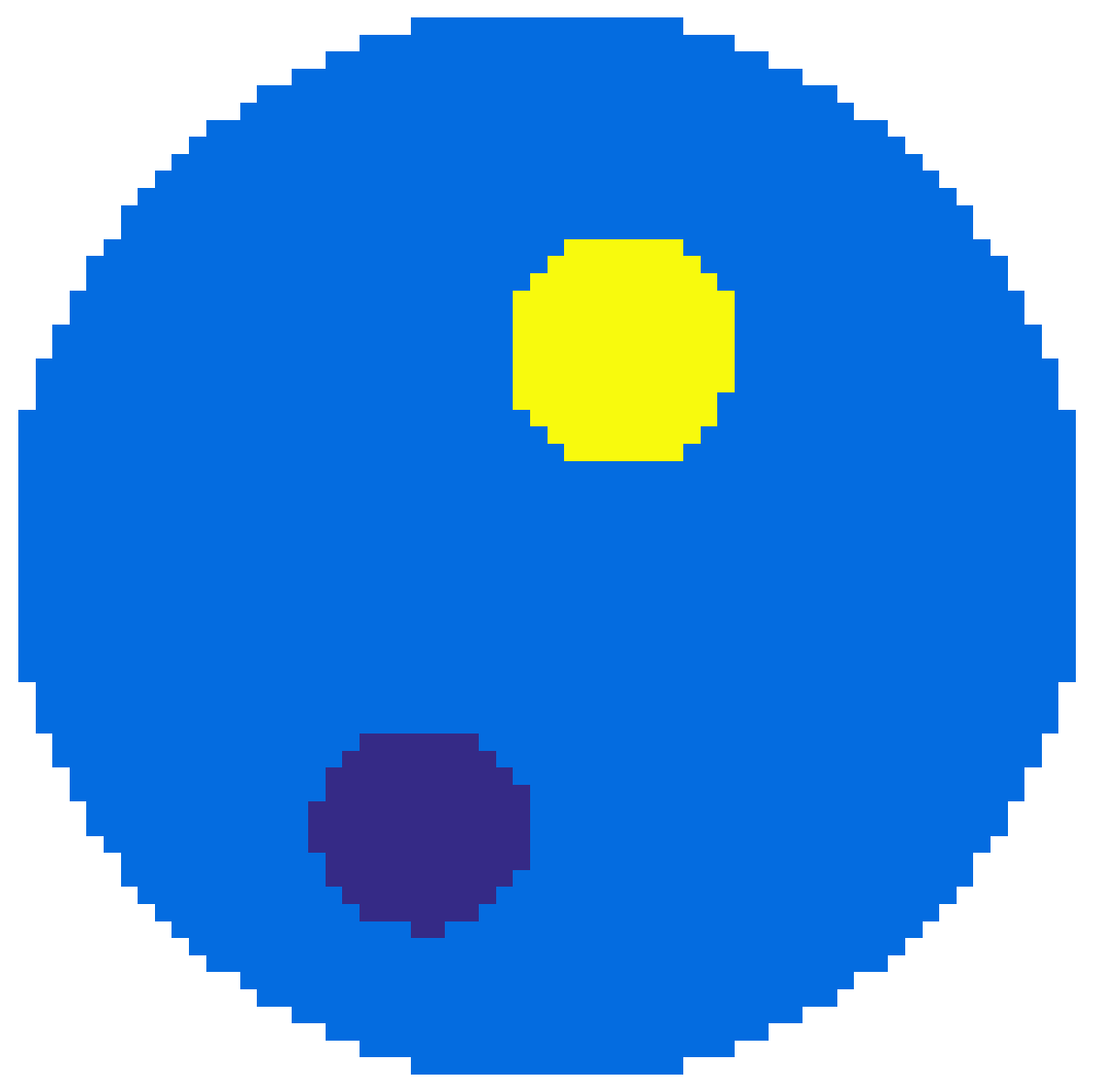}}
\put(85,85){\includegraphics[width=75pt]{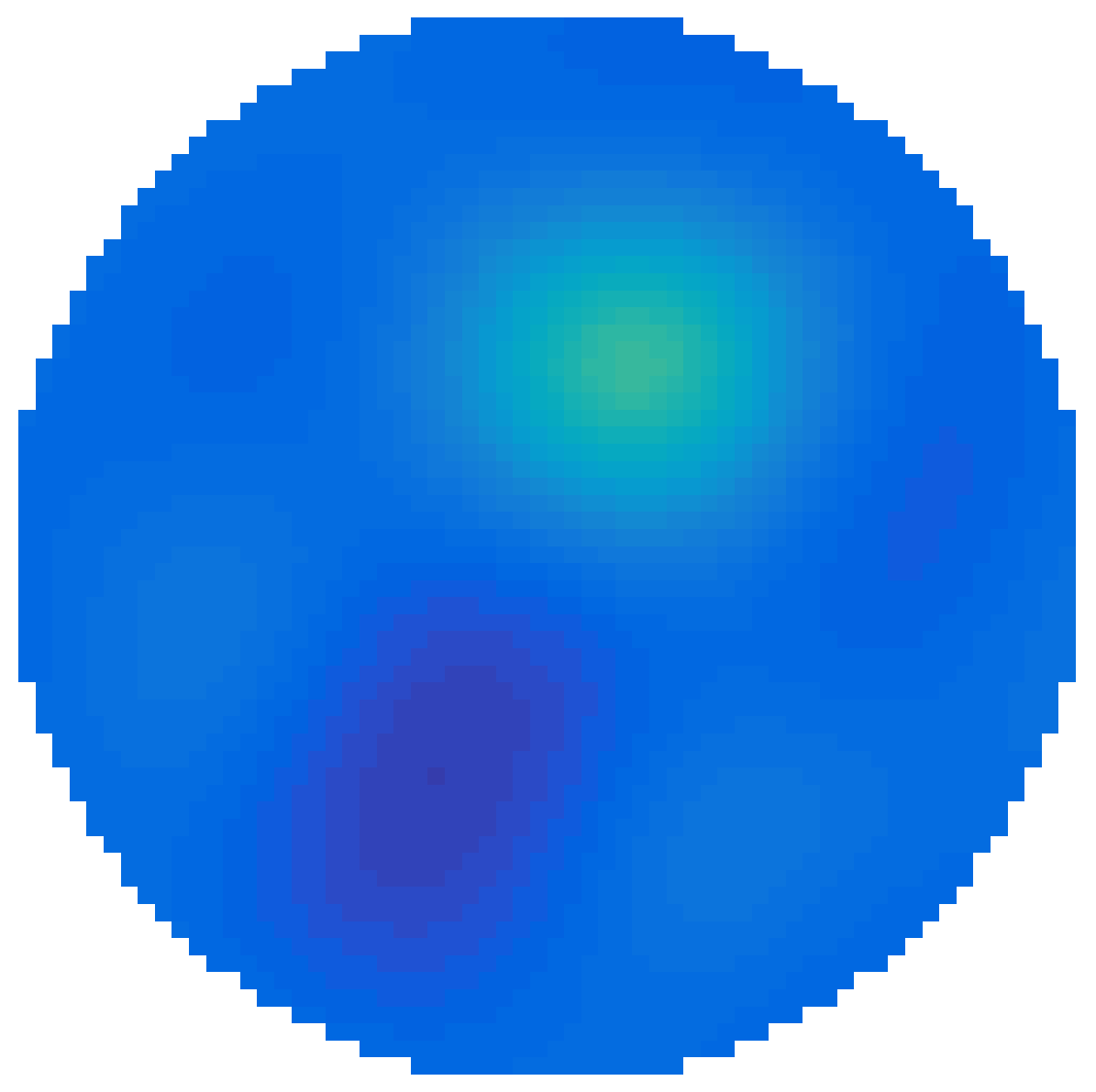}}
\put(160,85){\includegraphics[width=90pt]{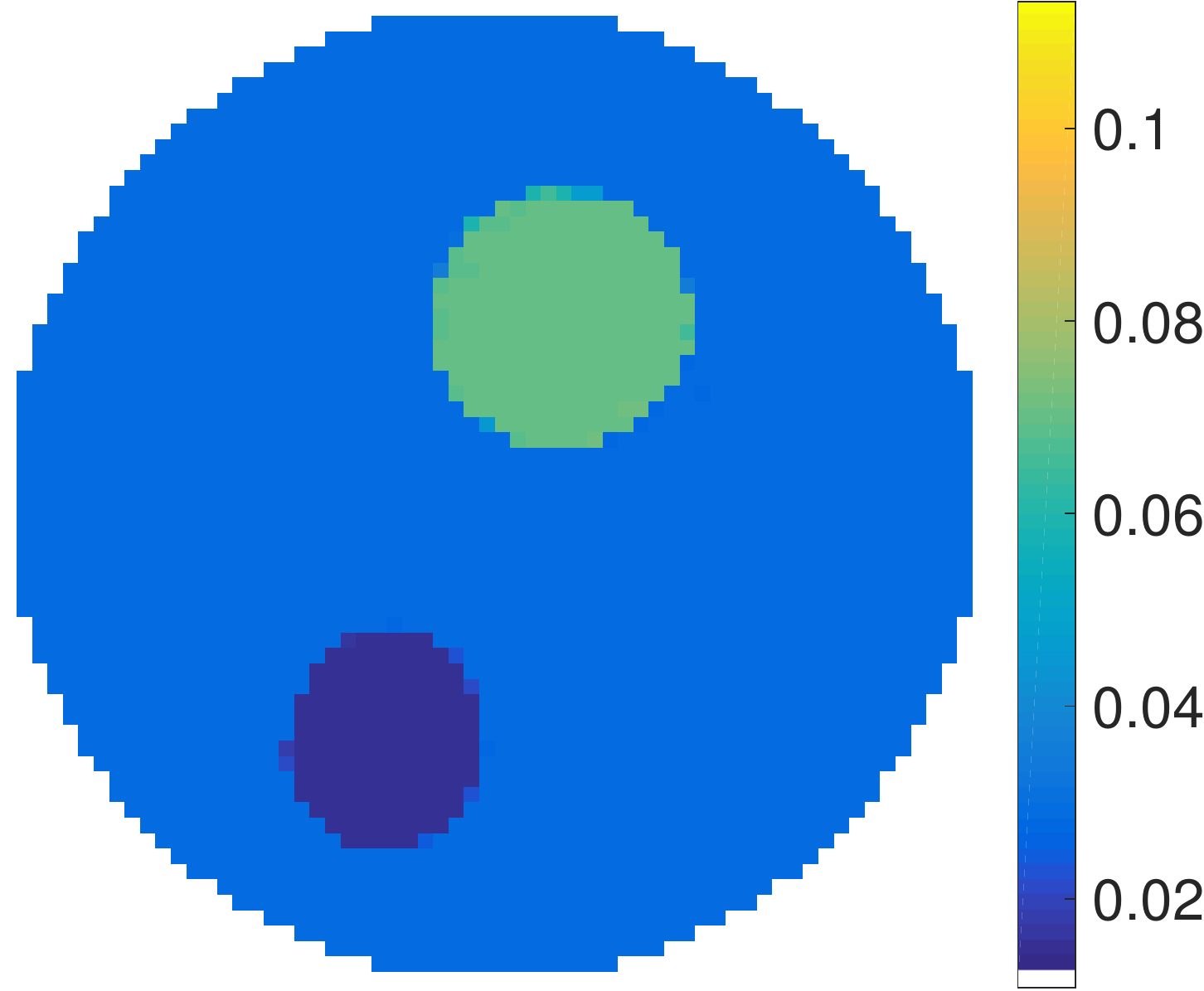}}
%
\put(10,0){\includegraphics[width=75pt]{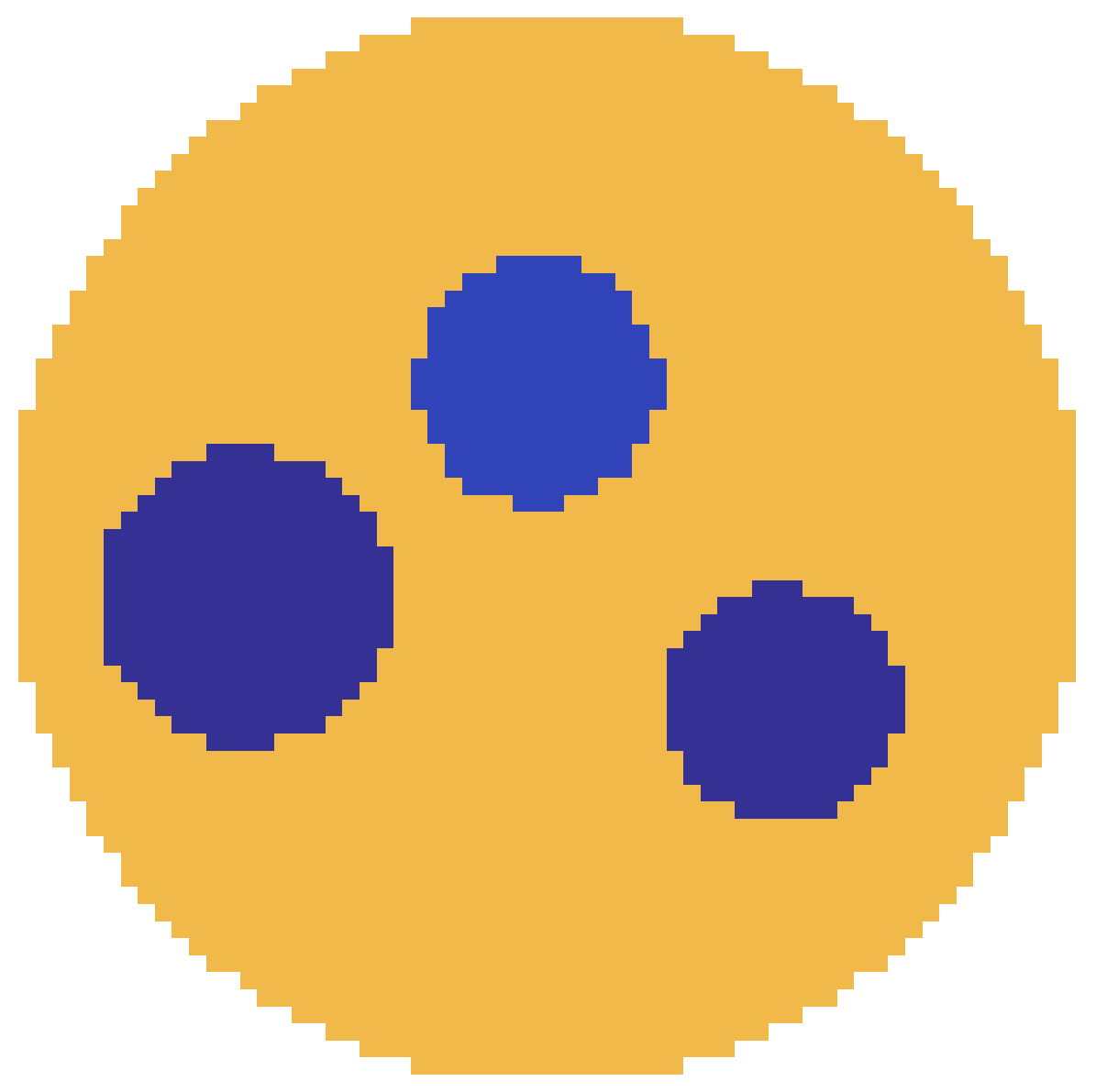}}
\put(85,0){\includegraphics[width=75pt]{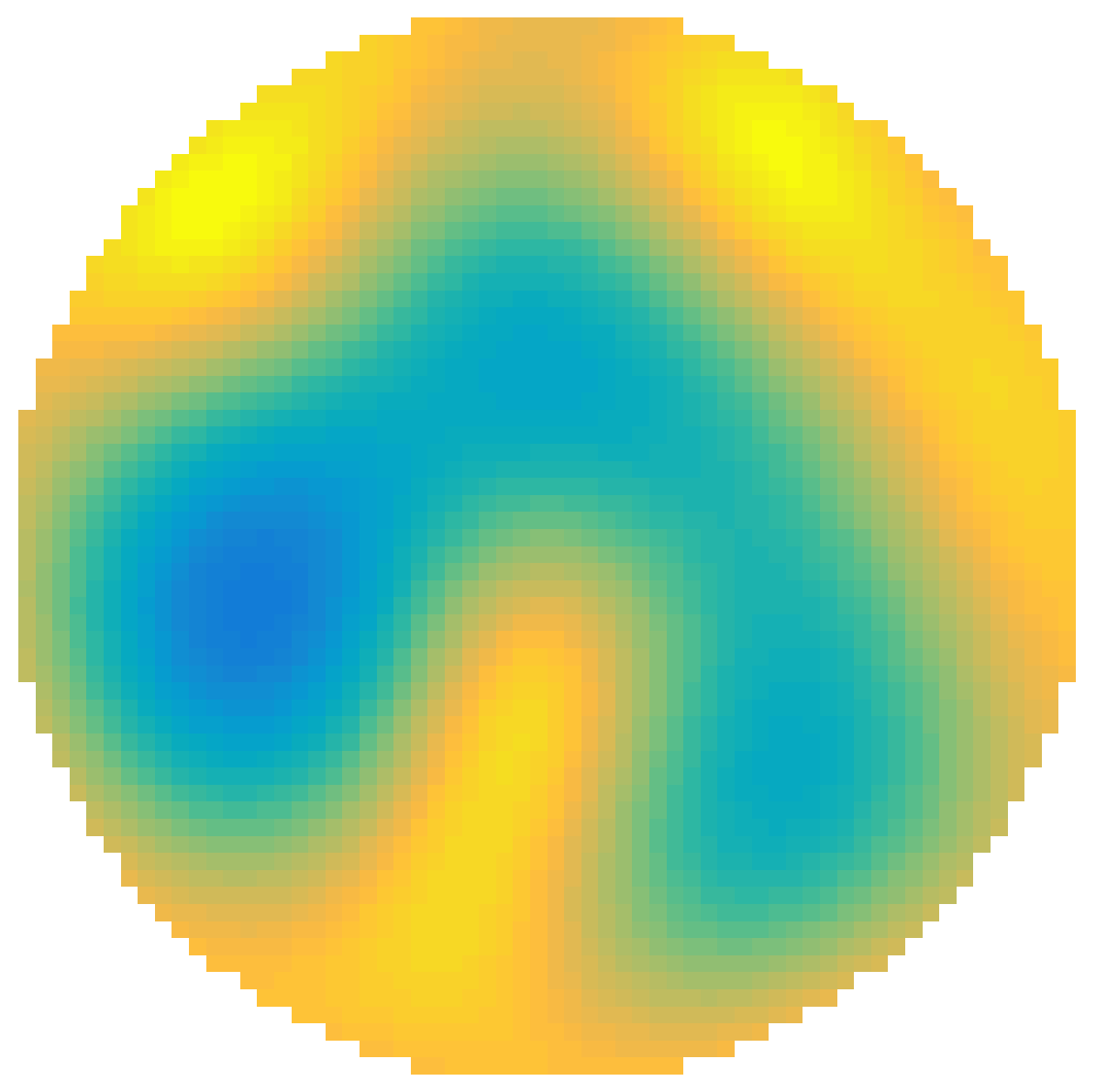}}
\put(160,0){\includegraphics[width=90pt]{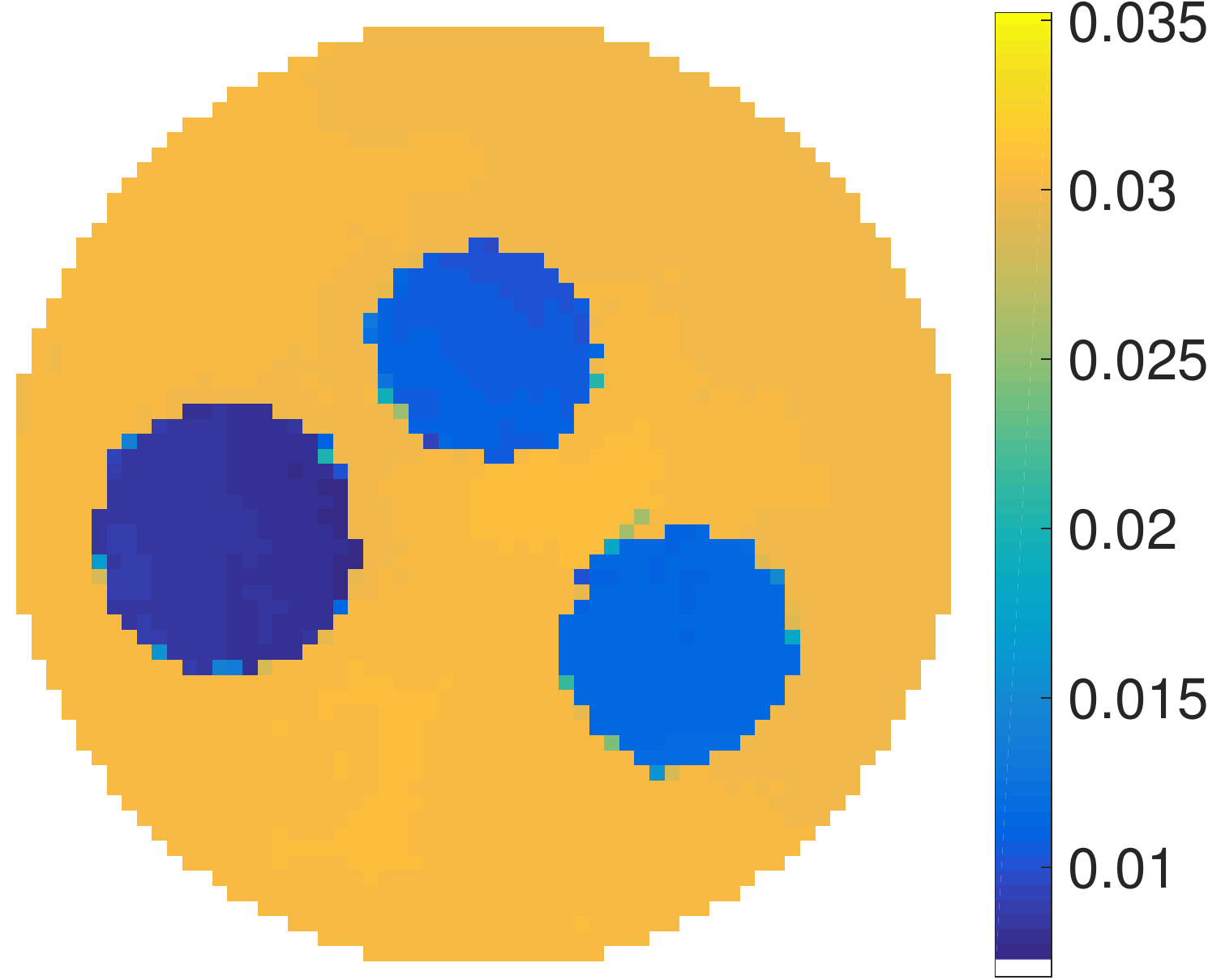}}

\put(25,255){\footnotesize \sc Phantom}

\put(108,265){\footnotesize \sc Low-pass}
\put(100,255){\footnotesize \sc D-bar Image}

\put(170,255){\footnotesize \sc  Deep D-bar Image}
\put(0,25){\rotatebox{90}{{\sc \footnotesize Sim 3}}}
\put(0,110){\rotatebox{90}{{\sc \footnotesize Sim 2}}}
\put(0,200){\rotatebox{90}{{\sc \footnotesize Sim 1}}}

\end{picture}
\caption{\label{fig:KIT4_Results_simData} Results for simulated test data from the KIT4 geometry. All phantom are drawn from the same distribution as the training data. The initial D-bar image is compared to the Deep D-bar image.  The D-bar images, on the full square are used as the `input' images for the CNN.  Images are displayed here on the circular geometry of the tank, for presentation only.  Each row is plotted on its own scale.}
\end{figure}
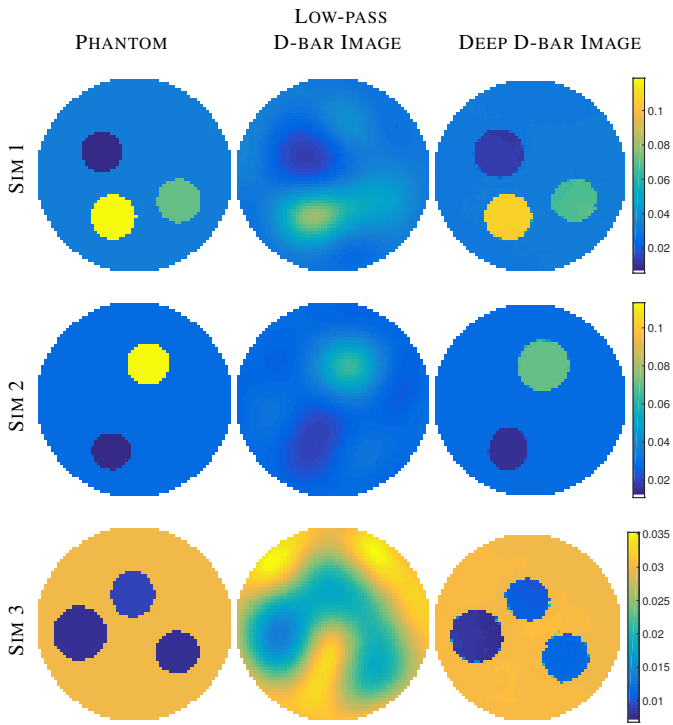

We begin with purely simulated data for the ACT4 and KIT4 examples.  For the ACT4 setting, we consider three scenarios, as shown in Figure~\ref{fig:ACT4_Results_simData}: one consistent with the training data but not used in the training (top), and two examples violating the training data - one with \trev{three} horizontally divided regions in the left lung (middle) and the final with a vertical division in the left lung (bottom).  Note that the `Low-pass D-bar' and `Deep D-bar' images are shown on the same scale for visual comparison.  The complete `input' and `output' of the CNN are on the unit square $[-1,1]^2$.  Reconstructions here are clipped to the disc for visualization purposes only.   Next, Figure~\ref{fig:KIT4_Results_simData} shows the results of the new algorithm on simulated data for the KIT4 example.  Three scenarios consistent with the training data, but not used in the training, are presented.  

\begin{figure}[h!]
\begin{picture}(100,190)
\put(40,10){\includegraphics[height=160pt]{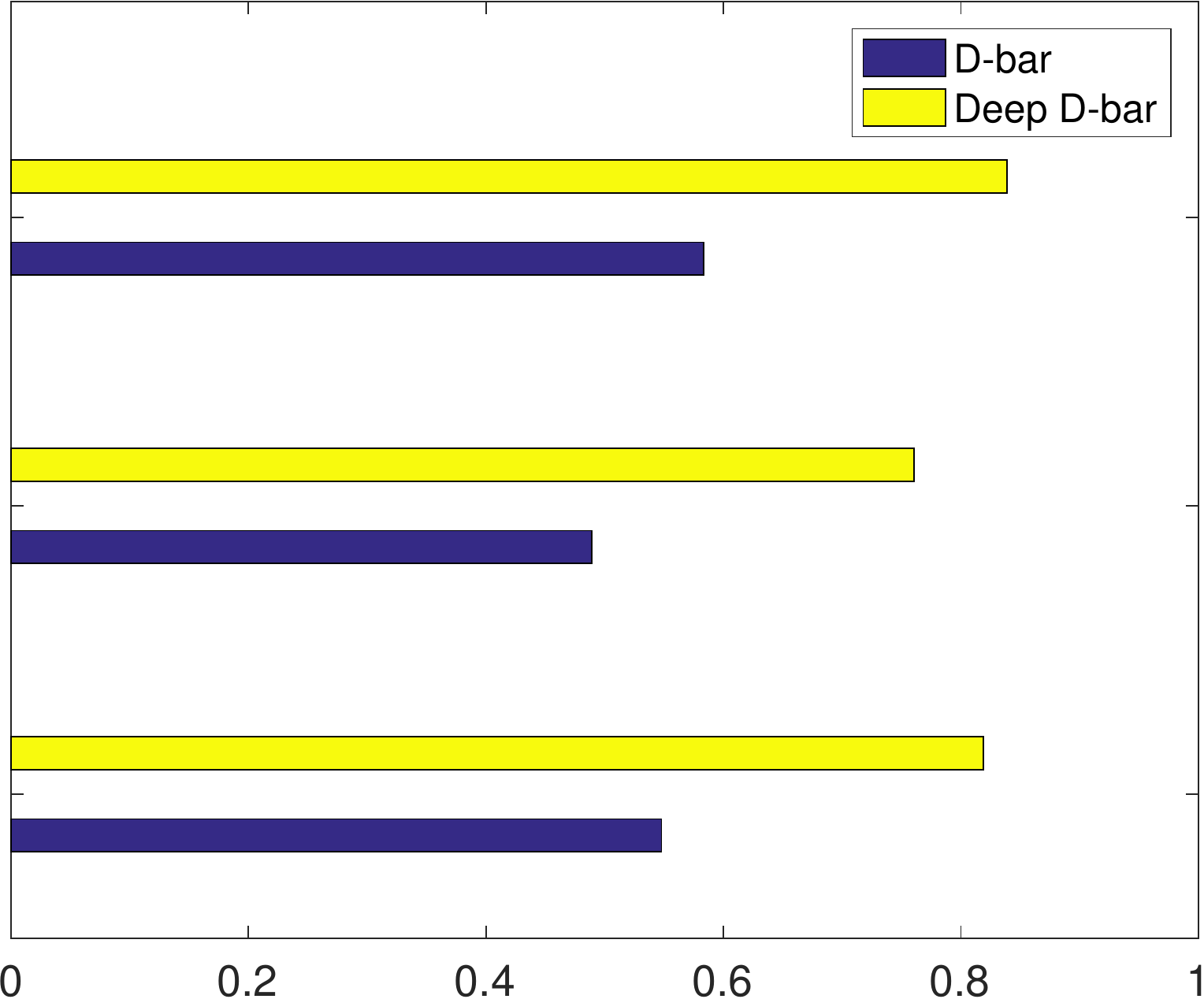}}

\put(130,0){\sc \footnotesize SSIM}
\put(50,180){\sc \small SSIM Comparison for ACT4 Sim. Recons }

\put(30,30){\rotatebox{90}{{\sc \scriptsize`Sim 3'}}}
\put(30,80){\rotatebox{90}{{\sc \scriptsize`Sim 2'}}}
\put(30,125){\rotatebox{90}{{\sc \scriptsize`Sim 1'}}}

\end{picture}
\caption{SSIM measurements are compared for the D-bar method and the new `Deep D-bar' method for the ACT4 reconstructions for the simulated data shown in Figure~\ref{fig:ACT4_Results_simData}.}  
\label{fig:act4_SSIMplot_simData}
\end{figure}
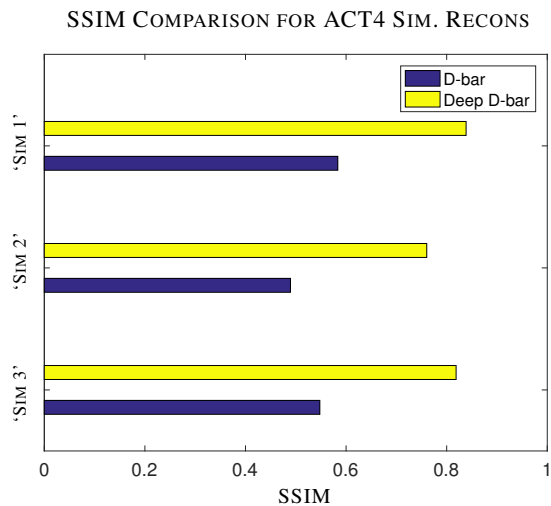

\begin{figure}[h!]
\begin{picture}(100,190)
\put(40,10){\includegraphics[height=160pt]{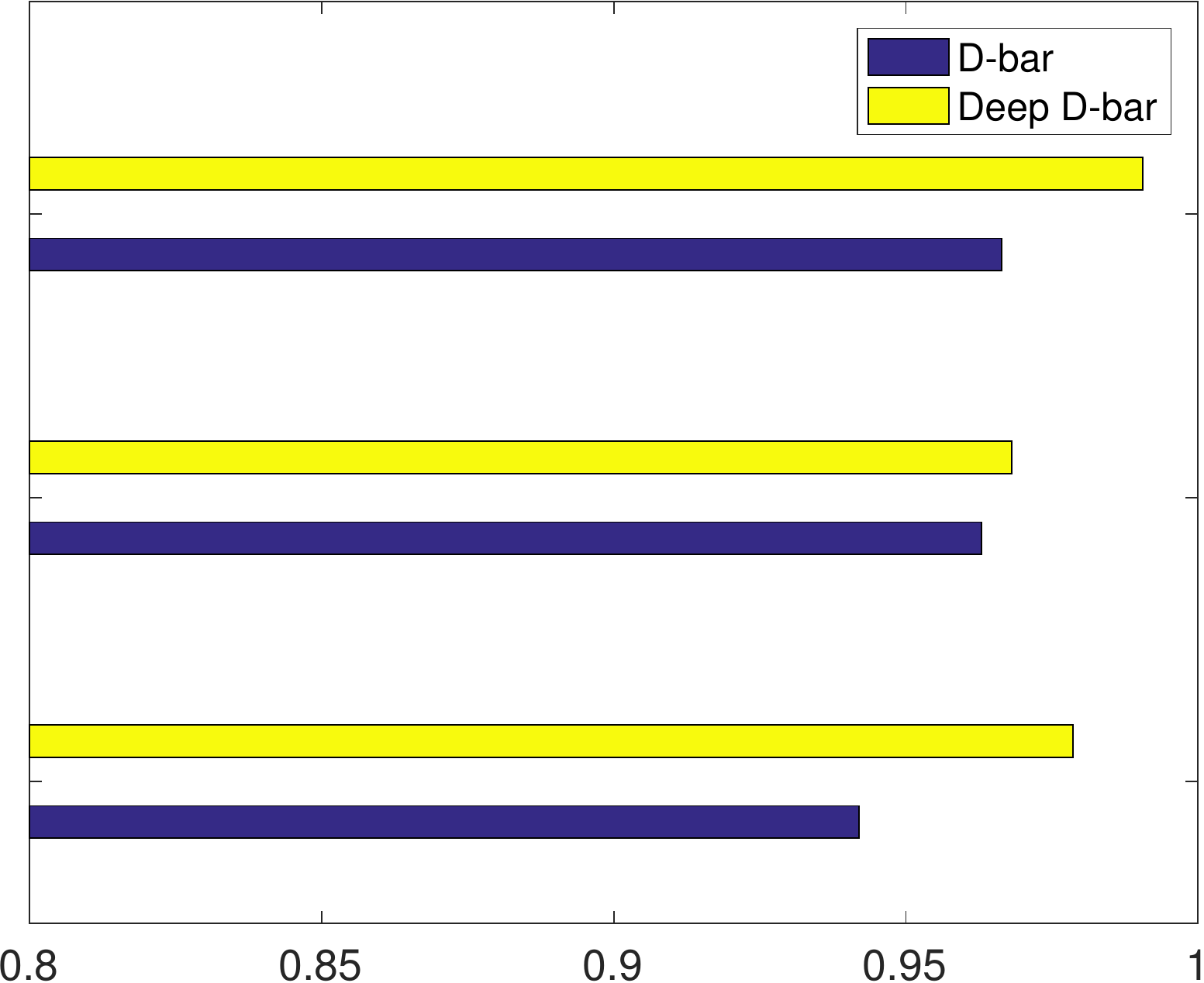}}

\put(130,0){\sc \footnotesize SSIM}
\put(50,180){\sc \small SSIM Comparison for KIT4 Sim. Recons }

\put(30,30){\rotatebox{90}{{\sc \scriptsize`Sim 3'}}}
\put(30,80){\rotatebox{90}{{\sc \scriptsize`Sim 2'}}}
\put(30,125){\rotatebox{90}{{\sc \scriptsize`Sim 1'}}}

\end{picture}
\caption{SSIM measurements are compared for the D-bar method and the new `Deep D-bar' method for the KIT4 reconstructions for the simulated data shown in Figure~\ref{fig:KIT4_Results_simData}.}  
\label{fig:kit4_SSIMplot_simData}
\end{figure}
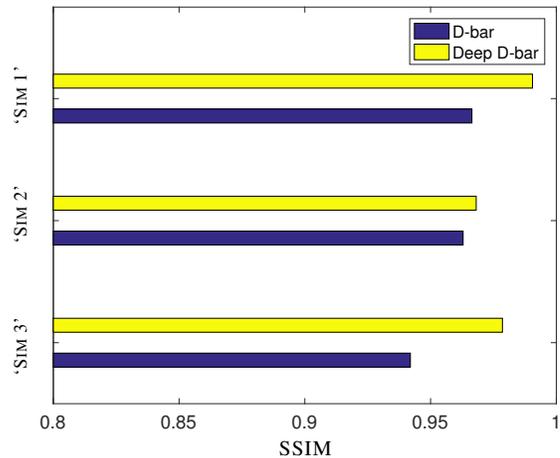

{\it Structural Similarity Indices} (SSIMs) computed for the \trev{simulated} ACT4 and KIT4 examples \trev{are} shown in Figures~\ref{fig:act4_SSIMplot_simData} and \ref{fig:kit4_SSIMplot_simData}, respectively. Additionally, we evaluated the  minimized $\ell^2$-loss by computing the mean relative error for a test set of 16 samples drawn from the same distribution as the training data.  The For the ACT4 simulations we improved from $28.05\%$ to $9.92\%$ and for the KIT4 test data from $16.82\%$ to $9.12\%$ relative $\ell^2$-error.

\subsection{Reconstructions from Experimental Data}
\begin{figure}[h!]
\centering
\begin{picture}(300,370)

\put(10,255){\includegraphics[width=75pt]{HLSA_saline_DeepDbar_DICOM_smaller.jpg}}
\put(85,255){\includegraphics[width=75pt]{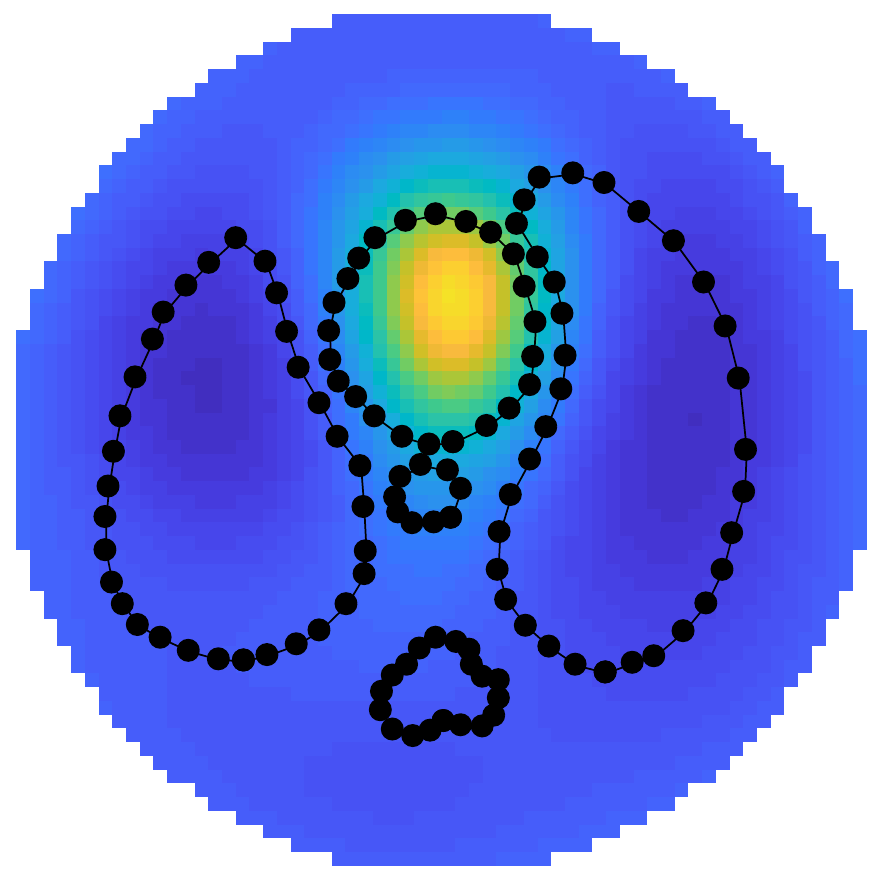}}
\put(160,255){\includegraphics[width=90pt]{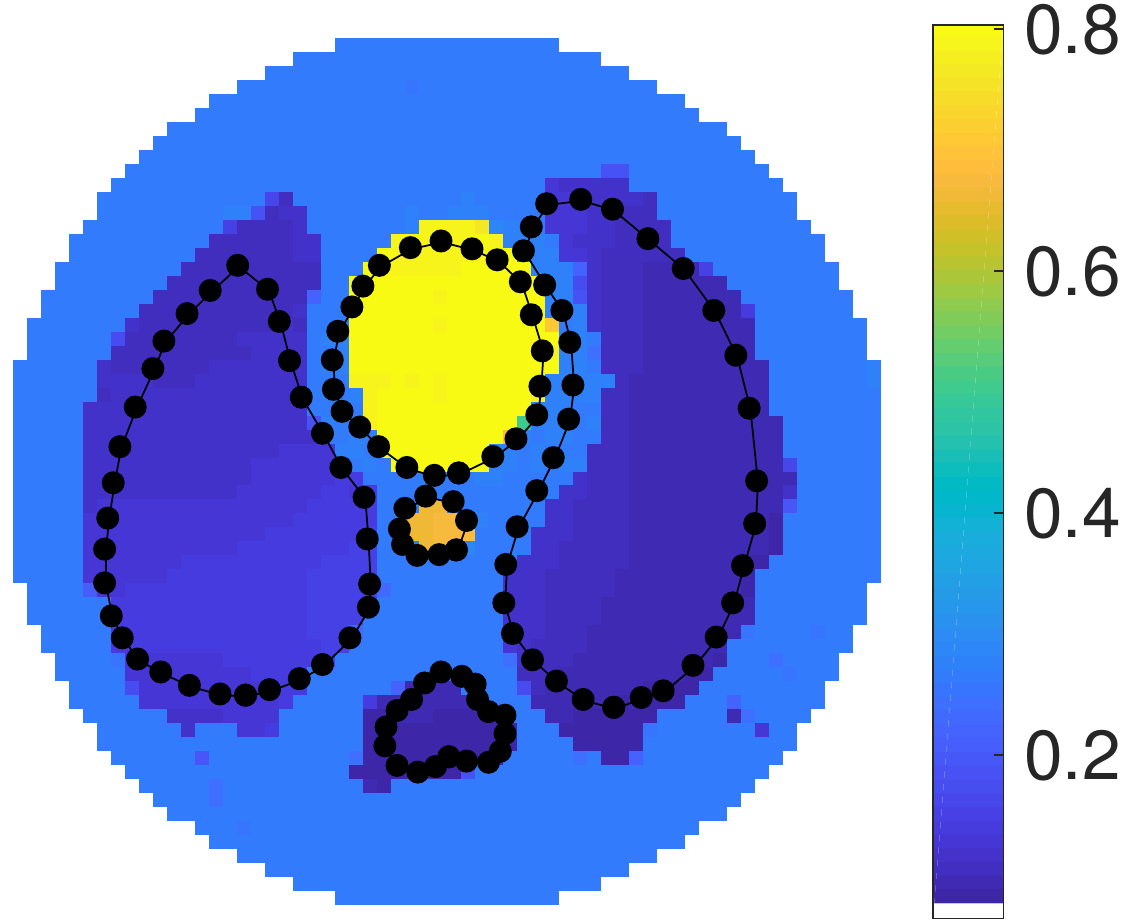}}

\put(10,170){\includegraphics[width=75pt]{HLSA_plueral_saline_DeepDbar_DICOM_smaller.jpg}}
\put(85,170){\includegraphics[width=75pt]{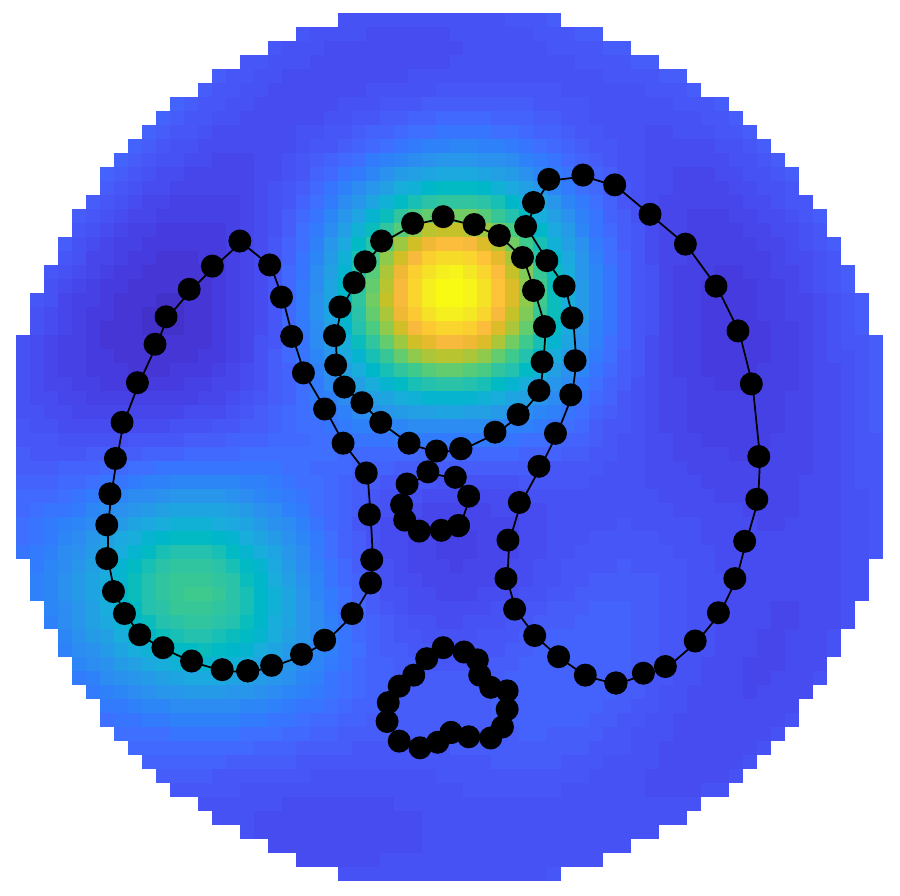}}
\put(160,170){\includegraphics[width=90pt]{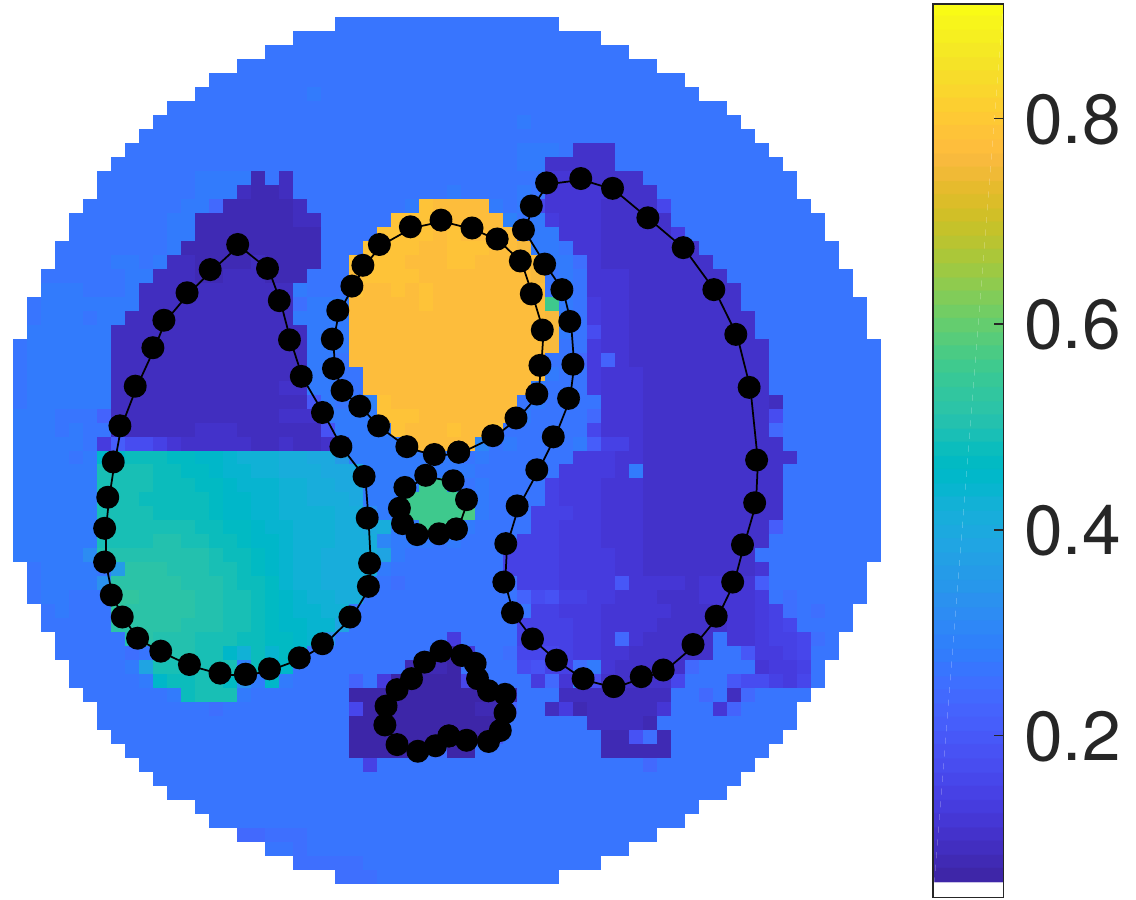}}

\put(10,85){\includegraphics[width=75pt]{HLSA_plastic_saline_DeepDbar_DICOM_smaller.jpg}}
\put(85,85){\includegraphics[width=75pt]{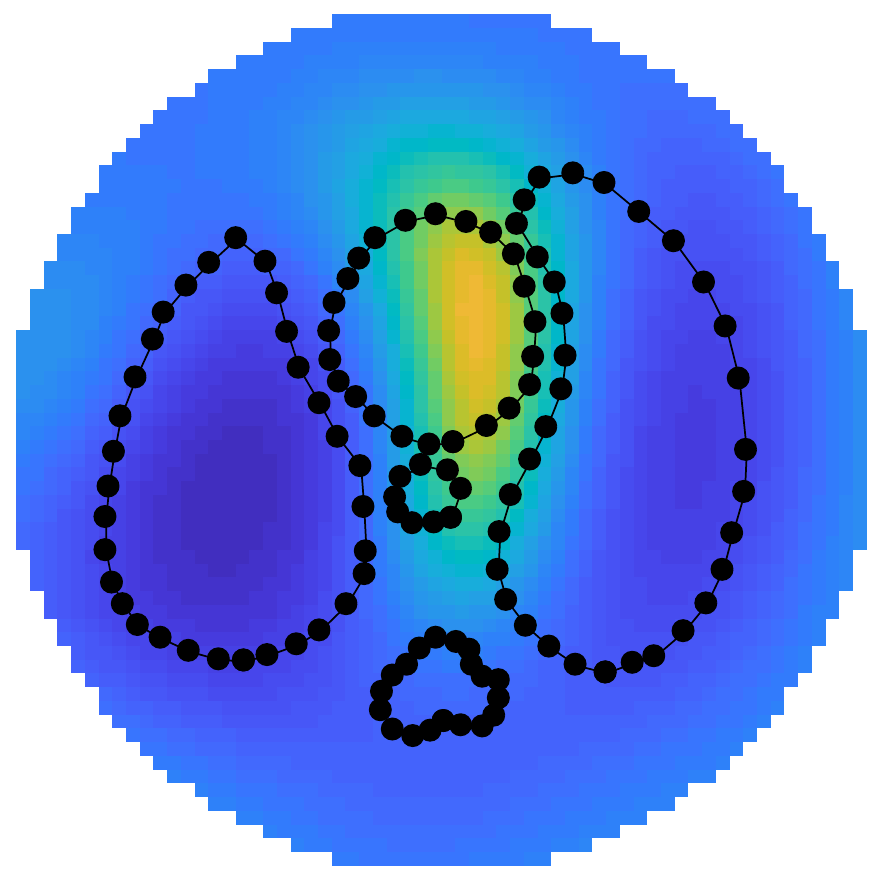}}
\put(160,85){\includegraphics[width=90pt]{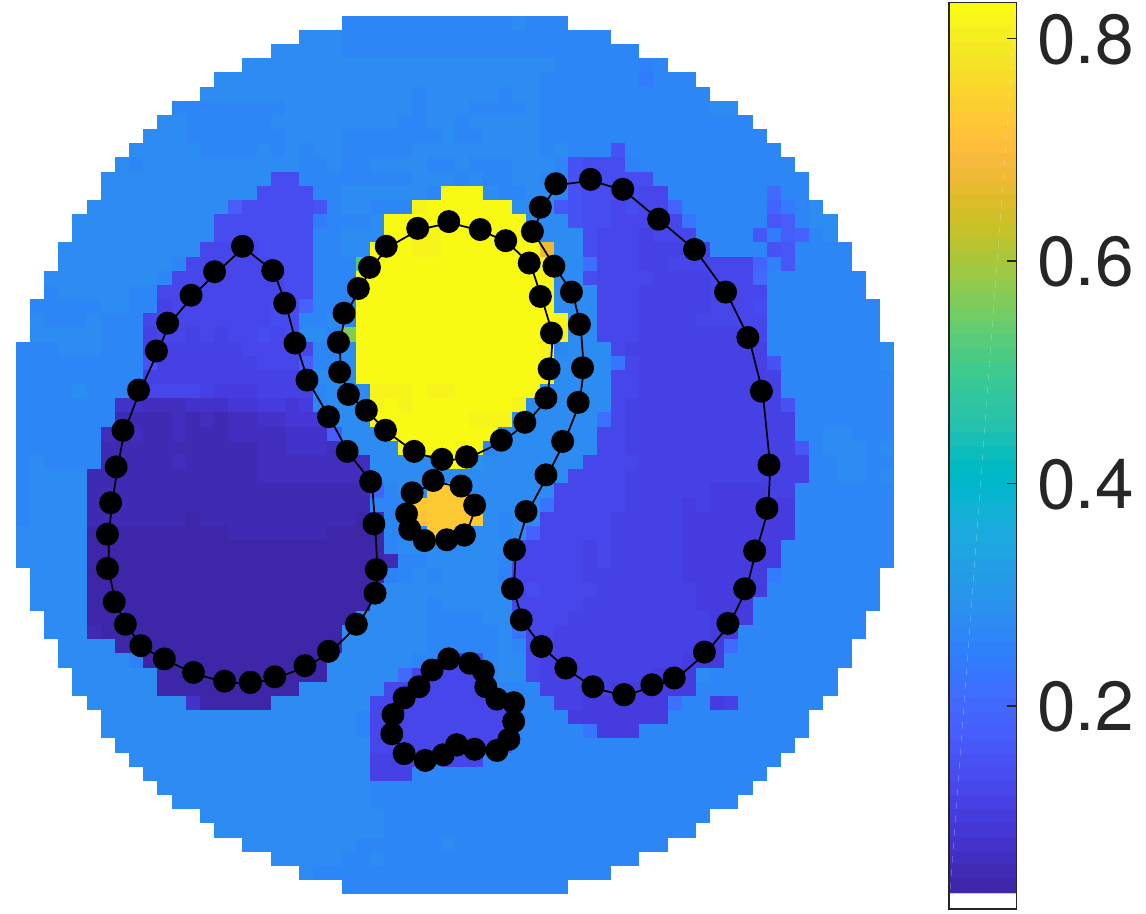}}

\put(10,0){\includegraphics[width=75pt]{ACT4_HLSA_wCopper_Correct_DICOM_smaller.jpg}}
\put(85,0){\includegraphics[width=75pt]{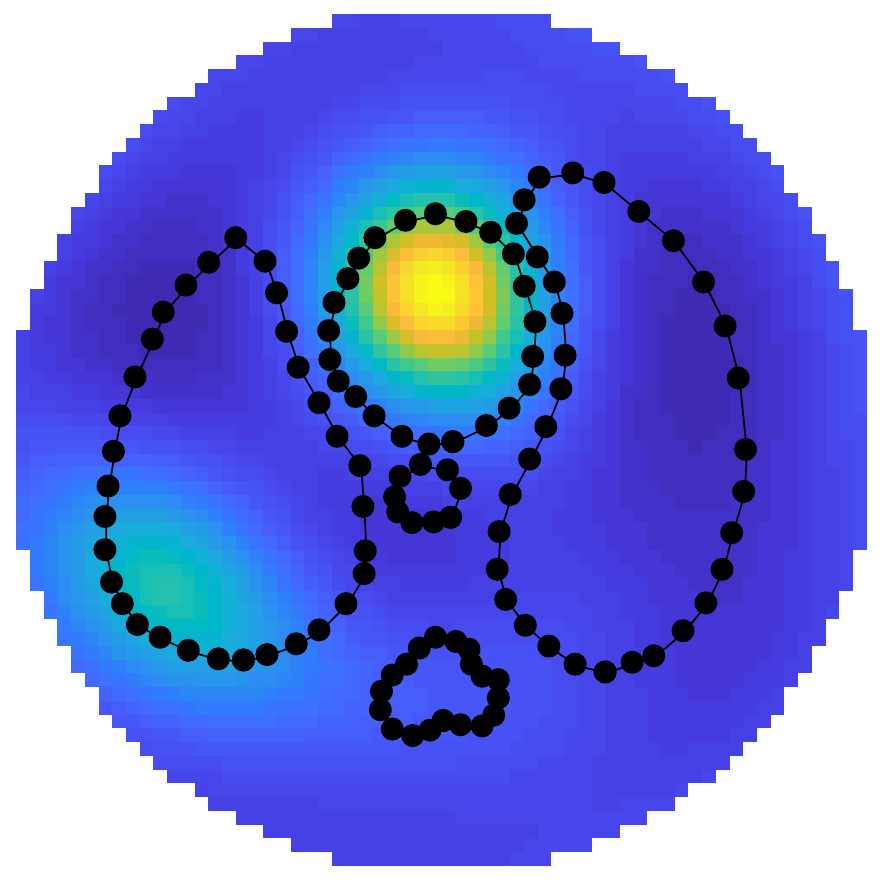}}
\put(160,0){\includegraphics[width=90pt]{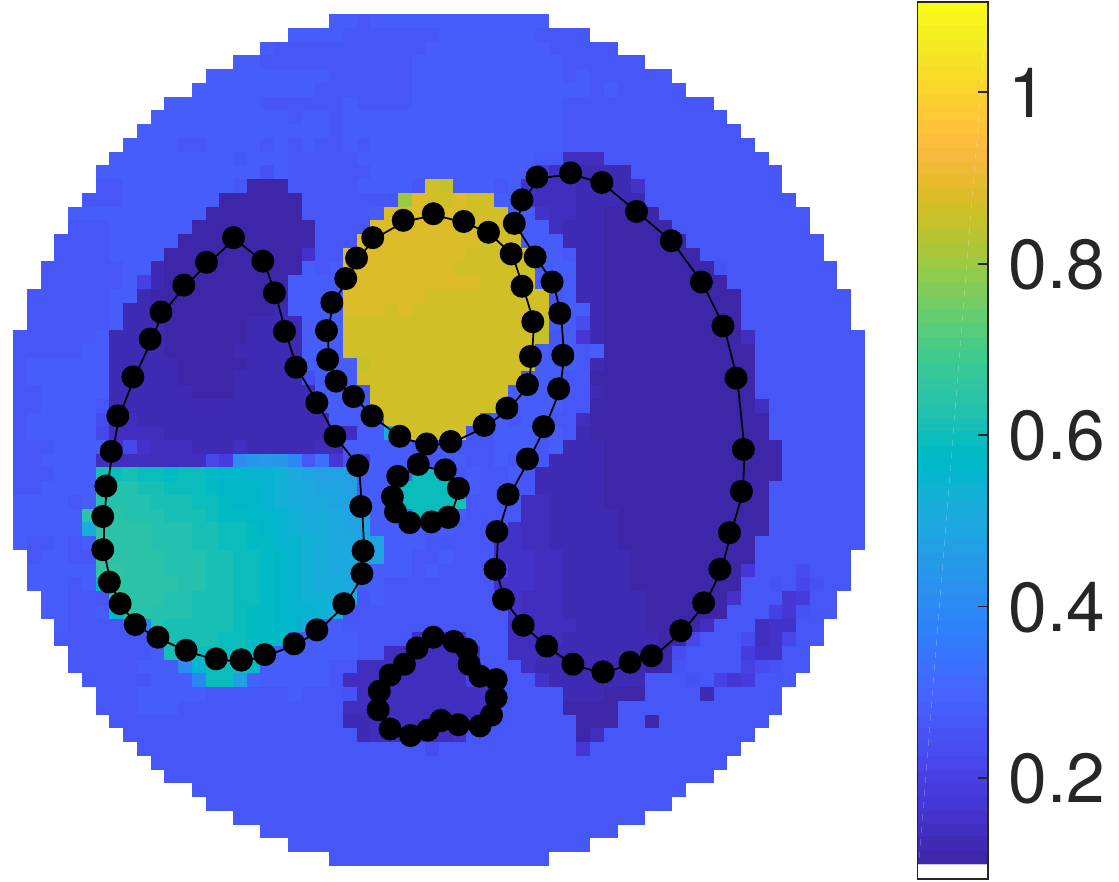}}

\put(25,345){\footnotesize \sc Experiment}

\put(108,355){\footnotesize \sc Low Pass}
\put(100,345){\footnotesize \sc D-bar Image}

\put(170,345){\footnotesize \sc  Deep D-bar Image}

\put(0,20){\rotatebox{90}{{\sc \footnotesize Injury 3}}}
\put(0,105){\rotatebox{90}{{\sc \footnotesize Injury 2}}}
\put(0,190){\rotatebox{90}{{\sc \footnotesize Injury 1}}}
\put(0,280){\rotatebox{90}{{\sc \footnotesize Healthy}}}

\end{picture}
\caption{\label{fig:ACT4_Results} ACT4 Results for the various test scenarios: Healthy, Injuries 1-3 corresponding to conductive agar, plastic tubes, and conductive copper tubes, respectively.  The initial D-bar image is compared to the Deep D-bar image.  The D-bar images, on the full square are used as the `input' images for the CNN.  Images are displayed here on the circular geometry of the tank, for presentation only.  Each row is plotted on its own scale.}
\end{figure}

\begin{figure}[h!]
\begin{picture}(180,210)

\put(20,10){\includegraphics[height=180pt]{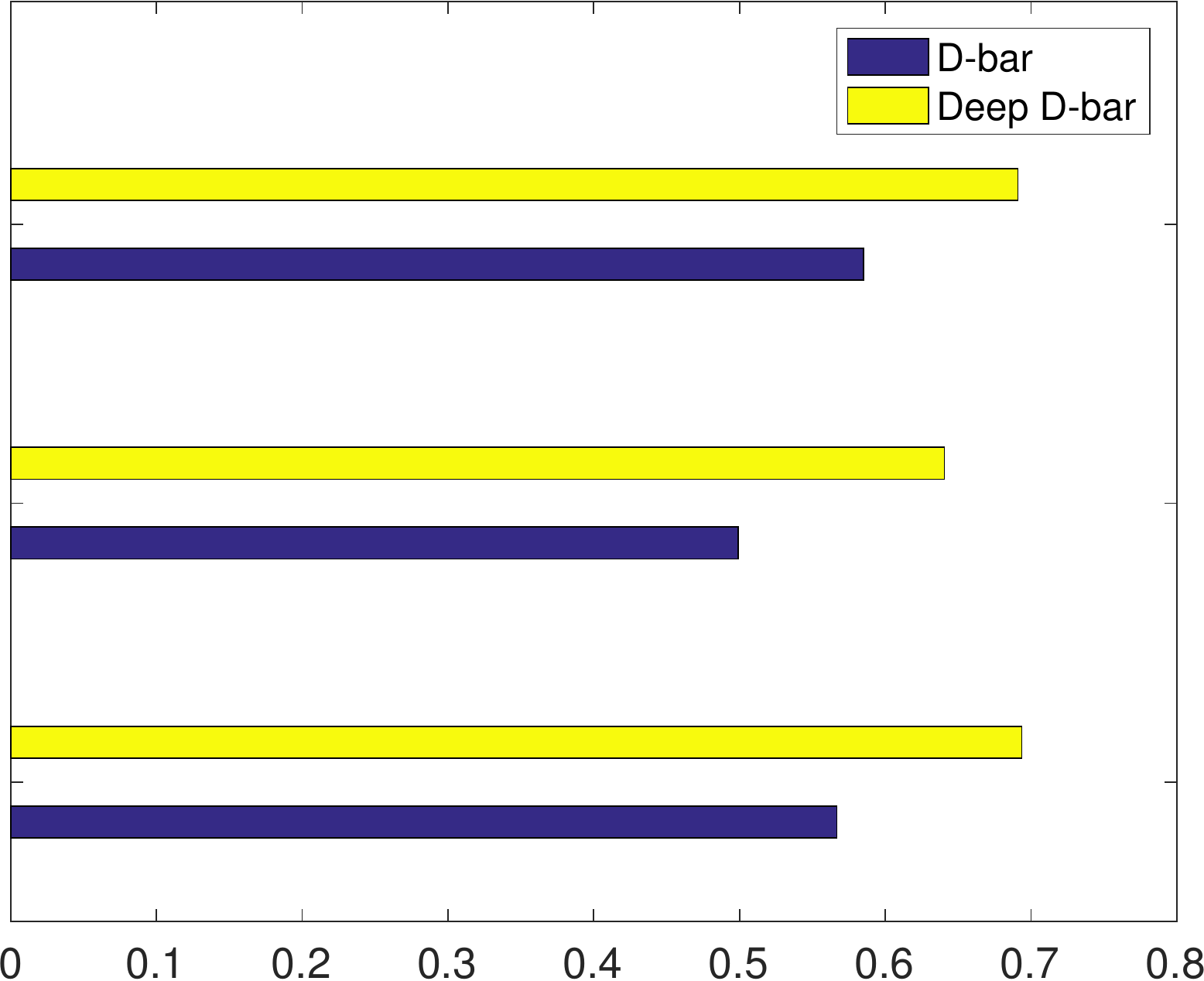}}

\put(120,0){\sc \footnotesize SSIM}
\put(50,200){\sc \small SSIM Comparison for ACT4 Recons}

\put(10,30){\rotatebox{90}{{\sc \scriptsize Injury 2}}}
\put(10,80){\rotatebox{90}{{\sc \scriptsize Injury 1}}}
\put(10,130){\rotatebox{90}{{\sc \scriptsize Healthy}}}

\end{picture}
\caption{SSIM measurements are compared for the D-bar method and the new `Deep D-bar' method for the ACT4 experimental data reconstructions shown in Figure~\ref{fig:ACT4_Results}.  Note that meaningful SSIMs could not be computed for `Injury 3' due to the copper inclusions which have infinite conductivity.}  
\label{fig:act4_SSIMplot}
\end{figure}

\trev{Next,} we proceed to reconstructions from experimental data.  Figure~\ref{fig:ACT4_Results} depicts the results of the {\it Deep D-bar} approach on four experiments with ACT4 data: {\sc Healthy} and {\sc Injuries 1-3} as shown in Figure~\ref{fig:RPI_phantoms}.  The black dots represent the approximate boundaries of the `healthy' organs, extracted from the photograph.  SSIMs (Figure~\ref{fig:act4_SSIMplot}) were computed for the experimental reconstructions with the exception of {\sc Injury~3}, which has the infinite conductors (copper tubes).  The SSIM comparisons used approximate `truth' images formed by assigning the measured conductivity values (Table~\ref{table:ACT4_setup}) in the respective regions.  

\begin{figure}[h!]
\centering
\begin{picture}(300,370)

\put(10,255){\includegraphics[width=75pt]{fantom_1_1.jpg}}
\put(85,255){\includegraphics[width=75pt]{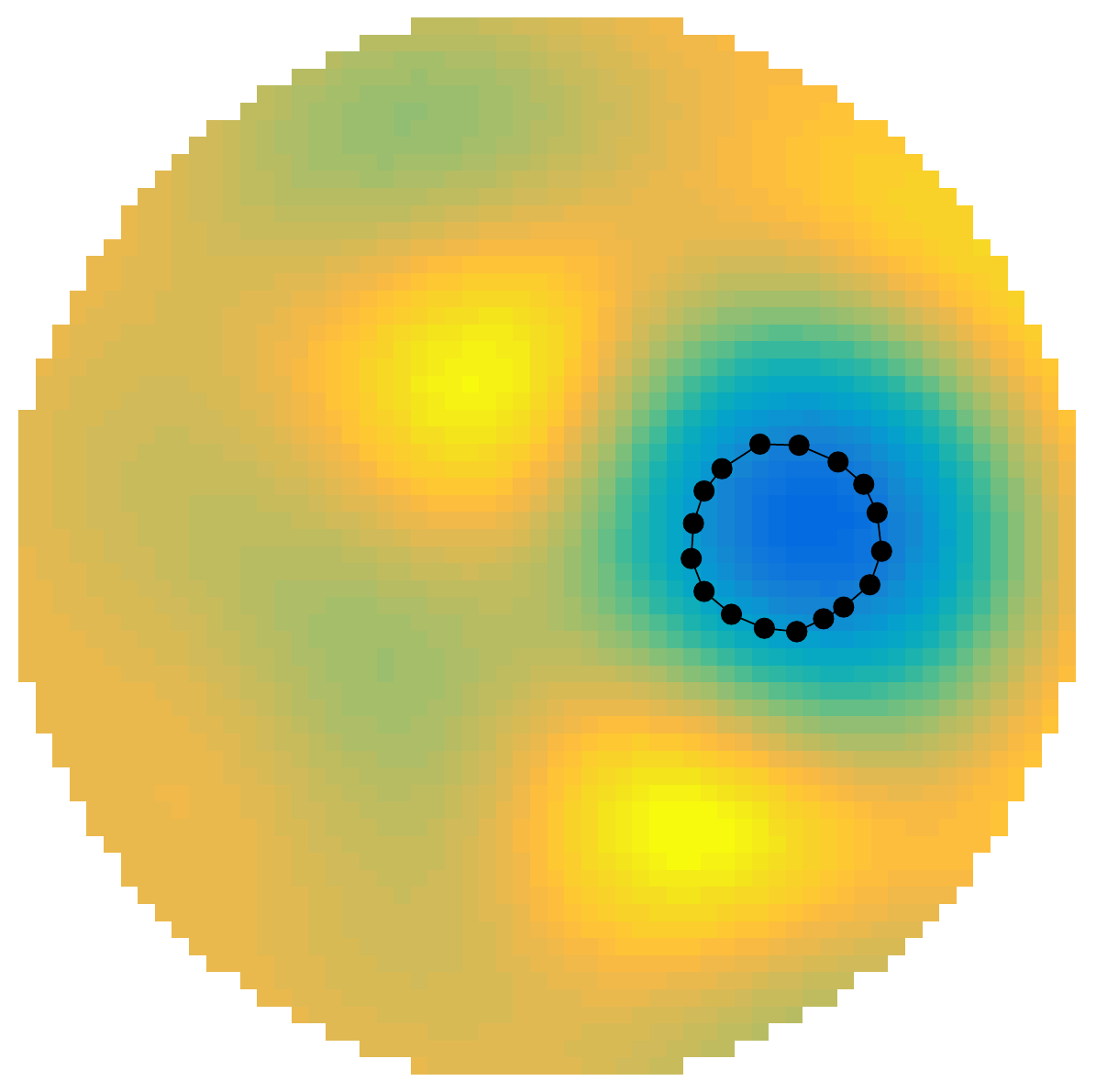}}
\put(160,255){\includegraphics[width=90pt]{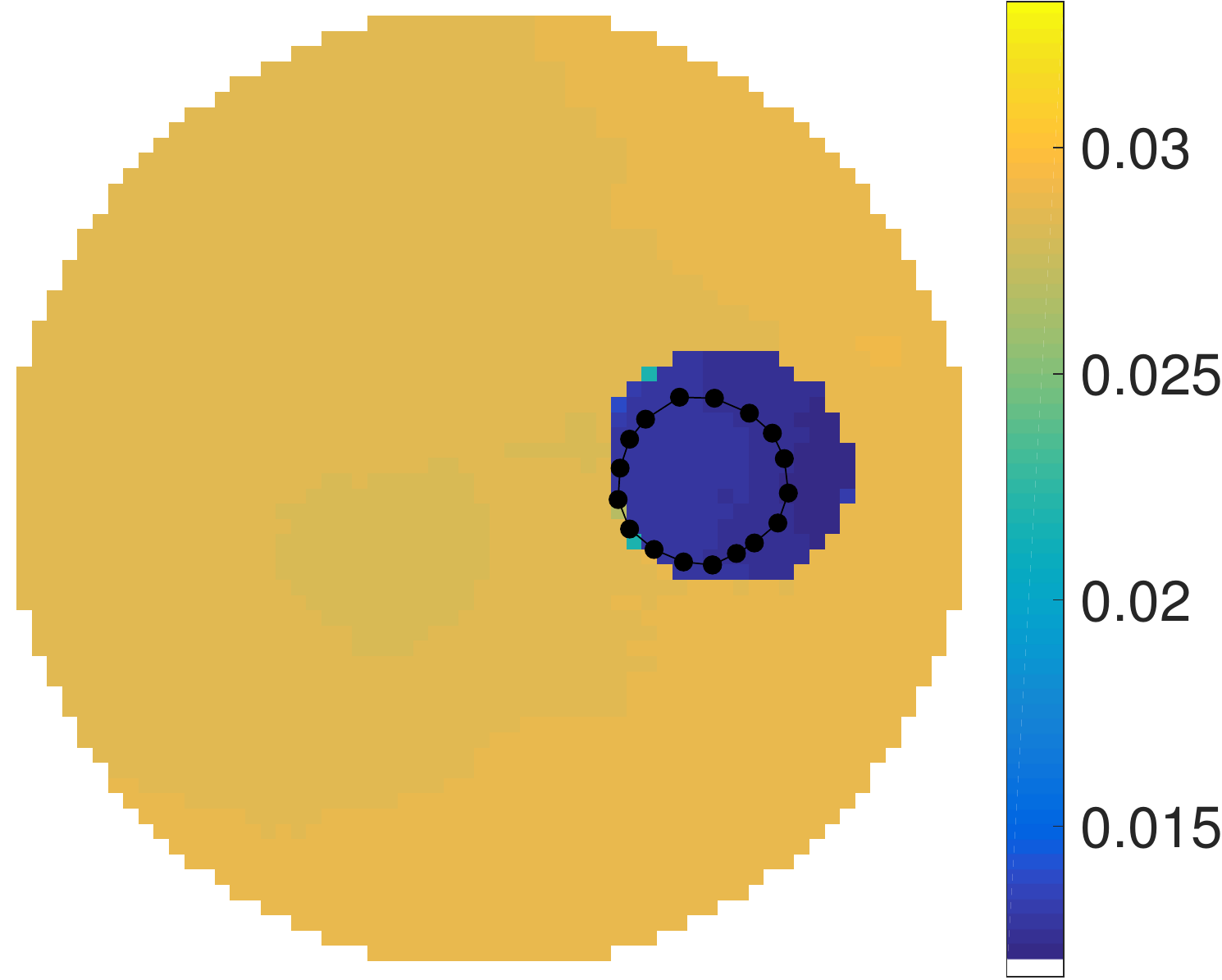}}

\put(10,170){\includegraphics[width=75pt]{fantom_2_2.jpg}}
\put(85,170){\includegraphics[width=75pt]{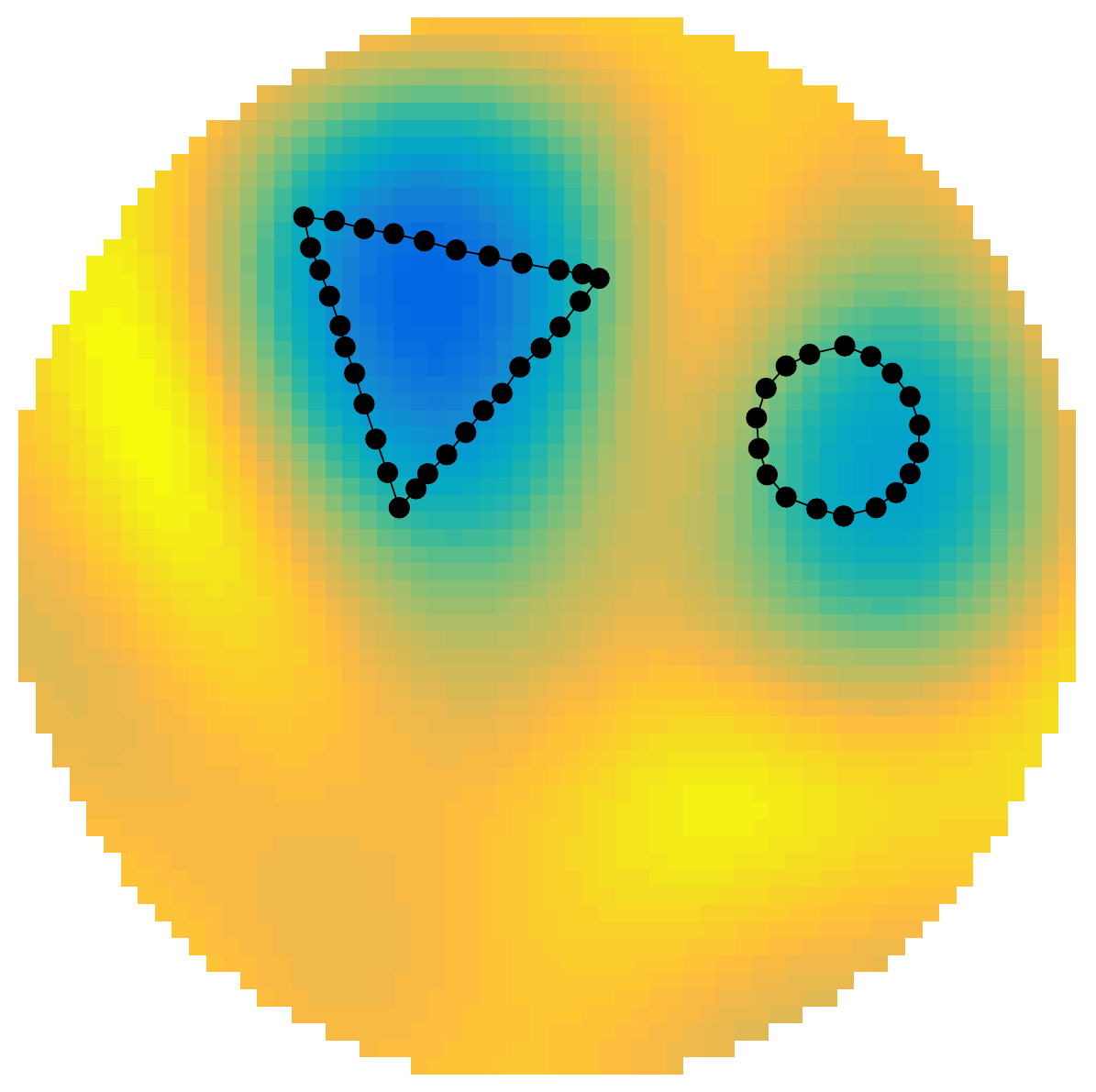}}
\put(160,170){\includegraphics[width=90pt]{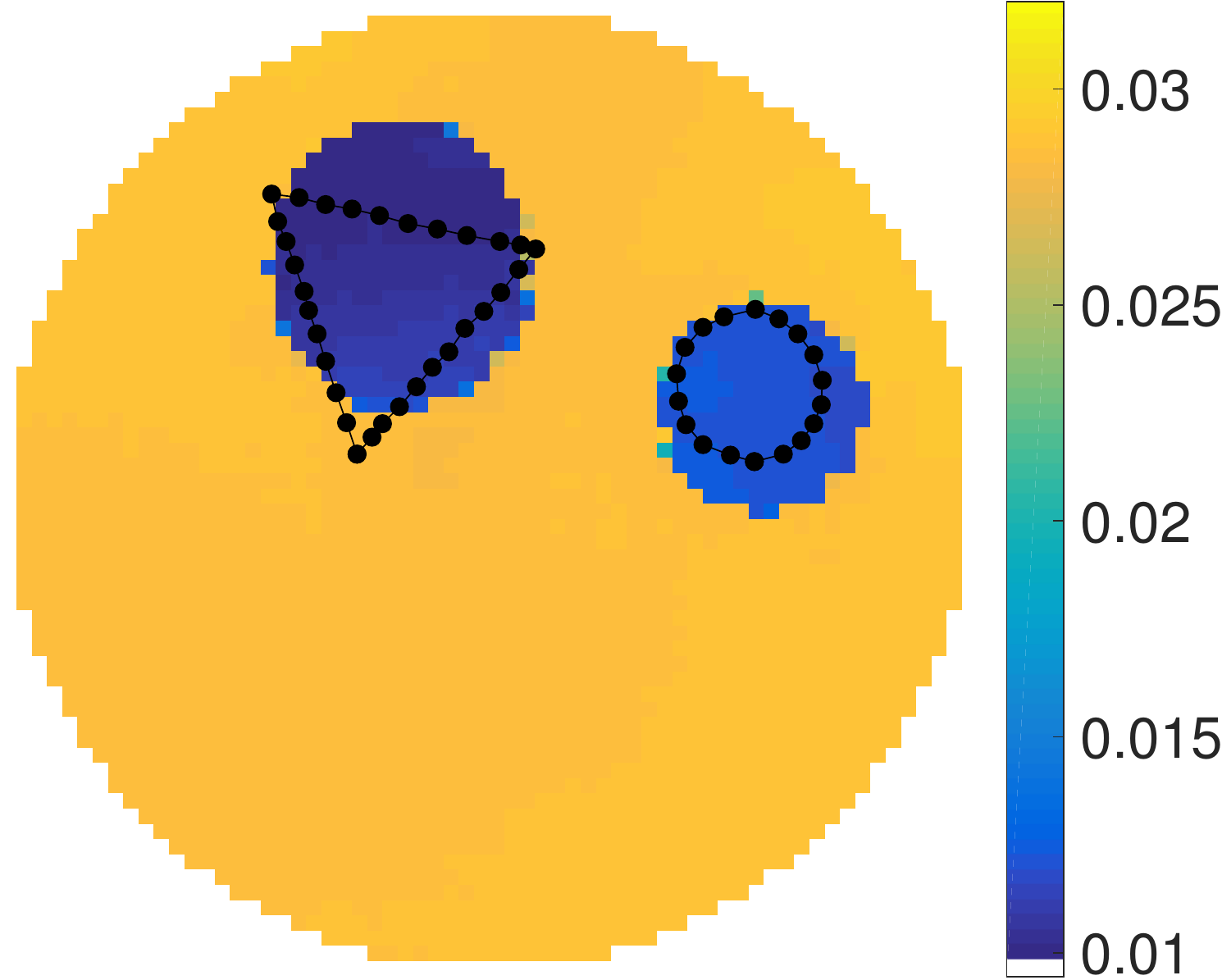}}

\put(10,85){\includegraphics[width=75pt]{fantom_3_4.jpg}}
\put(85,85){\includegraphics[width=75pt]{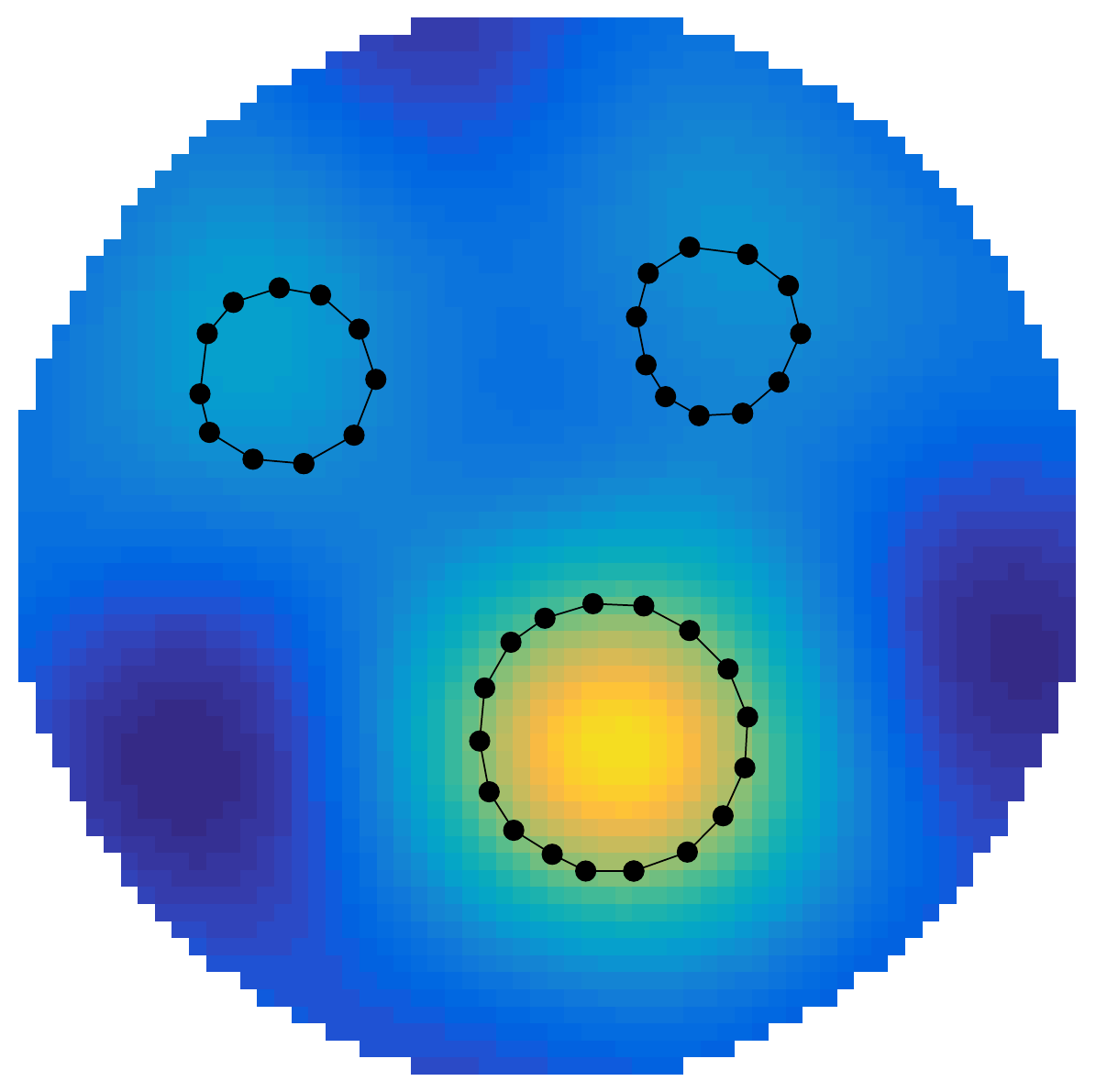}}
\put(160,85){\includegraphics[width=90pt]{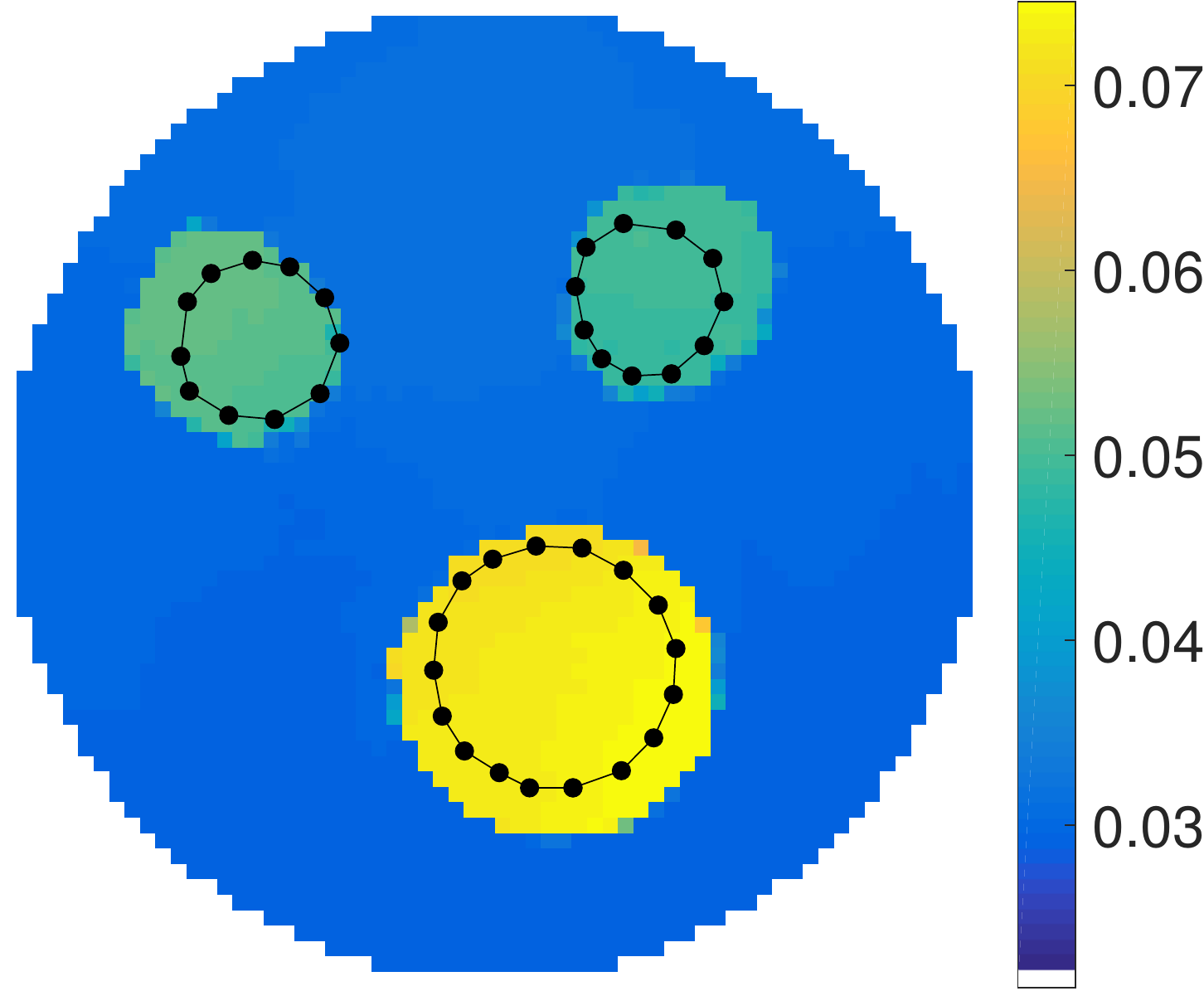}}

\put(10,0){\includegraphics[width=75pt]{fantom_4_4.jpg}}
\put(85,0){\includegraphics[width=75pt]{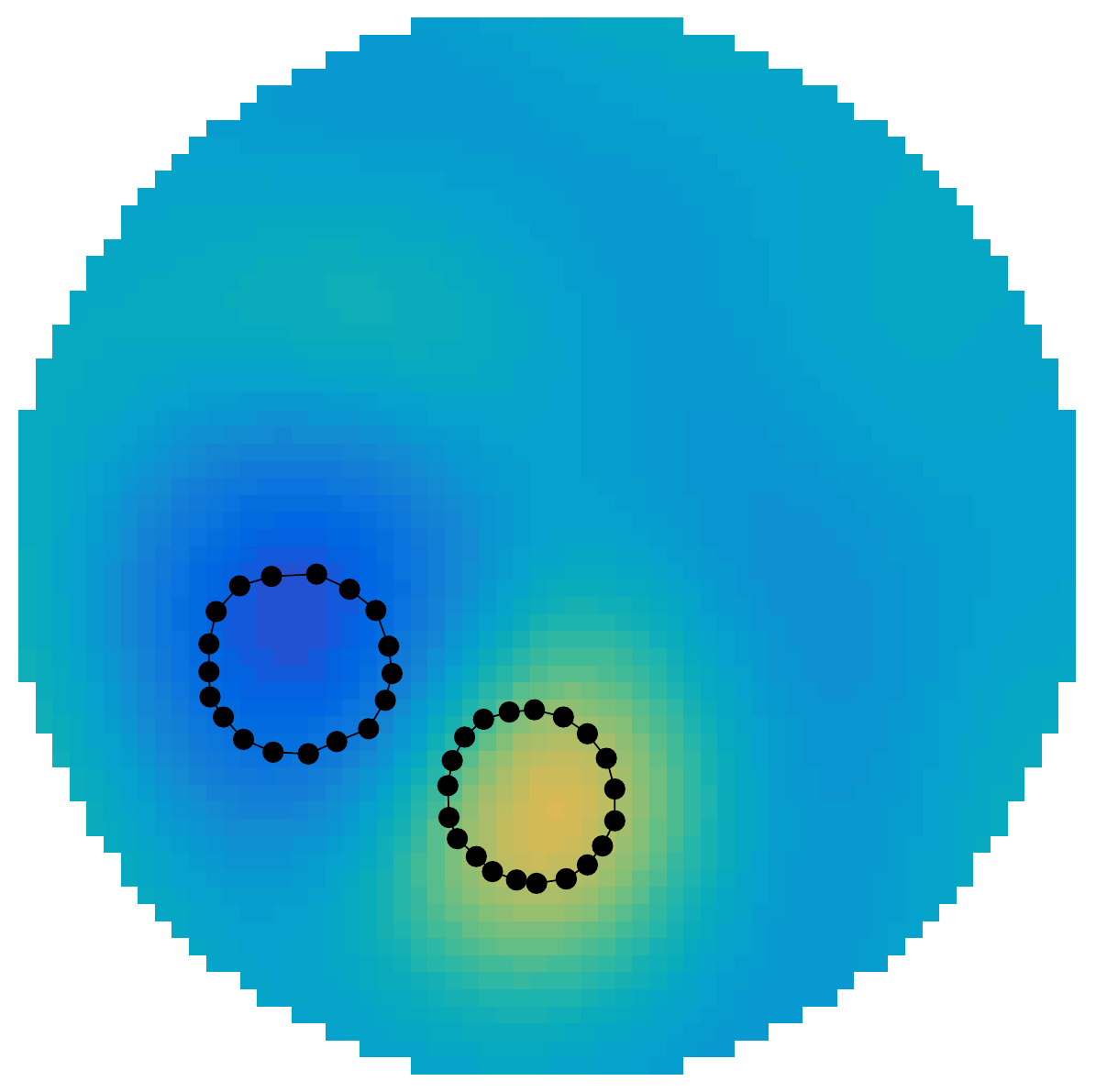}}
\put(160,0){\includegraphics[width=90pt]{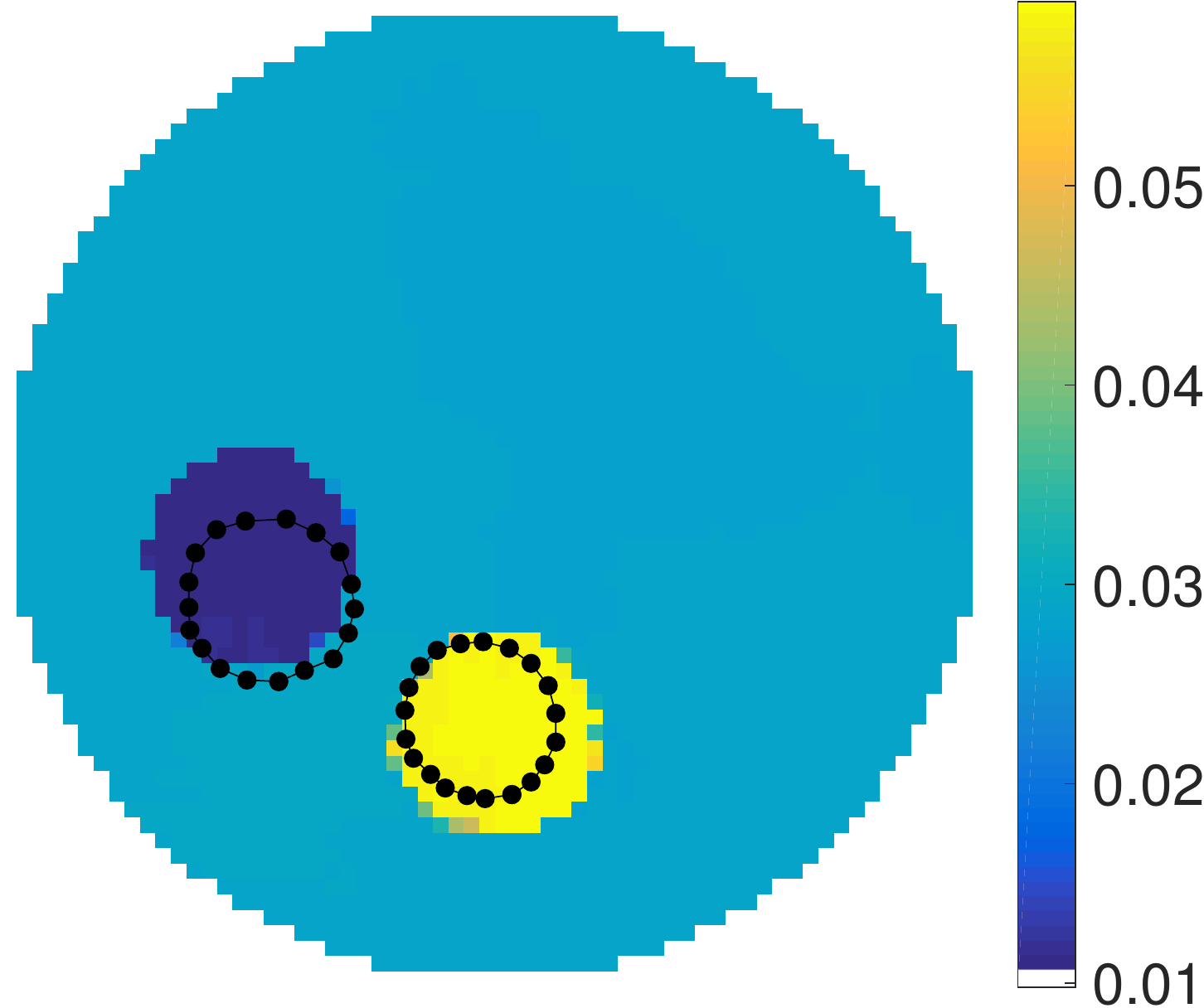}}

\put(25,345){\footnotesize \sc Experiment}

\put(108,355){\footnotesize \sc Low Pass}
\put(100,345){\footnotesize \sc D-bar Image}

\put(170,345){\footnotesize \sc  Deep D-bar Image}

\put(0,15){\rotatebox{90}{{\sc \footnotesize Phantom 4.4}}}
\put(0,100){\rotatebox{90}{{\sc \footnotesize Phantom 3.4}}}
\put(0,185){\rotatebox{90}{{\sc \footnotesize Phantom 2.2}}}
\put(0,275){\rotatebox{90}{{\sc \footnotesize Phantom 1.1}}}

\end{picture}
\caption{\label{fig:KIT4_Results} KIT4 Results for the various phantoms with conductive and/or resistive targets, as shown in the first column.  The initial D-bar image is compared to the Deep D-bar image.  The D-bar images, on the full square are used as the `input' images for the CNN.  Images are displayed here on the circular geometry of the tank, for presentation only.  Each row is plotted on its own scale.}
\end{figure}

Lastly, Figure~\ref{fig:KIT4_Results} shows results of the method on the four KIT4 scenarios shown in Figure~\ref{fig:KIT4_phantoms}.  The overlaid black dots depict the approximate `true' locations of the targets as extracted from their corresponding photographs.  No SSIMs were computed here since the objects are infinite conductors and resistors. \trevNew{For a comparison to an iterative method with a total variation prior \cite{Gonzalez2017}, performed on the same KIT4 data, we refer the reader to the documentation \cite{kit4data}.}

\section{Discussion}\label{sec:discussion}
The reconstructions shown in Figures~\ref{fig:ACT4_Results_simData}, \ref{fig:KIT4_Results_simData}, \ref{fig:ACT4_Results}, and \ref{fig:KIT4_Results} demonstrate that {\it Deep D-bar} provides superior reconstructions giving both visual and quantitative improvements.  In particular, the SSIMs (Figs. \ref{fig:act4_SSIMplot_simData}, \ref{fig:kit4_SSIMplot_simData}, and \ref{fig:act4_SSIMplot}\trev{)} show significant SSIM increases for the {\it Deep D-bar} vs. {\it Low-pass D-bar}.  Note that for the SSIM computation for ACT4 {\it Injury 2} (plastic tubes), the `truth' image was unrealistically set to zero in the lower portion of the right lung, even though the tubes do not entirely fill that region.

We remind the reader that no experimental (truth, reconstruction) pairs were used in training the network and no adaptation to the experimental system was necessary, apart from the number of electrodes in the system. The training was done purely with simulated data. In most applications, either a transfer training \cite{Hauptmann2018} or training with a golden standard from the same system must be performed.  This demonstrates the robustness of our approach.  Additionally, we expect further improvements in the ACT4 reconstructions if more complicated injuries are included in the training and remind the reader that the {\it Low-pass D-bar} and {\it Deep D-bar} reconstructions are shown on the same scale, which does mask the true dynamic range of the \trev{individual} images.

We review additional simplifications used in our process: 1) we used the {\it continuum electrode model} for the boundary conditions in the training data, 2) the FEM solver used to form $\mathbf{L}_1$ for the ACT4 and KIT4 experimental data examples was not finely tuned to  either EIT device (which is required for iterative minimization-based methods), and 3) the D-bar solver was not optimized for the respective ACT4/KIT4 data.  Rather it was used merely to provide the low-pass reconstructions used as inputs in the CNN.  These simplifications were used to demonstrate the robustness of the approach to both noise in the data and tolerance to modeling errors at multiple stages of the reconstruction process.  

\trev{Initial experiments performed with the original U-Net architecture, i.e. convolutional filters of size $3\times 3$, did not perform satisfactorily leading us to increase the filter size in this study to $5\times 5$. This of course leads to an increase in parameters from $8.6\cdot 10^6$ to $2.4\cdot 10^7$ resulting in longer training times.  No batch normalization was needed in our training processes and training times were only 4 hours per network, due to the rather small image size.}

The evaluation of the CNN is highly efficient on a GPU and \trev{took on} average $7.65$ms for a single sample, hence we expect Deep D-bar to be real-time capable. This can be done by combining the D-bar reconstruction, as outlined in \cite{Dodd2014}, with the application of the CNN in a unified framework to reduce overhead due to data transmission.


\subsection{Generalization}\label{sec:general}
\trev{An important aspect for medical imaging is the robustness and consistency of reconstructions. The successful transition to experimental data suggests that the proposed {\it Deep D-bar} method is robust enough for translational imaging.  Furthermore, Figures~\ref{fig:ACT4_Results_simData} and \ref{fig:ACT4_Results} (ACT4) illustrate that the network can handle reconstructions of phantoms that do not conform with the training data. However, while we were able to localize the inclusions in Figure~\ref{fig:KIT4_Results} (Phantom 2.2, KIT4), the sharp angular boundaries of the triangular target were not recovered when using only circular inclusion training data.  Our initial testing suggests that this can be improved upon by including significant training on triangular inclusions.  Challenges recovering triangular shapes have also been observed in \cite{Liu2018}. \trevNew{In terms of image quality, our Deep {D-bar} approach is comparable to results from (slower) iterative methods on similar data from the KIT4 system, see \cite{kit4data,Gonzalez2017,Liu2018}.  }
  Additionally, as the ACT4 injuries we simulated were elementary (only using a horizontal dividing line rather than the true diagonal cut and incomplete regional replacements), we expect that the reconstructions may improve further if more complex injuries were introduced.  For human targets, a larger database of training data could be employed and built from an anatomical atlas or collection of CT/MR scans both including and not including abnormalities/injuries.}

\begin{figure}[h!]
\centering
\begin{picture}(300,100)

\put(0,0){\includegraphics[width=75pt]{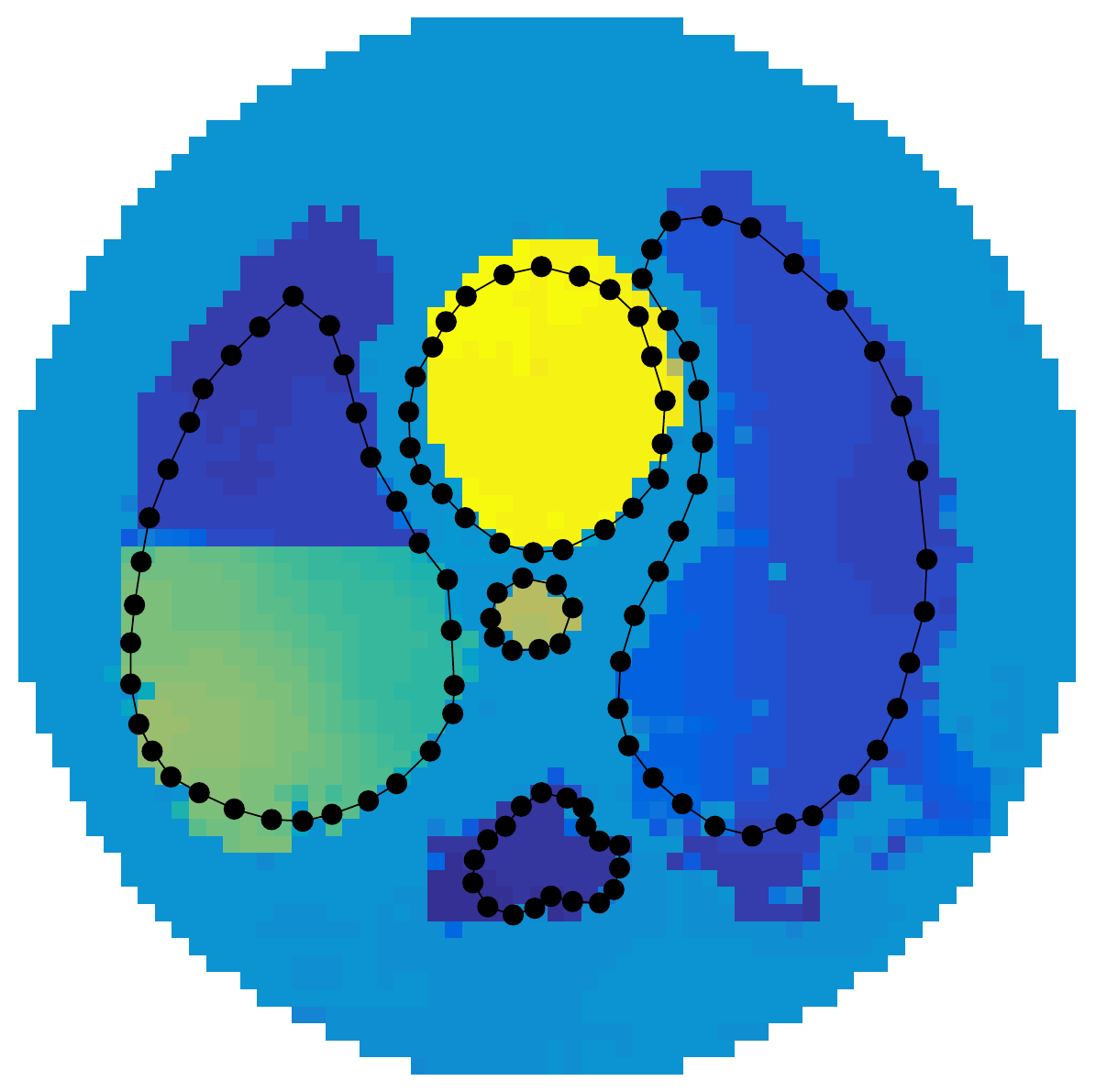}}
\put(80,0){\includegraphics[width=75pt]{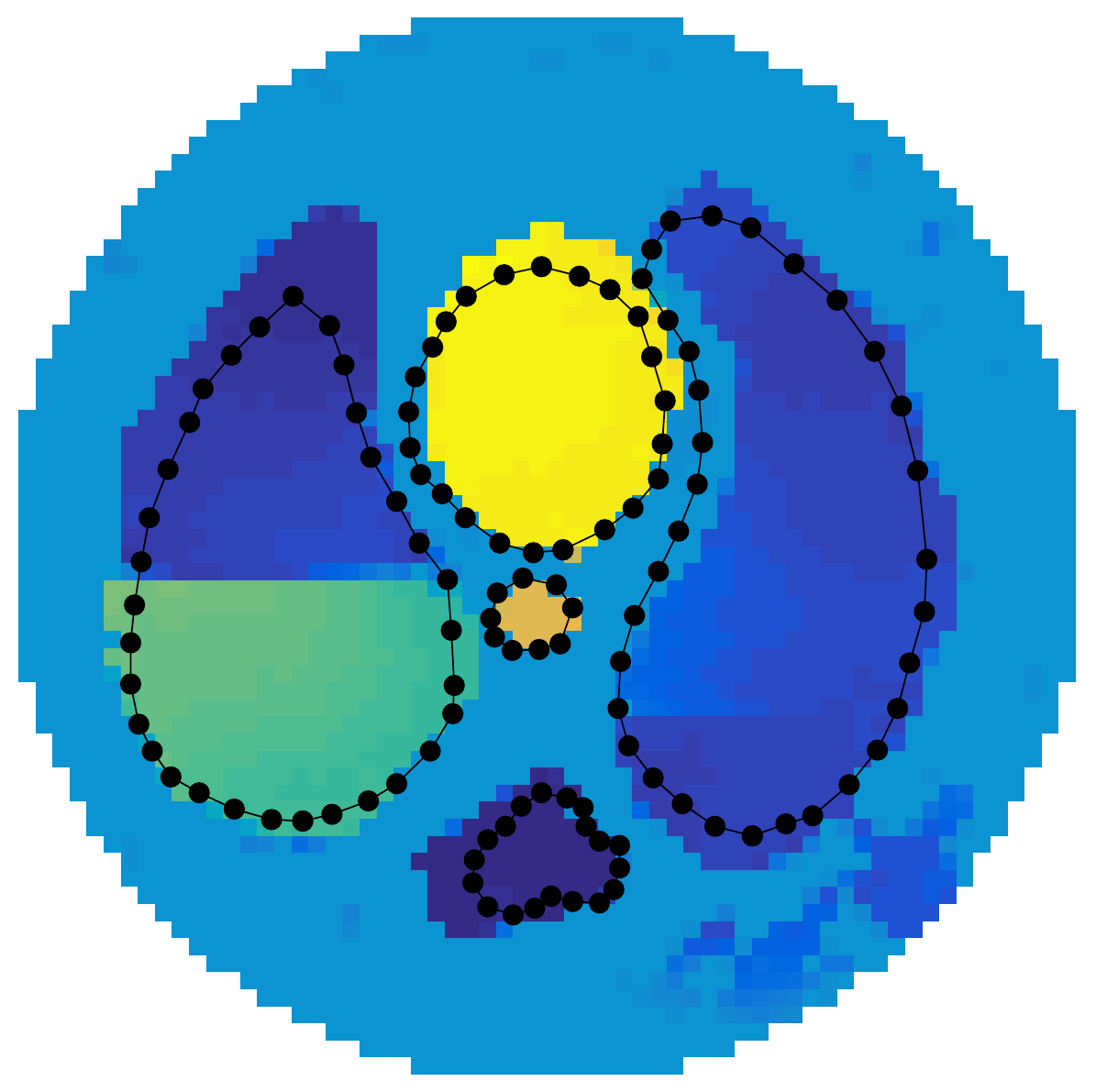}}
\put(160,0){\includegraphics[width=90pt]{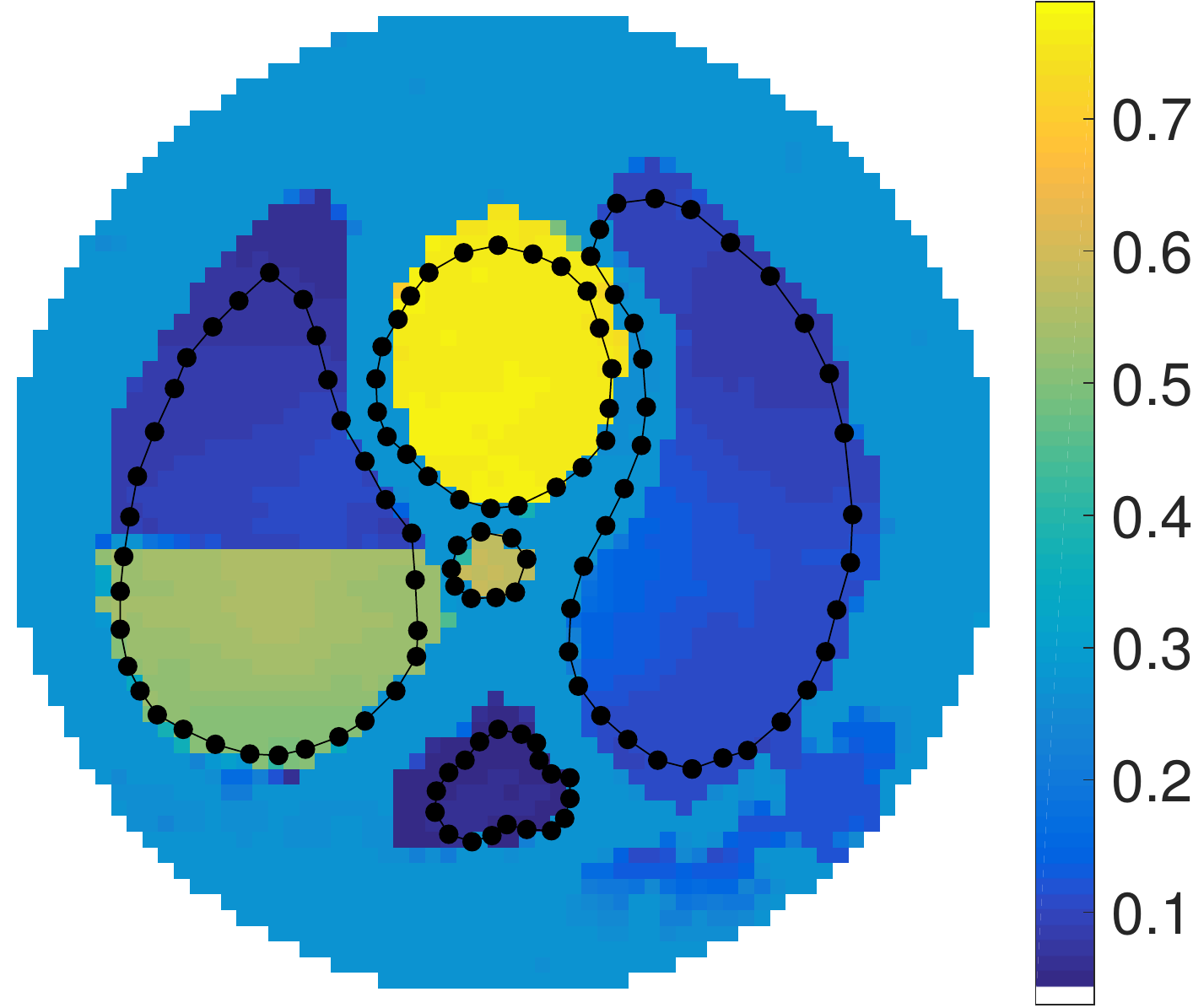}}

\put(10,90){\footnotesize \sc\underline{Old Network}}
\put(22,80){\footnotesize \sc$R=4.5$}

\put(130,90){\footnotesize \sc \underline{New Network}}
\put(105,80){\footnotesize \sc$R=4$}
\put(182,80){\footnotesize \sc$R=4.5$}
\end{picture}
\caption{\label{fig:varyingScatRad} \trev{Comparison of results for the ACT4 `Injury 1' (conductive agar in a lung) dataset from two different CNNs.  The `old' network denotes the network trained used a fixed cutoff radius for the scattering data ($R=4.5$) whereas the `new' network was trained with varying cut-off radii (randomized from the interval $[4,5]$). The result from the original, fixed radius $R=4.5$ network, is compared to results using $R=4$ and $R=4.5$ for the input image in the new network.  The SSIMs remained consistent: 0.6405, 0.6397, and 0.6459, from left to right.}}
\end{figure}
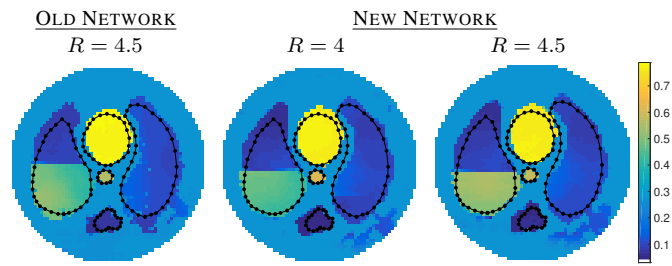

\trev{Crucial for the success of the post-processing network is consistency in the input reconstructions. In order to improve flexibility of the network, one can train the network on reconstructions from scattering data with varying cut-off radii. This allows the user to decide on the quality of the measured data at hand and adjust the cut-off radii as needed for the input reconstruction to the network. First tests have shown that this procedure indeed improves consistency and stability of the reconstructions as illustrated in Figure~\ref{fig:varyingScatRad}, where we have trained the ACT4 network with varying cut-off radii $R\in[4,5]$.  While the SSIMs remained consistent, the localization and recovered conductivity of the injury did improve with new variable radii network.}

\trev{While we chose, in this study, to match D-bar and CNNs due to the robustness of D-bar for absolute and time-difference imaging and the convolved nature of the D-bar reconstructions, alternative reconstruction methods for the input images could also be used.}

\section{Conclusion\trev{s}}\label{sec:conclusion}
The D-bar method for 2D EIT provides reliable reconstructions of the conductivity but suffers from a blurring due to a low-pass filtering of the scattering data.  Sharp improvements in \trev{absolute EIT} image quality can be achieved by coupling the D-bar reconstruction method with a convolutional neural network.  We demonstrated that a CNN can effectively learn the deblurring using only simulated data and still transition to experimental data without including any experimental data in the training itself.  As the training can be done offline ahead of time, and the D-bar method provides real-time conductivity reconstructions \cite{Dodd2014}, the post-processing step by the trained CNN adds minimal time to the overall image recovery process, due to the highly efficient evaluation on a GPU.  While this work is shown in 2D, we expect the approach to extend to 3D once the D-bar computational framework has been further developed.

\section*{Acknowledgments}
We gratefully acknowledge the support of NVIDIA Corporation with the donation of the Titan Xp GPU used for this research.
The authors would like to thank the Isaac Newton Institute for
Mathematical Sciences, Cambridge, for support and hospitality during
the programme `Variational methods and effective algorithms for imaging and vision' where work on this paper was
undertaken, EPSRC grant no EP/K032208/1. AH acknowledges support from EPSRC project EP/M020533/1, `Medical image computing for next-generation healthcare technology'. We additionally thank the EIT group at RPI \footnote{\url{https://www.ecse.rpi.edu/homepages/saulnier/eit/eit.html}} for their assistance and for providing the {\sc ACT4} tank data.

\bibliographystyle{unsrt}
\small
\bibliography{Inverse_problems_references_2017}

\begin{thebibliography}{10}

\bibitem{Cinnella2015}
Gilda Cinnella, Salvatore Grasso, Pasquale Raimondo, Davide D’Antini, Lucia
  Mirabella, Michela Rauseo, and Michele Dambrosio.
\newblock Physiological effects of the open lung approach in patients with
  early, mild, diffuse acute respiratory distress syndromean electrical
  impedance tomography study.
\newblock {\em The Journal of the American Society of Anesthesiologists},
  123(5):1113--1121, 2015.

\bibitem{Grant2011}
C.A. Grant, T.~Pham, J.~Hough, T.~Riedel, C.~Stocker, and A.~Schibler.
\newblock Measurement of ventilation and cardiac related impedance changes with
  electrical impedance tomography.
\newblock {\em Critical Care}, 15(1):R37, 2011.

\bibitem{Karagiannidis2015}
Christian Karagiannidis, Andreas~D. Waldmann, Carlos Ferrando~Ortol{\'a},
  Manuel Mu{\~n}oz~Martinez, Anxela Vidal, Arnoldo Santos, Peter~L. R{\'o}ka,
  Manuel Perez~M{\'a}rquez, Stephan~H. Bohm, and Fernando Suarez-Spimann.
\newblock Position-dependent distribution of ventilation measured with
  electrical impedance tomography.
\newblock {\em European Respiratory Journal}, 46(suppl 59), 2015.

\bibitem{Pesenti2016}
Antonio Pesenti, Guido Musch, Daniel Lichtenstein, Francesco Mojoli, Marcelo
  B.~P. Amato, Gilda Cinnella, Luciano Gattinoni, and Michael Quintel.
\newblock Imaging in acute respiratory distress syndrome.
\newblock {\em Intensive Care Medicine}, 42(5):686--698, 2016.

\bibitem{Reinius2015}
H.~Reinius, J.~B. Borges, F.~Fred\'{e}n, L.~Jideus, E.~D. L.~B. Camargo,
  M.~B.~P. Amato, G.~Hedenstierna, Larsson A., and F.~Lennmyr.
\newblock Real-time ventilation and perfusion distributions by electrical
  impedance tomography during one-lung ventilation with capnothorax.
\newblock {\em Acta Anaesthesiologica Scandinavica}, 59(3):354--368, 2015.

\bibitem{Schlibler2013}
A.~Schlibler, T.M.Y. Pham, A.A. Moray, and C.~Stocker.
\newblock Ventilation and cardiac related impedance changes in children
  undergoing corrective open heart surgery.
\newblock {\em Physiological Measurement}, 34:1319--1327, 2013.

\bibitem{Knudsen2009}
K.~Knudsen, M.~Lassas, J.L. Mueller, and S.~Siltanen.
\newblock Regularized {D}-bar method for the inverse conductivity problem.
\newblock {\em Inverse Problems and Imaging}, 3(4):599--624, 2009.

\bibitem{Zhou2015}
Zhou Zhou, Gustavo~Sato dos Santos, Thomas Dowrick, James Avery, Zhaolin Sun,
  Hui Xu, and David~S Holder.
\newblock Comparison of total variation algorithms for electrical impedance
  tomography.
\newblock {\em Physiological measurement}, 36(6):1193, 2015.

\bibitem{Gonzalez2017}
Gerardo Gonz{\'a}lez, Ville Kolehmainen, and Aku Sepp{\"a}nen.
\newblock Isotropic and anisotropic total variation regularization in
  electrical impedance tomography.
\newblock {\em Computers \& Mathematics with Applications}, 2017.

\bibitem{Darde2013}
J{\'e}r{\'e}mi Dard{\'e}, N~Hyv{\"o}nen, A~Sepp{\"a}nen, and Stratos Staboulis.
\newblock Simultaneous reconstruction of outer boundary shape and admittivity
  distribution in electrical impedance tomography.
\newblock {\em SIAM Journal on Imaging Sciences}, 6(1):176--198, 2013.

\bibitem{Nissinen2011}
A.~Nissinen, V.~Kolehmainen, and J.~P. Kaipio.
\newblock Compensation of modelling errors due to unknown domain boundary in
  electrical impedance tomography.
\newblock {\em IEEE Transaction on Medical Imaging}, 30(2):231--242, 2011.

\bibitem{Murphy2009}
E.~K. Murphy and J.~L. Mueller.
\newblock Effect of domain-shape modeling and measurement errors on the 2-d
  {D}-bar method for electrical impedance tomography.
\newblock {\em IEEE Transactions on Medical Imaging}, 28(10):1576--1584, 2009.

\bibitem{Kolehmainen2008b}
V.~Kolehmainen, M.~Lassas, and P.~Ola.
\newblock Electrical impedance tomography problem with inaccurately known
  boundary and contact impedances.
\newblock {\em Medical Imaging, IEEE Transactions on}, 27(10):1404 --1414, oct.
  2008.

\bibitem{Jin2017}
Kyong~Hwan Jin, Michael~T McCann, Emmanuel Froustey, and Michael Unser.
\newblock Deep convolutional neural network for inverse problems in imaging.
\newblock {\em IEEE Transactions on Image Processing}, 26(9):4509--4522, 2017.

\bibitem{Kang2017}
Eunhee Kang, Junhong Min, and Jong~Chul Ye.
\newblock A deep convolutional neural network using directional wavelets for
  low-dose x-ray ct reconstruction.
\newblock {\em Medical Physics}, 44(10), 2017.

\bibitem{Sandino2017}
Christopher~M Sandino, Neerav Dixit, Joseph~Y Cheng, and Shreyas~S Vasanawala.
\newblock Deep convolutional neural networks for accelerated dynamic magnetic
  resonance imaging.
\newblock In {\em NIPS 2017, Medical Imaging Meets NIPS Workshop.}, 2017.

\bibitem{Antholzer2017}
Stephan Antholzer, Markus Haltmeier, and Johannes Schwab.
\newblock Deep learning for photoacoustic tomography from sparse data.
\newblock {\em arXiv preprint arXiv:1704.04587}, 2017.

\bibitem{Hauptmann2018}
A.~Hauptmann, F.~Lucka, M.~Betcke, N.~Huynh, J.~Adler, B.~Cox, P.~Beard,
  S.~Ourselin, and S.~Arridge.
\newblock Model based learning for accelerated, limited-view 3d photoacoustic
  tomography.
\newblock {\em IEEE Transactions on Medical Imaging}, 2018.

\bibitem{Adler2017}
Jonas Adler and Ozan {\"O}ktem.
\newblock Solving ill-posed inverse problems using iterative deep neural
  networks.
\newblock {\em Inverse Problems}, 33(12):124007, 2017.

\bibitem{Delbary2012}
Fabrice Delbary, Per~Christian Hansen, and Kim Knudsen.
\newblock Electrical impedance tomography: 3d reconstructions using scattering
  transforms.
\newblock {\em Applicable Analysis}, 91(4):737--755, 2012.

\bibitem{Nachman1996}
A.~I. Nachman.
\newblock Global uniqueness for a two-dimensional inverse boundary value
  problem.
\newblock {\em Annals of Mathematics}, 143:71--96, 1996.

\bibitem{Beals1985}
Richard Beals and Ronald~R. Coifman.
\newblock Multidimensional inverse scatterings and nonlinear partial
  differential equations.
\newblock In {\em Pseudodifferential operators and applications (Notre Dame,
  Ind., 1984)}, pages 45--70. Amer. Math. Soc., Providence, RI, 1985.

\bibitem{Alessandrini1988}
G.~Alessandrini.
\newblock Stable determination of conductivity by boundary measurements.
\newblock {\em Applicable Analysis}, 27:153--172, 1988.

\bibitem{Isaacson2004}
D.~Isaacson, J.~L. Mueller, J.~C. Newell, and S.~Siltanen.
\newblock Reconstructions of chest phantoms by the {D}-bar method for
  electrical impedance tomography.
\newblock {\em IEEE Transactions on Medical Imaging}, 23:821--828, 2004.

\bibitem{Dodd2014}
Melody Dodd and Jennifer~L Mueller.
\newblock A real-time {D}-bar algorithm for 2-{D} electrical impedance
  tomography data.
\newblock {\em Inverse problems and imaging}, 8(4):1013--1031, 2014.

\bibitem{Hamilton2017_PhysMeas2}
S.~J. Hamilton, J.~L. Mueller, and T.~R. Santos.
\newblock Robust computation of 2d {EIT} absolute images with d-bar methods.
\newblock (under review, arXiv preprint: \url{1712.00379}).

\bibitem{Henkin2010}
Gennadi Henkin and Matteo Santacesaria.
\newblock On an inverse problem for anisotropic conductivity in the plane.
\newblock {\em Inverse Problems}, 26(9):095011, 2010.

\bibitem{Hamilton2014a}
S.~J. Hamilton, M.~Lassas, and S.~Siltanen.
\newblock A {D}irect {R}econstruction {M}ethod for {A}nisotropic {E}lectrical
  {I}mpedance {T}omography.
\newblock {\em Inverse Problems}, 30:(075007), 2014.

\bibitem{Ronneberger2015}
Olaf Ronneberger, Philipp Fischer, and Thomas Brox.
\newblock U-net: Convolutional networks for biomedical image segmentation.
\newblock In {\em International Conference on Medical Image Computing and
  Computer-Assisted Intervention}, pages 234--241. Springer, 2015.

\bibitem{Martin2017}
S{\'e}bastien Martin and Charles~TM Choi.
\newblock A post-processing method for three-dimensional electrical impedance
  tomography.
\newblock {\em Scientific Reports}, 7, 2017.

\bibitem{Wiatowski2017}
Thomas Wiatowski and Helmut B{\"o}lcskei.
\newblock A mathematical theory of deep convolutional neural networks for
  feature extraction.
\newblock {\em IEEE Transactions on Information Theory}, 2017.

\bibitem{Liu2005}
N.~Liu, G.J. Saulnier, J.C. Newell, D.~Isaacson, and T-J. Kao.
\newblock Act4: a high-precision, multi-frequency electrical impedance
  tomograph.
\newblock Presented at 6th Conference on Biomedical Applications of Electrical
  Impedance Tomography, June 2005.
\newblock London, U.K.

\bibitem{Saulnier2007}
G.~Saulnier, N.~Liu, C.~P. Tamma, H.~Xia, T.-J. Kao, J.~Newell, and
  D.~Isaacson.
\newblock An electrical impedance spectroscopy system for breast cancer
  detection.
\newblock {\em Proc. 29th Int. Conf. IEEE Eng. Med. Biol. Soc}, (1):4154--4157,
  Aug. 2007.

\bibitem{kit4data}
Andreas Hauptmann, Ville Kolehmainen, Nguyet~Minh Mach, Tuomo Savolainen, Aku
  Sepp\"anen, and Samuli Siltanen.
\newblock Open 2d {E}lectrical {I}mpedance {T}omography data archive.
\newblock {\em arXiv:1704.01178}, 2017.

\bibitem{Hauptmann2017a}
Andreas Hauptmann.
\newblock Approximation of full-boundary data from partial-boundary electrode
  measurements.
\newblock {\em Inverse Problems (accepted), arXiv preprint arXiv:1703.05550},
  2017.

\bibitem{Somersalo1992}
Erkki Somersalo, Margaret Cheney, and David Isaacson.
\newblock Existence and uniqueness for electrode models for electric current
  computed tomography.
\newblock {\em SIAM Journal on Applied Mathematics}, 52(4):1023--1040, 1992.

\bibitem{Hyvoenen2009}
N.~Hyv{\"o}nen.
\newblock Approximating idealized boundary data of electric impedance
  tomography by electrode measurements.
\newblock {\em Mathematical Models and Methods in Applied Sciences},
  19(07):1185--1202, 2009.

\bibitem{Mueller2012}
J.~Mueller and S.~Siltanen.
\newblock {\em Linear and Nonlinear Inverse Problems with Practical
  Applications}, volume~10 of {\em Computational Science and Engineering}.
\newblock SIAM, 2012.

\bibitem{DeAngelo2010}
M.~DeAngelo and J.~L. Mueller.
\newblock 2d {D}-bar reconstructions of human chest and tank data using an
  improved approximation to the scattering transform.
\newblock {\em Physiological Measurement}, 31:221--232, 2010.

\bibitem{Vainikko2000}
Gennadi Vainikko.
\newblock Fast solvers of the {L}ippmann-{S}chwinger equation.
\newblock In {\em Direct and inverse problems of mathematical physics
  ({N}ewark, {DE}, 1997)}, volume~5 of {\em Int. Soc. Anal. Appl. Comput.},
  pages 423--440. Kluwer Acad. Publ., Dordrecht, 2000.

\bibitem{Liu2018}
D.~Liu, A.~K. Khambampati, and J.~Du.
\newblock A parametric level set method for electrical impedance tomography.
\newblock {\em IEEE Transactions on Medical Imaging}, 37(2):451--460, Feb 2018.

\end{thebibliography}
\end{document}